\documentclass[a4,12pt]{article}

\usepackage[british]{babel} %%% 'french', 'german', 'spanish', 'danish', etc.
\usepackage{ucs}                % Unicode support in LaTeX
\usepackage[utf8x]{inputenc}    % We use UTF-8 encoded files only 
\usepackage[T1]{fontenc}     % for the right font encoding
\usepackage{lmodern}         % fonts for good pdf files
\usepackage{float}

\usepackage[linesnumbered,lined,boxed,commentsnumbered]{algorithm2e}

\usepackage[singlespacing]{setspace}
\usepackage[
twoside,
top=2.5cm,
bottom=2.5cm,
inner=2.0cm,
outer=2.0cm,
includehead,
heightrounded
]{geometry}

\usepackage[sfdefault=cmbr,OMLmathsans]{isomath}
\usepackage{latexsym}
\usepackage{epsfig,xspace}

\usepackage{amscd}
\usepackage{upgreek}  % upright greek letters
\usepackage{tensor}   % for tensor index notation
\usepackage{refdef}
\usepackage{float}    % for [H] position of figures. 

\usepackage{fancyhdr}
\usepackage[export]{adjustbox}
\usepackage{amsthm}
\usepackage{amsfonts}
\usepackage{longtable}
\usepackage{microtype}
\usepackage[parfill]{parskip}
\usepackage{authblk}
\usepackage[tiny]{titlesec}
\usepackage{hyperref}
\hypersetup{colorlinks=true,linkcolor=blue,filecolor=magenta}
\usepackage{graphicx}
\usepackage{sidecap}
\usepackage[thinc]{esdiff}   % math symbols for derivative
\usepackage{tabularx}
\usepackage{color, colortbl}
\usepackage[usernames,dvipsnames]{xcolor}
\usepackage{amsmath}
\usepackage{amssymb}
\usepackage{epstopdf}
\usepackage{lipsum}
\usepackage{caption}
\usepackage{subcaption}
\usepackage{ifontdef} % define some fonts

\newcommand{\correction}[1]{\noindent \color{Black} #1\normalcolor}

\newcommand{\ignore}[1]{}

\newcommand{\vek}[1]{\mathchoice{\displaystyle\boldsymbol{#1}}
{\textstyle\boldsymbol{#1}}{\scriptstyle\boldsymbol{#1}}
{\scriptscriptstyle\boldsymbol{#1}}}

\newcommand{\ops}[1]{\mathchoice{\displaystyle\mathsf{#1}}
{\textstyle\mathsf{#1}}{\scriptstyle\mathsf{#1}}
{\scriptscriptstyle\mathsf{#1}}}

\newcommand{\tns}[1]{\mathchoice{\displaystyle\mathsans{#1}}
{\textstyle\mathsans{#1}}{\scriptstyle\mathsans{#1}}
{\scriptscriptstyle\mathsans{#1}}}

\newcommand{\di}{\mathrm{d}}

\def\Aop{\operatornamewithlimits{%
		\mathchoice{\vcenter{\hbox{\huge A}}}
		{\vcenter{\hbox{\Large A}}}
		{\mathrm{A}}
		{\mathrm{A}}}}

\title{Crack propagation simulation without crack tracking algorithm: embedded discontinuity formulation with incompatible modes}
\author[1,*]{A.~Stanic}
\author[2]{B.~Brank} 
\author[3]{A.~Ibrahimbegovic} 
\author[4]{H.~G.~Matthies} 

\affil[1]{\footnotesize{\emph{University of Twente, Faculty of Engineering Technology, The Netherlands}}} 
\affil[2]{\footnotesize{\emph{University of Ljubljana, Faculty of Civil and Geodetic Engineering, Slovenia}}}
\affil[3]{\footnotesize{\emph{Universit\'{e} de Technologie de Compi\`{e}gne - Alliance Sorbonne Universit\'{e}s, Laboratoire de M\'{e}canique, France}}}
\affil[4]{\footnotesize{\emph{Technische Universit\"at Braunschweig, Institute of Scientific Computing, Germany}}} 
\affil[*]{corresponding author, e-mail: \ttt{a.stanic@utwente.nl}}

\begin{document}
	
\maketitle
	
	\begin{abstract}		
		We show that for the simulation of crack propagation in quasi-brittle, two-dimensional solids, very good results can be obtained with an embedded strong discontinuity quadrilateral finite element that has incompatible modes. Even more importantly, we demonstrate that these results can be obtained without using a crack tracking algorithm. Therefore, the simulation of crack patterns with several cracks, including branching, becomes possible. The avoidance of a tracking algorithm is mainly enabled by the application of a novel, local (Gauss-point based) criterion for crack nucleation, which determines the time of embedding the localisation line as well as its position and orientation. We treat the crack evolution in terms of a thermodynamical framework, with softening variables describing internal dissipative mechanisms of material degradation. As presented by numerical examples, many elements in the mesh may develop a crack, but only some of them actually open and/or slide, dissipate fracture energy, and eventually form the crack pattern. The novel approach has been implemented for statics and dynamics, and the results of computed difficult examples (including Kalthoff's test) illustrate its very satisfying performance. It effectively overcomes unfavourable restrictions of the standard embedded strong discontinuity formulations, namely the simulation of the propagation of a single crack only. Moreover, it is computationally fast and straightforward to implement. Our numerical solutions match the results of experimental tests and previously reported numerical results in terms of crack pattern, dissipated fracture energy, and load-displacement curve.
 	\end{abstract}


	\section{Introduction}
	\label{sec:intro}
	The strong discontinuity approach (SDA) is a well established method (\cite{Dvorkin Cuitino Gioia 1990-art, Dvorkin Assanelli 1991-art, Klisinski et al 1991-art, Mosler 2004-art, Oliver et al 2006-art, Manzoli et al 2006-art, Linder and Armero 2007-art, Brancherie Ibrahimbegovic 2009-art, Linder and Armero 2009-art, Dujc et al 2010-art, Oliver et al 2012-art, Saksala et al 2015-art, Saksala et al 2016-art, Lloberas-Valls et al 2016-art, Nikolic et al 2018-art,Stanic Brank Brancherie 2020-art, Dujc Brank 2020-art}) for modelling crack propagation in quasi-brittle materials (such as concrete, high strength steel, and rock). In this work, we focus on the development of a corresponding quadrilateral finite element with a statically and kinematically optimal non-symmetric formulation \cite{Jirasek 2000-art}.

	When the SDA is applied for a two-dimensional finite element (\cite{Dvorkin Cuitino Gioia 1990-art, Dvorkin Assanelli 1991-art, Klisinski et al 1991-art, Mosler 2004-art, Oliver et al 2006-art, Linder and Armero 2007-art, Brancherie Ibrahimbegovic 2009-art, Dujc et al 2010-art, Oliver et al 2012-art, Lloberas-Valls et al 2016-art, Stanic Brank Brancherie 2020-art, Dujc Brank 2020-art}), the discontinuity is embedded as a straight line that crosses the element and encompasses the localised material failure. According to \cite{Willam 2000-art}, material failure is a sequence of events that starts with the \textcolor{black}{loss of material stability} in a point and leads to the degradation of a continuum into a discontinuum. In this respect, one distinguishes between two groups of criteria: failure criteria and fracture criteria \cite{Willam 2000-art}. The failure criterion signals the initiation of the material failure process. For quasi-brittle materials, a representative example is Rankine’s condition stating that the material failure starts when the maximum principal stress reaches the material tensile strength. In the SDA, it is often used in connection with the criterion for crack nucleation, which determines when and how the localisation line is embedded into the finite element. The stress intensity factors and the fracture energy criterion belong to the group of fracture criteria that identify the existence of a fully open crack or notched defect in a sense of fracture mechanics, and can be used for monitoring the crack propagation \cite{Willam 2000-art}.

	The efficiency of the SDA mostly depends on the correct prediction of the discontinuity path \cite{Oliver et al 2006-art}. For some problems, the possible crack geometry is known, which allows for prescribing the crack path prior to the analysis. In most cases, however, the crack path has to be found by computing the crack propagation according to the current stress state in the solid. In this sense, more recent publications on quadrilaterals with embedded strong discontinuity advocate the application of crack tracking algorithm as an indispensable tool for computing the crack propagation, see e.g.\ \cite{Oliver et al 2002-art, Manzoli et al 2006-art, Linder and Armero 2007-art, Linder and Armero 2009-art, Dujc et al 2010-art, Dujc CC 2010-art, Dujc et al 2013-art, Dujc Brank 2020-art, Stanic-Brank-Korelc 2016-art, Stanic Brank Brancherie 2020-art}. For a particular crack path, the crack tracking algorithm determines the front element, which is the next element to be added to the string of cracked elements. When the front element fulfills the crack nucleation criterion, the algorithm inserts a discontinuity line by imposing the geometric continuity of element-wise cracks across the finite element mesh. Both the specific crack tracking algorithm and the specific crack nucleation criterion play an important role in the efficiency of the SDA.

	The first quadrilaterals with embedded strong discontinuity, presented in \cite{Dvorkin Cuitino Gioia 1990-art, Dvorkin Assanelli 1991-art, Klisinski et al 1991-art}, were used for computing simple examples without relying on a crack tracking algorithm. In \cite{Klisinski et al 1991-art}, the shear test was simulated, and the crack nucleation criterion was fulfilled when the stress state in the element reached Rankine’s condition in all four integration points.	The same crack nucleation criterion was used in \cite{Dvorkin Assanelli 1991-art} for analysing a deep concrete beam under four-point bending with mixed plane stress quadrilaterals, which had a crack line (opening in mode I fashion) crossing the centre of the element. As already mentioned, more recent SDA formulations with embedded discontinuity quadrilaterals rely strongly on crack tracking algorithms. This was used in \cite{Manzoli et al 2006-art}, with the discontinuity embedded when Rankine's condition was satisfied at the centroid of the quadrilateral element --- with its direction perpendicular to the corresponding maximum principal stress direction. In \cite{Linder and Armero 2007-art}, the applied crack tracking algorithm allows for crack propagation through as many elements as have fulfilled the crack nucleation criterion. There, the latter is defined for Rankine's failure condition to be met in all four integration points, with the discontinuity being perpendicular to the average of the principal directions in the Gauss-points. In \cite{Dujc et al 2010-art, Stanic-Brank-Korelc 2016-art, Stanic Brank Brancherie 2020-art}, the crack tracking algorithm was essential for establishing the convergence of the formulations. The maximum principal stress and the corresponding principal direction were computed from the average stress state over all four integration points, and the crack line was embedded when the maximum principal stress exceeded the material tensile strength. In \cite{Linder and Armero 2009-art}, a special criterion based on the crack tip velocity was used to detect branching of the crack path.

	From the above mentioned references it can be concluded that the current state-of-the-art for the embedded discontinuity quadrilateral formulations relies on various crack nucleation criteria and various crack tracking algorithms. However, it is shown below that it is possible to perform reliable simulations of complex crack propagation problems by relying only on an appropriate crack nucleation criterion, without using any additional external crack tracking algorithm. Note that similar computations with the SDA have already been reported with limited success for simple triangular embedded strong discontinuity finite elements in \cite{Brancherie Ibrahimbegovic 2009-art, Saksala et al 2015-art, Do et al 2017-art, Nikolic et al 2018-art}, and for quadrilateral finite elements in \cite{Ibrahimbegovic Brancherie 2003-art}, though in the latter case only for the shear separation mode.

	Here the crack opening and the connected evolution of \emph{fractured surface} is interpreted not so much as a geometrical condition, but as a phenomenological thermodynamic variable belonging to a dissipative mechanism governing the fracture energy dissipation through a dissipation pseudo-potential \cite{Ibrahimbegovic-book} in the sense of standard generalised materials. This interpretation of the crack evolution is the reason why it is not necessary to apply an algorithm for tracking the crack path and thereby enforcing continuous cracks. The validity of this approach is demonstrated with the numerical evidence provided below in the section with numerical examples. As is usual in such models, the softening variables are located at the integration points on the embedded discontinuity line of the element, and have no continuity requirements across element boundaries. For this reason, it is argued that it is not necessary for the discontinuities to be geometrically connected into a continuous line on the mesh level. Rather, the embedded discontinuity is viewed as an internal mechanism for the dissipation of fracture energy, which is used as a measure of material degradation that may also describe complete material failure.

	In this work, we apply the embedded SDA with constant crack opening and sliding on the well-known Q6 quadrilateral with incompatible modes \cite{Wilson Ibrahimbegovic 1991-art, Ibrahimbegovic-book}. \correction{Our goal is to compute crack propagation with fairly coarse meshes. For this reason, we prefer Q6 to the standard Q4 quadrilateral. It is well known that Q6 can alleviate the excessively stiff response of Q4 for coarse meshes. Our numerical examples show that the incompatible modes have a favorable influence on the fracture analysis based on embedded SDA. The differences between the Q6 and Q4 embedded discontinuity formulations stand out in the vicinity of propagating crack. We show this in two instances in our examples, namely in four point bending test in Section \ref{sec:FourPB} and in Kalthoff's test in Section \ref{sec:Kalthoff}.} Quasi-brittle materials with infinitesimal linear (plane stress or plane strain) elasticity as the material model for the bulk of the element  are considered. For the description of cohesion and energy dissipation in the element's internal interface, a softening cohesive model that takes into account mixed mode crack opening is used to govern the fracture energy release for the crack evolution \cite{Brancherie Ibrahimbegovic 2009-art}. In particular, we choose two uncoupled, non-associative, damage models with Rankine's material failure condition in order to deal with crack opening and sliding related to fracture modes I and II, respectively \correction{(for coupled model we refer to \cite{Kucerova et al 2009-art, Do Ibrahimbegovic 2018-art})}. To display purely the capabilities of this approach, no plasticity in the sense of ductile behaviour is considered, either in the bulk or in the crack. Thus, at complete unloading of the solid no residual displacements occur. It is worth mentioning that the crack acts as a rigid interface before the onset of discrete damage and does not allow inter-penetration even if the crack closes.

	The main novelty of the present work is the demonstration that the derived embedded discontinuity (ED) Q6 quadrilateral with incompatible modes can be effectively applied for the simulation of crack propagation in quasi-brittle solids without relying on crack tracking algorithms to trace the propagation of cracks. This allows a highly efficient formulation and, when compared to other formulations, very rapid computations. To the best of our knowledge, this is for the first time that an ED formulation for a quadrilateral that does not use any crack tracking algorithm has proved successful for different crack separation modes. We see the main reason for the success of the present formulation in the application of the novel crack nucleation criterion. The crack nucleation criterion is enforced only locally, on the Gauss-point level, whilst in the other work cited above this is done on the element level, or with averaging the stresses on an even wider area. It is worth mentioning that in this work we use the term ``crack'' in a loose sense, because ``crack'' is called any embedded discontinuity, with or without cohesion. This does not strictly comply with Griffith's definition of a crack \cite{Griffith 1921-art}, namely that ``a crack is assumed to be a traction-free surface''.

	The rest of the paper is organised so that Section~\ref{sec:model_description} briefly describes the ED formulation with incompatible modes, where Subsection~\ref{sec:ED_Criterion} explains the nucleation criterion and the computation of the crack orientation. The solution procedure for statics and dynamics is described in Section~\ref{sec:EqSys}, while Section~\ref{sec:NumExamples} contains several numerical examples. The conclusions are given in the final Section~\ref{sec:Conclusions}.

	\section{Model description}
	\label{sec:model_description}
	This section briefly describes the formulation for the embedded strong discontinuity Q6 quadrilateral with incompatible modes that can be used for static and dynamic simulations of crack propagation in quasi-brittle two-dimensional solids. Finally, the nucleation criterion and computation of the crack direction and position are explained in detail.

	\subsection{Separation modes} 
\begin{figure}[!htb]
		\centering
		\includegraphics[height=0.3\textwidth]{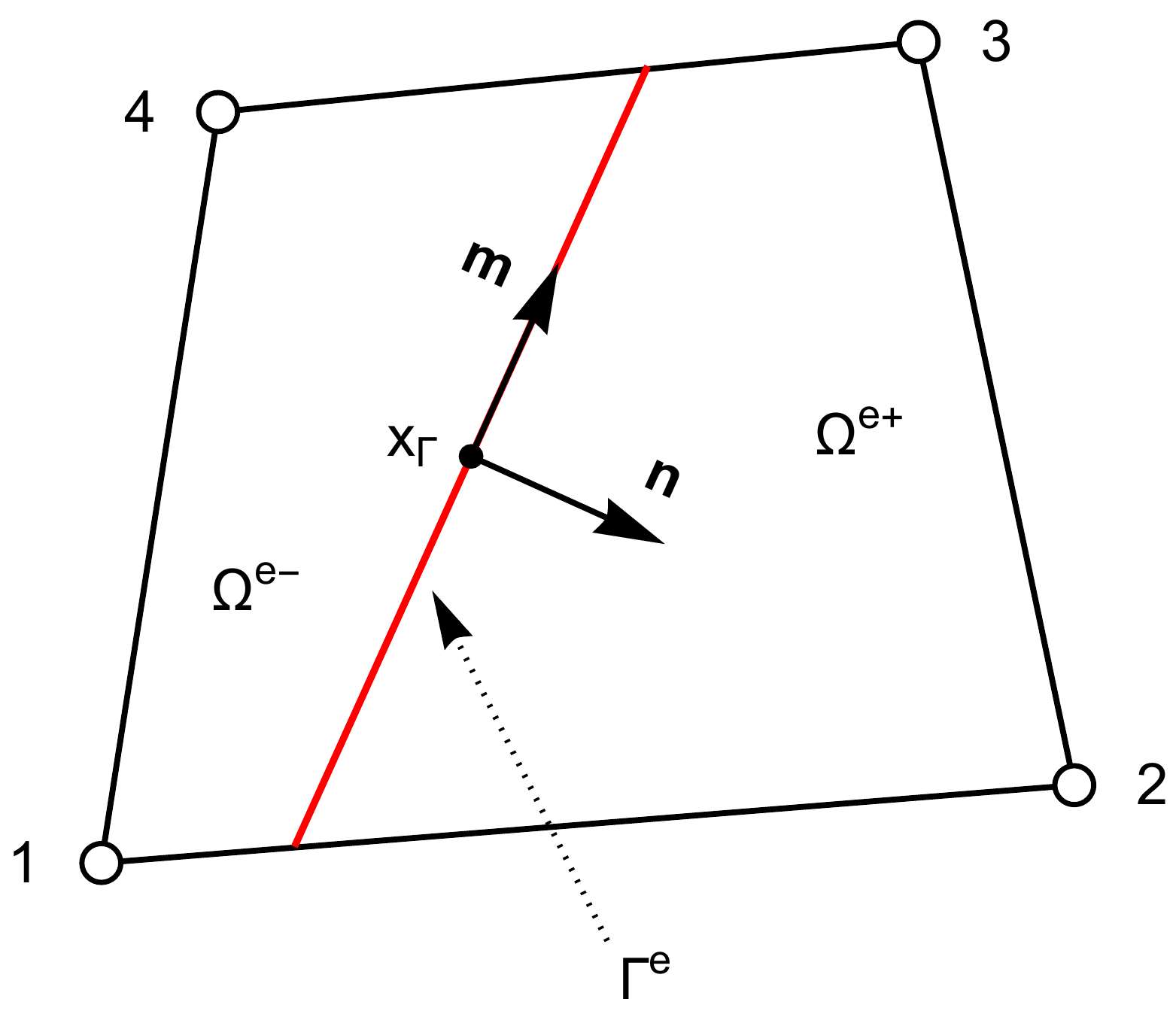}
		\caption{Quadrilateral element with inner interface.}
		\label{fig:Q4}
	\end{figure}

\begin{figure}[!htb]
		\centering
		\includegraphics[height=0.35\textwidth]{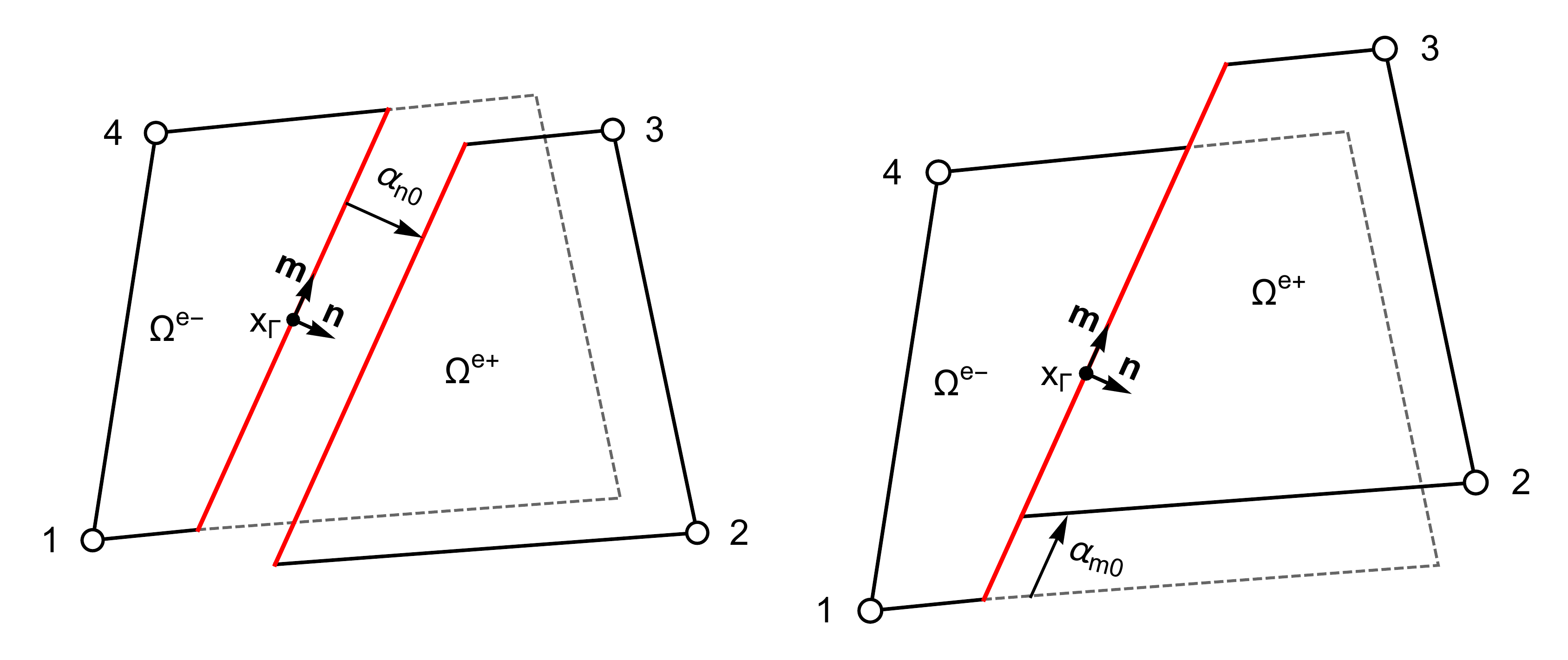}
		\caption{Separation modes ($\alpha_{1}=\alpha_{n0}$, $\alpha_{2}=\alpha_{m0}$).}
		\label{fig:Q4_SepModes}
	\end{figure}

	Let us consider a quadrilateral finite element as shown in Fig.~\ref{fig:Q4}, divided by an inner line (interface) ${\Gamma}^e$ of length $l_{\Gamma}$ into ${\Omega}^{e+}$ and ${\Omega}^{e-}$ (with ${\Omega}^e = {\Omega}^{e+} \cup {\Gamma}^e \cup {\Omega}^{e-}$ and ${\Omega}^{e+} \cap {\Omega}^{e-} = \emptyset $). Its geometry is defined by the bi-linear mapping $\boldsymbol{\xi} \to \boldsymbol{x}$ with $\boldsymbol{x} (\boldsymbol{\xi}) = \sum^4_{a=1}{N_a (\boldsymbol{\xi})}{\boldsymbol{x}}_a$, ${\boldsymbol{x}}_a={\left[x_a, y_a\right]}^T$, where $\boldsymbol{x}={\left[x,y\right]}^T \in {\Omega}^e$, $\boldsymbol{\xi} = (\xi, \eta) \in \left[{-1,1}\right] \times \left[{-1,1}\right]$, ${\boldsymbol{x}}_a$ are the coordinates of node $a$, and $N_a (\boldsymbol{\xi}) = \frac{1}{4} (1+{\xi }_a\xi) (1+{\eta}_a\eta)$ are Lagrange interpolation functions. In order to model the opening and sliding between ${\Omega}^{e-}$ and ${\Omega}^{e+}$, parameters ${\alpha}_i, i=1,2,$ are introduced, which are amplitudes of two separation modes presented in Fig. \ref{fig:Q4_SepModes}. The first separation mode describes opening, and the second one describes sliding. % between ${\Omega}^{e-}$ and ${\Omega}^{e+}$.
	
	\subsection{Coupling incompatible modes with separation modes}
	Let the above introduced separation modes be incorporated into a Q6 quadrilateral finite element already enriched with the standard incompatible modes. Moreover, let the separation modes be treated in the same manner as the incompatible modes (for the theoretical and implementation details on incompatible modes we refer to \cite{Ibrahimbegovic-book}).

	The interpolation for the displacements $\boldsymbol{u} = {\left[u_x,u_y\right]}^T$ 
	inside the element in terms of the nodes is:
	\begin{equation}
		\label{eq:displ}
		\boldsymbol{u}(\boldsymbol{\xi},\Gamma^{e}) = 
		\sum_{a=1}^{4} {N_a(\boldsymbol{\xi}) {\boldsymbol{d}}_a} + 
		\sum_{b=1}^{2} {M_b(\boldsymbol{\xi}) {\boldsymbol{\rho}}_b} + 
		\sum_{i=1}^{2} {{\boldsymbol{p}}_i (\boldsymbol{\xi}, {\Gamma}^{e}) {\alpha }_i} ,
	\end{equation}
	where ${\boldsymbol{d}}_a={\left[u_{xa},u_{ya}\right]}^T$ are the displacements of node $a$, $M_1\left(\xi \right)=1-{\xi }^2$ and $M_2\left(\eta \right)=1-{\eta }^2$ are the interpolation functions for the incompatible modes $\vek{\rho}_b$, while ${\boldsymbol{p}}_i={\left[p_{ix},p_{iy}\right]}^T$ are the interpolation functions for the crack opening parameters ${\alpha}_i$. The crack opening interpolation functions ${\boldsymbol{p}}_i$ are obtained from the condition that ${\Omega}^{e-}$ should remain still when ${\Omega}^{e+}$ is rigidly displaced from it for a basic separation mode (at assuming ${\boldsymbol{\rho}}_b=\boldsymbol{0}$), see \cite{Dujc et al 2010-art}:
	\begin{equation}
		\label{eq:ED Interp Func}
		{\boldsymbol{p}}_1 = \left(\mrm{H}_{\Gamma}\left(\boldsymbol{x}\right) - \sum_{a \in {\Omega}^{e+}} {N_a} \right) \boldsymbol{n} {, \ \ \ \ \ \ \ } 
		{\boldsymbol{p}}_2 = \left(\mrm{H}_{\Gamma}\left(\boldsymbol{x}\right) - \sum_{a \in {\Omega}^{e+}} {N_a} \right) \boldsymbol{m} .
	\end{equation}
	Here, $\boldsymbol{n}={\left[n_x,n_y\right]}^T$ and $\boldsymbol{m}={\left[m_x,m_y\right]}^T={\left[{-n}_y,n_x\right]}^T$ are the unit normal and unit tangent to ${\Gamma}^{e}$, respectively, and $\mrm{H}_{\Gamma}$ is the Heaviside function, equal to unity for
	$\boldsymbol{x} \in {\Omega}^{e+}$ and vanishing otherwise, which inserts a strong discontinuity in the displacements at ${\Gamma}^e$.
	
	%\noindent 
	The infinitesimal strains are the symmetric gradient of \eqref{eq:displ}, thus written in vector form:
	\begin{equation}
		\label{eq:strains}
		\boldsymbol{\epsilon} \left(\boldsymbol{\xi},{\Gamma}^{e} \right) = 
		\sum_{a=1}^{4} {{\boldsymbol{B}}_a \left(\boldsymbol{\xi} \right){\boldsymbol{d}}_a} + 
		\sum_{b=1}^{2} {{\widetilde{\boldsymbol{G}}}_b \left(\boldsymbol{\xi} \right) {\boldsymbol{\rho}}_b} + 
		\sum_{i=1}^{2} {{\boldsymbol{G}}_i\left(\boldsymbol{\xi},{\Gamma}^{e}\right){\alpha }_i},
	\end{equation}
	where 
	\begin{equation}
		\label{eq:Goper}
		{{\boldsymbol{B}}_a = \left[\begin{array}{cc}
			\frac{\partial N_a}{\partial x} & 0 \\ 
			0 & \frac{\partial N_a}{\partial y} \\ 
			\frac{\partial N_a}{\partial y} & \frac{\partial N_a}{\partial x} 
			\end{array}
			\right]}
		,\ \ \ %\quad
		{\widetilde{\boldsymbol{G}}}_b = {\boldsymbol{G}}_b - \frac{1}{A_{\Omega^e}} \int_{\Omega^e}{{\boldsymbol{G}}_b}{ \, \di \Omega} 
		,\ \ \
		{\boldsymbol{G}}_b = \left[\begin{array}{cc}
		\frac{\partial M_b}{\partial x} & 0 \\ 
		0 & \frac{\partial M_b}{\partial y} \\ 
		\frac{\partial M_b}{\partial y} & \frac{\partial M_b}{\partial x} \end{array}
		\right] ,
	\end{equation}
	with $A_{\Omega^e}$ the area of the element. 
	Note that ${\widetilde{\boldsymbol{G}}}_b$ is designed so that $\int_{\Omega^e}{{\widetilde{\boldsymbol{G}}}_b}{\di\mathrm{\Omega}}=\boldsymbol{0}$, which enforces element-wise orthogonality of the enhanced strains (due to incompatible modes) to constant stresses over ${\Omega}^e$. The expressions for ${\boldsymbol{G}}_i
	% = {\left[\frac{\partial p_{ix}}{\partial x},\frac{\partial p_{iy}}{\partial y},\frac{\partial p_{ix}}{\partial y}\mathrm{+}\frac{\partial p_{iy}}{\partial x}\right]}^T
	$ in \eqref{eq:strains} are
	\begin{equation}
		\label{eq:GoperExt}
		{\boldsymbol{G}}_1 = 
		\underbrace{-\sum_{a \in 	{\Omega}^{e+}}{{\boldsymbol{B}}_a\boldsymbol{n}}}_{{\overline{\boldsymbol{G}}}_1} + 
		\underbrace{\ops{\updelta}_{\Gamma}\left(\boldsymbol{x}\right) {\boldsymbol{B}}_n\boldsymbol{n}}_{{\overline{\overline{\boldsymbol{G}}}}_1}
		,\ \
		{\boldsymbol{G}}_2 = 
		\underbrace{-\sum_{a\in {\Omega}^{e+}}{{\boldsymbol{B}}_a\boldsymbol{m}}}_{{\overline{\boldsymbol{G}}}_2} + 
		\underbrace{\ops{\updelta}_{\Gamma}\left(\boldsymbol{x}\right){\boldsymbol{B}}_n\boldsymbol{m}}_{{\overline{\overline{\boldsymbol{G}}}}_2}	
		,\ \
		{{\boldsymbol{B}}_n = \left[\begin{array}{cc}
			n_x & 0 \\ 
			0 & n_y \\ 
			n_y & n_x 
			\end{array}
			\right]} ,
	\end{equation}
	where $\ops{\updelta}_{\Gamma} \left(\boldsymbol{x}\right)$ is the Dirac-$\updelta$ sitting on the set ${\Gamma}^e$ and defined by $\frac{\partial \mrm{H}_{\Gamma}\left(\boldsymbol{x}\right)}{\partial x}=\ops{\updelta}_{\Gamma} \left(\boldsymbol{x}\right)n_x$ and $\frac{\partial \mrm{H}_{\Gamma}\left(\boldsymbol{x}\right)}{\partial y}=\ops{\updelta}_{\Gamma} \left(\boldsymbol{x}\right)n_y$. 
	By taking into account the additive decomposition of ${\boldsymbol{G}}_i$ 
	(see \eqref{eq:GoperExt}) the strains in \eqref{eq:strains} can be rewritten as
	\begin{equation}
		\label{eq:strainsExt}
		\boldsymbol{\epsilon} = 
		\underbrace{
			\sum_{a=1}^{4}{{\boldsymbol{B}}_a {\boldsymbol{d}}_{a}} + 
			\sum_{b=1}^{2}{{\widetilde{\boldsymbol{G}}}_b\ {\boldsymbol{\rho}}_b} + 
			\sum_{i=1}^{2}{{\overline{\boldsymbol{G}}}_i{\alpha }_i}}_{\overline{\boldsymbol{\epsilon}}
		} + 
		\underbrace{
			\sum_{i=1}^{2}{{\overline{\overline{\boldsymbol{G}}}}_i{\alpha}_i}}_{\overline{\overline{\boldsymbol{\epsilon}}}
		} ,
	\end{equation}
	where $\overline{\boldsymbol{\epsilon}}$ are the regular bounded bulk strains, and 
	$\overline{\overline{\boldsymbol{\epsilon}}}$ are singular unbounded strains at ${\Gamma}^e$. Note that the area integral of the Heaviside function derivatives is (can be checked with e.g.\ 
	Mathematica \cite{Wolfram -man})
	\begin{equation}
		\label{eq:DeltaFun}
		\int_{\Omega^e}{\frac{\partial \mrm{H}_{\Gamma}\left(\boldsymbol{x}\right)}{\partial x_j}\di\Omega} = n_{x_j} \int_{\Omega^e}{{\ops{\updelta}}_{\Gamma}\left(\boldsymbol{x}\right)\di\Omega} = n_{x_j}\int_{\Gamma^e}{\di \Gamma} = n_{x_j} l_{\Gamma}
		,\quad % \ \ \ \ \ 
		j=1,2,\ \ x_1=x,\ \ x_2=y ,
	\end{equation}
	which %implies the following conclusion on 
	gives for the area integral of a function $f\left(\boldsymbol{x}\right)$ multiplied by $\ops{\updelta}_{\Gamma} \left(\boldsymbol{x}\right)$ the usual distribution definition of
	the \emph{line} Dirac-$\updelta$:
	\begin{equation}
		\label{eq:DeltaFun2}
		\int_{\Omega^e}{{\ops{\updelta}}_{\Gamma}\left(\boldsymbol{x}\right) f\left(\boldsymbol{x}\right)\di\Omega = 
			\int_{\Gamma^e}{f\left(\boldsymbol{x}\right) \di\Gamma}} .
	\end{equation}
	
	%\noindent 
	For the virtual strains, we adopt the same interpolation as for the real strains \eqref{eq:strains}, except for the part due to the separation modes:
	\begin{equation}
		\label{eq:virtStrains}
		\widehat{\boldsymbol{\epsilon}}\left(\boldsymbol{\xi},{\Gamma}^{e}\right) = 
		\underbrace{
			\sum_{a=1}^{4}{{\boldsymbol{B}}_a \left(\boldsymbol{\xi}\right) {\widehat{\boldsymbol{d}}}_a}}_{{\widehat{\overline{\boldsymbol{\epsilon}}}}^d} + 
		\underbrace{
			\sum_{b=1}^{2}{{\widetilde{\boldsymbol{G}}}_b\left(\boldsymbol{\xi}\right){\widehat{\boldsymbol{\rho}}}_b} + 
			\sum_{i=1}^{2}{{\widehat{\boldsymbol{G}}}_i\left(\boldsymbol{{\xi}},\Gamma^{e}\right){\widehat{\alpha}}_i}
		}_{{\widehat{\overline{\boldsymbol{\epsilon}}}}^{\beta}} .
	\end{equation}
	Here, ${\widehat{\overline{\boldsymbol{\epsilon}}}}^d$ are virtual strains due to virtual nodal displacements ${\widehat{\boldsymbol{d}}}_a$, ${\widehat{\overline{\boldsymbol{\epsilon}}}}^{\beta}$ are virtual strains due to the enhancement with both incompatible and separation modes, ${\widehat{\boldsymbol{\rho}}}_b$ are virtual incompatible-mode parameters, and ${\widehat{\alpha}}_i$ are virtual separation-mode parameters. The ${\widehat{\boldsymbol{G}}}_i$ in \eqref{eq:virtStrains} is designed so that $\int_{\Omega^e}{{\widehat{\boldsymbol{G}}}_i}{{\di\Omega}} = \boldsymbol{0}$ for constant stresses over $\Omega^e$:
	\begin{equation}
		\label{eq:virtGoper}
		{\widehat{\boldsymbol{G}}}_i={\boldsymbol{G}}_i - \frac{1}{A_{{\Omega}^e}} \int_{\Omega^e}{{\boldsymbol{G}}_i}{{\di\Omega}} 
		\  \  \  \Rightarrow  \  \  \
		\int_{\Omega^e}{{\widehat{\boldsymbol{G}}}_i}{{\di\Omega}}=\boldsymbol{0} ,
	\end{equation}
	leading to (see \eqref{eq:GoperExt})
	\begin{align}
		\label{eq:virtGoper1Ext}
		{\widehat{\boldsymbol{G}}}_1 &= \underbrace{{\overline{\boldsymbol{G}}}_1 - \frac{1}{A_{\Omega^e}}\int_{\Omega^e}{{\overline{\boldsymbol{G}}}_1}{{\di\Omega}} - \frac{l_{\Gamma}}{A_{\Omega^e}} {\boldsymbol{B}}_n \boldsymbol{n}}_{{\overline{\widehat{\boldsymbol{G}}}}_1} + \underbrace{\ops{\updelta}_{\Gamma} \left(\boldsymbol{x}\right){\boldsymbol{B}}_n\boldsymbol{n}}_{{\overline{\overline{\widehat{\boldsymbol{G}}}}}_1} , \\
%	\end{equation}
%and
%	\begin{equation}
		\label{eq:virtGoper2Ext}
		{\widehat{\boldsymbol{G}}}_2 &= \underbrace{{\overline{\boldsymbol{G}}}_2 - \frac{1}{A_{{\Omega}^e}} \int_{\Omega^e}{{\overline{\boldsymbol{G}}}_2}{{\di\Omega}}-\frac{l_{\Gamma}}{A_{\Omega^e}}{\boldsymbol{B}}_n\boldsymbol{m}}_{{\overline{\widehat{\boldsymbol{G}}}}_2}+\underbrace{\ops{\updelta}_{\Gamma} \left(\boldsymbol{x}\right){\boldsymbol{B}}_n\boldsymbol{m}}_{{\overline{\overline{\widehat{\boldsymbol{G}}}}}_2} .
	\end{align}
	From \eqref{eq:virtGoper1Ext}--\eqref{eq:virtGoper2Ext} it can be observed that 
	${\widehat{\boldsymbol{G}}}_i$ can be also additively decomposed, with ${\overline{\widehat{\boldsymbol{G}}}}_i$ being regular and bounded, and ${\overline{\overline{\widehat{\boldsymbol{G}}}}}_i$ including the singular $\ops{\updelta}_{\Gamma} \left(\boldsymbol{x}\right)$. The design of ${\widehat{\boldsymbol{G}}}_i$ in \eqref{eq:virtGoper} enforces element-wise orthogonality of virtual enhanced strains (due to separation modes) and  constant stresses over $\Omega^e$.

	\subsection{Stationary point equations}
	The variational equations are derived from the modified Hu-Washizu functional that takes into account enhanced strains \eqref{eq:virtStrains} and orthogonality of stresses and enhanced strains, see e.g.\ \cite{Ibrahimbegovic-book}. We will omit the derivation and elaborate only on the resulting variational (i.e.\ stationary point) equations. 
	
	%\noindent 
	The first variational equation for a discretised 2d solid with $N_e$ elements is 
	\begin{equation}
		\Aop_{e=1}^{N_e}\left(G^{int,e}-G^{ext,e}\right) = 0 ,
	\end{equation}
	where $\Aop$ is the finite element assembly operator, $G^{ext,e}=\sum_{a=1}^{4}{{\widehat{\boldsymbol{d}}}^T_a{\boldsymbol{f}}^{ext,e}_a}$ is the virtual work of external forces
	on the element, and $G^{int,e}=t^e\int_{\Omega^e}{{\widehat{\overline{\boldsymbol{\epsilon }}}}^{d,T}\vek{C}\,\overline{\boldsymbol{\epsilon}}\,\di\Omega}$ is the virtual work of internal forces in the element. By taking into account \eqref{eq:strainsExt} and \eqref{eq:virtStrains}, 
	and assuming linearly elastic material in the bulk, the latter can be rewritten as:
	\begin{equation}
		\label{eq:intWork}
		G^{int,e} = t^e \int_{\Omega^e}
		{
			{\left(\sum_{a=1}^{4}{{\boldsymbol{B}}_a{\widehat{\boldsymbol{d}}}_a}\right)}^T \vek{C} \left(
			\sum_{c=1}^{4}{{\boldsymbol{B}}_c{\boldsymbol{d}}_c} + 
			\sum_{b=1}^2{{\widetilde{\boldsymbol{G}}}_b\ {\boldsymbol{\rho }}_b} + 
			\sum_{i=1}^{4}{{\overline{\boldsymbol{G}}}_i{\alpha }_i}\right)\, \di\Omega  ,
		}
	\end{equation}
	where $\vek{C}$ is the constitutive matrix for linear elasticity for the bulk, and $t^e$ is the element thickness, assumed to be constant. Thus, the first variational equation yields the following set of equations
	\begin{equation}
		\label{eq:GlobalEq}
		\Aop_{e=1}^{N_e}\underbrace{\left({\boldsymbol{f}}^{int,e} - {\boldsymbol{f}}^{ext,e}\right)}_{{\boldsymbol{R}}^e_d} = \boldsymbol{0} ,
	\end{equation}
	where ${\vek{f}}^{ext,e} = {\left[{\vek{f}}^{ext,e,T}_a\right]}_{(a=1,\dots,4)}^T$ and 
	$\vek{f}^{int,e} = {\left[{\vek{f}}^{int,e,T}_a\right]}_{(a=1,\dots,4)}^T$, 
	with the structure of the latter being obvious from \eqref{eq:intWork}. 
	
	%\noindent 
	The second variational equation is 	$G^{e,2}=t^e\int_{\Omega^e}{{\left({\widehat{\overline{\boldsymbol{\epsilon}}}}^{\beta }\right)}^T  \vek{C} \,\overline{\boldsymbol{\epsilon}}\,\di\Omega}=0$; and with 	 ${\widehat{\overline{\boldsymbol{\epsilon }}}}^{\beta}$ and $\overline{\boldsymbol{\epsilon}}$ given in \eqref{eq:virtStrains} and \eqref{eq:strainsExt}, respectively, this yields:
	\begin{align}
		\label{eq:intWork2i}
%		\begin{split}
		G^{e,2} &= t^e \int_{\Omega^e}{\left(
				\sum_{b=1}^{2}{{\widetilde{\boldsymbol{G}}}_b\ {\widehat{\boldsymbol{\rho}}}_b} + \sum_{i=1}^{4}{\left({\overline{\widehat{\boldsymbol{G}}}}_i + {\overline{\overline{\widehat{\boldsymbol{G}}}}}_i\right){\widehat{\alpha}}_i}\right)}^T
				\vek{\sigma} \,\di\Omega =0, \quad \text{ with} \\
				\label{eq:intWork2-sig}
			\vek{\sigma} &= \vek{C} \left(
				\sum_{c=1}^{4} {{\boldsymbol{B}}_c{\boldsymbol{d}}_c} + 
				\sum_{b=1}^{2} {{\widetilde{\boldsymbol{G}}}_b\ {\boldsymbol{\rho}}_b} + 
				\sum_{i=1}^{4} {{\overline{\boldsymbol{G}}}_i{\alpha}_i}\right)				 
		\quad \text{ the elastic stress}.
%		\end{split}
	\end{align}
	In the modified Hu Washizu functional, the stresses vanish due to the assumed orthogonality	between stresses and enhanced strains. However, we will assume that the element stresses 	$\boldsymbol{\sigma}$ can still be computed as in \eqref{eq:intWork2-sig}. This enables us to elaborate on the second term of \eqref{eq:intWork2i}. By applying \eqref{eq:DeltaFun2}, \eqref{eq:virtGoper1Ext}, and \eqref{eq:virtGoper2Ext}, one obtains
	\begin{align}
		\label{eq:LocalEq1} 
		t^e\int_{\Omega^e}{{\overline{\widehat{\boldsymbol{G}}}}^T_1 \boldsymbol{\sigma}\, \di\Omega}
		 + t^e\int_{\Gamma^e}{{\boldsymbol{n}}^T{\boldsymbol{B}}^T_n \boldsymbol{\sigma}|_{\Gamma^e}
		 \, \di\Gamma}&=0 \\
%	\end{equation}
%	
%	\begin{equation}
		\label{eq:LocalEq2} 
		t^e\int_{\Omega^e}{{\overline{\widehat{\boldsymbol{G}}}}^T_2 \boldsymbol{\sigma}\, \di\Omega}
		+ t^e\int_{\Gamma^e}{{\boldsymbol{m}}^T{\boldsymbol{B}}^T_n \boldsymbol{\sigma}|_{\Gamma^e}
		\, \di\Gamma}&=0 .
	\end{align}
	Let us recall that the rules for the transformation of planar stresses are ${\sigma}_{nn} = {\boldsymbol{n}}^T {\boldsymbol{B}}^T_n \boldsymbol{\sigma} = {\boldsymbol{n}}^T \boldsymbol{Sn}$ and ${\sigma}_{nm} = {\boldsymbol{m}}^T {\boldsymbol{B}}^T_n \boldsymbol{\sigma} = {\boldsymbol{m}}^T \boldsymbol{Sn}$, where $\boldsymbol{S} = \left[{\left[{\sigma }_{xx},\, {\sigma }_{xy}\right]}^T,\, {\left[{\sigma}_{xx},\ {\sigma}_{yy}\right]}^T\right]$ is the stress tensor in matrix form. On the edge of $\Omega^{e-}$ (with unit normal $\boldsymbol{n}$ and unit tangent $\boldsymbol{m}$), the bulk normal and shear stresses are equal to the boundary tractions on that edge: ${\sigma}_{nn}|_{\Omega^{e-}edge} = {\boldsymbol{n}}^T {\boldsymbol{B}}^T_n 	 \boldsymbol{\sigma}|_{\Omega^{e-}edge}$ and ${\sigma}_{nm}|_{\Omega^{e-}edge} = 	 {\boldsymbol{m}}^T {\boldsymbol{B}}^T_n \boldsymbol{\sigma}|_{\Omega^{e-}edge}$. Similarly, on the edge of $\Omega^{e+}$ (with unit normal $\check{\boldsymbol{n}} = -\boldsymbol{n}$ and unit tangent $\check{\boldsymbol{m}} = -\boldsymbol{m}$), one has ${\sigma }_{\check{n}\check{n}}|_{\Omega^{e+}edge} = \left({\boldsymbol{-}\boldsymbol{n}}^T\right)\left({\boldsymbol{-}\boldsymbol{B}}^T_n\right)	 \boldsymbol{\ }\boldsymbol{\sigma}|_{\Omega^{e+}edge}$ and ${\sigma }_{\check{n}\check{m}}|_{\Omega^{e+}edge} = \left({\boldsymbol{-}\boldsymbol{m}}^T\right)	 \left({\boldsymbol{-}\boldsymbol{B}}^T_n\right) \boldsymbol{\sigma }|_{\Omega^{e+}edge}$. The tractions on both sides of the crack are equal, and are denoted as $t_n$ and $t_m$. By taking the above into account, it can be concluded that \eqref{eq:LocalEq1} and \eqref{eq:LocalEq2} are local equilibrium equations that project (in a weak sense) the bulk stresses onto tractions in the crack: 
	\begin{align}
		\label{eq:LocalEqh1}
		\underbrace{t^e \int_{\Omega^e}{{\overline{\widehat{\boldsymbol{G}}}}^T_1 \boldsymbol{\sigma}
		 \,\di\Omega} + t^e \int_{\Gamma^e}{t_n\, \di\Gamma}}_{h^e_1}&=0 , \text{ and}\\
%	\end{equation}
%	
%	\begin{equation}
		\label{eq:LocalEqh2}
		\underbrace{t^e \int_{\Omega^e}{{\overline{\widehat{\boldsymbol{G}}}}^T_2 
		\boldsymbol{\sigma}\, \di\Omega} + t^e \int_{\Gamma^e}{t_m \,\di\Gamma}}_{h^e_2}&=0 .
	\end{align}
	Moreover, the first sum in \eqref{eq:intWork2i} gives
	\begin{equation}
		\label{eq:LocalEq3}
		\underbrace{t^e \int_{{\Omega}^e}{{\widetilde{\boldsymbol{G}}}^T_b 
		\boldsymbol{\sigma} d\mathrm{\Omega}}}_{{\boldsymbol{R}}_{\rho,b}} = \boldsymbol{0},
		\quad  b=1,2 ,
	\end{equation}
	and for the latter use we introduce the notation ${\boldsymbol{R}}_{\rho}^{e T} = 
	{\left[{\boldsymbol{R}}^T_{\rho,1},{\boldsymbol{R}}^T_{\rho ,2}\right]}$. 

	%\noindent
	Thus, to summarise: the governing equations for the embedded-strong-discontinuity Q6 elements with incompatible modes are: the global equilibrium equation \eqref{eq:GlobalEq}, the local (i.e.\ element) equilibrium equations \eqref{eq:LocalEqh1} and \eqref{eq:LocalEqh2}, and the additional condition due to the incompatible modes \eqref{eq:LocalEq3}.

	\subsection{Constitutive relations} 
	\label{sec:Cons_rel}
	For the bulk, we have assumed elastic material with $\vek{C}$ being the plane stress or plane strain constitutive matrix. In order to account for the dissipative phenomena on the embedded cohesive interface (i.e.\ the crack), a traction-separation law is introduced. To that end, consider two uncoupled damage laws for cohesion in the interface, related to opening in $\boldsymbol{n}$ direction and sliding in $\boldsymbol{m}$ direction (see Fig.~\ref{fig:Q4_SepModes} and \cite{Bude 2015-phd}), formulated in the framework of thermodynamics of interfaces.
	
	%\noindent 
	As for mode I opening, the corresponding Helmholtz free energy is taken as
	\begin{equation}
		{\psi }_n\left({\overline{\overline{u}}}_n,\overline{\overline{Q}}_n,
		{\overline{\overline{\xi}}}_n\right) = 
		\frac{1}{2}{\overline{\overline{Q}}}^{-1}_n{\overline{\overline{u}}}^2_n + 
		\overline{\overline{{\Xi}}}\left({\overline{\overline{\xi}}}_n\right),
	\end{equation} 
	where ${\overline{\overline{Q}}}_n$ (an internal phenomenological variable) and ${\overline{\overline{u}}}_n$ are the compliance and opening, respectively, the softening potential $\overline{\overline{{\Xi}}}\left({\overline{\overline{\xi}}}_n\right)$ controls softening and is used implicitly in equation \eqref{eq:StressLikeVar}, and ${\overline{\overline{\xi}}}_n$ is an internal phenomenological variable for softening. The failure function for softening, which deals with the unilaterality of mode I failure (only ${\overline{\overline{u}}}_n\ge 0$ is allowed), is defined as:
	\begin{equation}
		\label{eq:FailureFunModeI}
		{\overline{\overline{\phi}}}_n\left(t_n,{\overline{\overline{q}}}_n\right) = 
		t_n - \left({\sigma}_{un} - {\overline{\overline{q}}}_n\right)\le 0 ,
	\end{equation}
	where ${\sigma}_{un}$ is the tensile strength, and ${\overline{\overline{q}}}_n$ is 
	the softening traction (thermodynamic force) associated with ${\overline{\overline{\xi}}}_n$ as
	\begin{equation}
		\label{eq:StressLikeVar}
		{\overline{\overline{q}}}_n\left({\overline{\overline{\xi}}}_n\right) = 
		-\frac{\partial {\psi}_n}{\partial {\overline{\overline{\xi}}}_n} = {
		\sigma}_{un}\left(1-\exp\left[{-\frac{{\sigma}_{un}}{G_{fn}}{\overline{\overline{\xi}}}_n}
		\right]\right),
	\end{equation}
	where $G_{fn}$ denotes the mode I fracture energy. The state equation for the normal traction 
	(the thermodynamic force associated to ${\overline{\overline{u}}}_n$) is:
	\begin{equation}
		\label{eq:TractionStress}
		t_n = \frac{\partial {\psi}_n}{\partial {\overline{\overline{u}}}_n} = 
		{\overline{\overline{Q}}}^{-1}_n{\overline{\overline{u}}}_n .
	\end{equation}
	We choose a non-associative model for the evolution of ${\overline{\overline{\xi}}}_n$
	and ${\overline{\overline{Y}}}_n$ (see \eqref{eq:EvolEqDamCompl}), 
	with the evolution equations given through the dissipation pseudo-potential $F_n$:
	\begin{equation}
		\label{eq:DissPotential}
		F_n\left(t_n,{\overline{\overline{q}}}_n\left({\overline{\overline{\xi }}}_n\right)\right) 
		= t_n\left(2-\ln\left[\frac{t_n}{{\sigma}_{un}}\right]\right) - 
		\left({\sigma }_{un}-{\overline{\overline{q}}}_n
		\left({\overline{\overline{\xi}}}_n\right)\right),
	\end{equation}
	so that the evolution equations for the internal variables are:
\correction{	\begin{align}
		\label{eq:EvolEqCrackOp}
		{\dot{\overline{\overline{\xi}}}}_n &= {\dot{\overline{\overline{\gamma}}}}_n\frac{\partial F_n}{\partial {\ \overline{\overline{q}}}_n} = {\dot{\overline{\overline{\gamma}}}}_n , \\
%	\end{equation}
%	
%	\begin{equation}
		\label{eq:EvolEqDamCompl}
		{\dot{\overline{\overline{Q}}}}_n &= {\dot{\overline{\overline{\gamma}}}}_n\frac{\partial F_n}{\partial {\overline{\overline{Y}}}_n} = {\dot{\overline{\overline{\gamma}}}}_n\frac{\partial F_n}{\partial t_n}\frac{\partial t_n}{\partial {\overline{\overline{Y}}}_n} = \frac{{\dot{\overline{\overline{\gamma}}}}_n}{t_n}\left(1+\frac{{\sigma}_{un}}{G_{fn}}{\dot{\overline{\overline{\xi}}}}_n\right),
	\end{align}}
	where \correction{${\dot{\overline{\overline{\gamma}}}}_n$} is the Lagrange damage multiplier for the normal direction, and ${\overline{\overline{Y}}}_n$ is the thermodynamically dual variable associated with the compliance ${\overline{\overline{Q}}}_n$ as ${\overline{\overline{Y}}}_n = -\frac{\partial {\psi}_n}{\partial {\overline{\overline{Q}}}_n} = \frac{1}{2}t^2_n$. These equations come with the Kuhn-Tucker optimality and consistency condition:
\correction{	\begin{equation}
		\label{eq:KuhnTuckerCond}
		{\dot{\overline{\overline{\gamma}}}}_n\ge 0, \quad
		{\overline{\overline{\phi}}}_n\le 0, \quad
		{\dot{\overline{\overline{\gamma}}}}_n{\overline{\overline{\phi}}}_n=0,\;\text{ and }\;
		{\dot{\overline{\overline{\gamma}}}}_n{\dot{\overline{\overline{\phi}}}}_n = 0 .
	\end{equation}}
	%	respectively. 
	
\begin{figure}[!htb]
		\centering
		\includegraphics[height=0.35\textwidth]{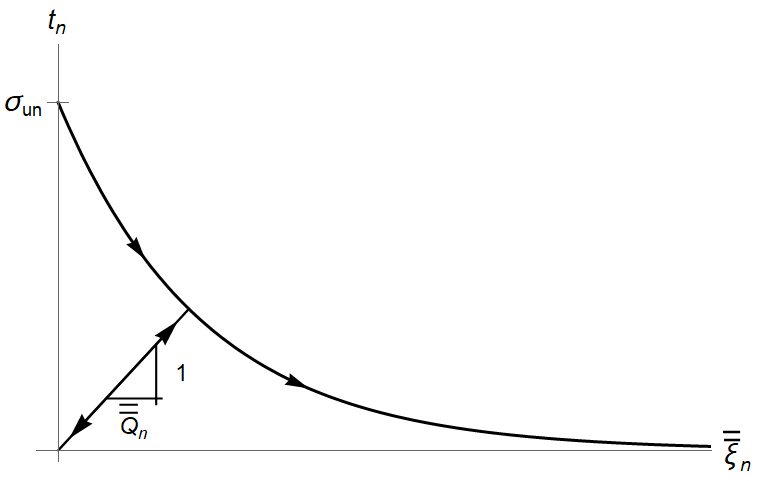}
		\caption{Rigid-damage exponential softening in normal direction.} 
		\label{fig:ExpSoft}
	\end{figure}

		%\noindent 
	A graphic representation of the resulting traction-separation damage law for mode I opening is shown in Fig.~\ref{fig:ExpSoft}.
	
	%\noindent 
	For the tangent direction $\boldsymbol{m}$, which is related to mode II sliding, a very similar damage model is considered. In contrast with the normal direction, where the unilaterality of mode I opening has to be taken into account, this is not needed for mode II, and the failure criterion for the sliding mode II is:
	\begin{equation}
		\label{eq:FailureFunModeII}
		{\overline{\overline{\phi}}}_m\left(t_m,{\overline{\overline{q}}}_m\right) = \left|t_m\right|-\left({\sigma}_{um} - {\overline{\overline{q}}}_m\right)\le 0 .
	\end{equation}
	Note that for mode II opening, all previously presented relations \eqref{eq:FailureFunModeI} -- \eqref{eq:KuhnTuckerCond} still hold, provided that index $n$ denoting the normal direction $\boldsymbol{n}$ is replaced by the index $m$ denoting the tangent direction $\boldsymbol{m}$. An exception is \eqref{eq:FailureFunModeI}, which is replaced by \eqref{eq:FailureFunModeII}. The local displacement jumps ${\overline{\overline{u}}}_n$ and ${\overline{\overline{u}}}_m$ equal the element crack parameters:
	\begin{equation}
		\label{eq:CrackPar}
		{\overline{\overline{u}}}_n = {\alpha}_1,\quad {\overline{\overline{u}}}_m = {\alpha}_2 .
	\end{equation}

	\subsection{Crack nucleation}
	\label{sec:ED_Criterion}
	The criterion for embedding a discontinuity in the finite element $e$ under consideration is based on the maximum principal stress. The procedure is summarised in Algorithm~\ref{alg:EmbeddingDisc}. After the solution for the loading increment (or time increment) has converged, one computes the major principal stress for all four integration points in the bulk:
	\begin{align}
	\label{eq:PrincipStress}
	\sigma_{p}^{e,i}&=\max \left(\sigma_{p1}^{e,i},\sigma_{p2}^{e,i}\right), \quad i = 1,\dots,4 , \\
%	\end{equation}
%	\begin{equation}
	\label{eq:PrincipStress2}
	{\sigma_{p1,p2}^{e,i}} &= \frac{{\sigma}_{xx}^{e,i}+{\sigma}_{yy}^{e,i}}{2} 
	\pm \sqrt{\left( \frac{{\sigma}_{xx}^{e,i}-{\sigma}_{yy}^{e,i}} {2} \right)^2 + 
	\left( {\sigma}_{xy}^{e,i} \right)^2} .
	\end{align}
	
	Having obtained the major principal stresses in the bulk Gauss-points, the maximum
	\begin{equation}
	\label{eq:MaxPrincipStress}
	\sigma_{p}^{e} = \underset{i = 1,\dots,4}{\max} \sigma_{p}^{e,i}
	\end{equation}
	is compared with the tensile strength $\sigma_{un}$. If $\sigma_{p}^{e} \geq \sigma_{un}$, the discontinuity is embedded into the finite element $e$.

	\begin{figure}[!hbt]
		\centering
		\includegraphics[width=14cm]{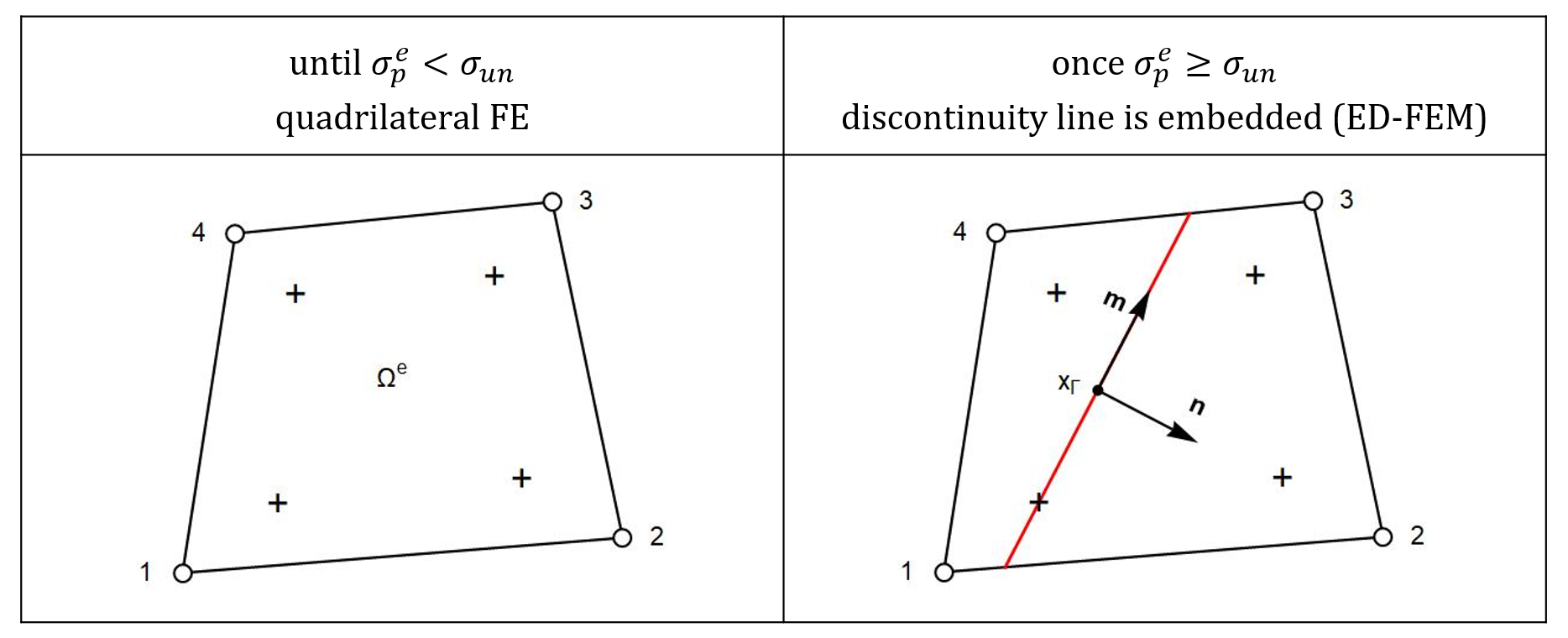}
		\caption{Crack nucleation criterion.} 
		\label{fig:CrackNuclCrit}
	\end{figure}
	The position and orientation of the discontinuity line are determined by using point $\boldsymbol{x}_{ED}^{e}$ and angle $\alpha_{ED}^{e}$ described below. Namely, the discontinuity line passes through $\boldsymbol{x}_{ED}^{e}$, and the inclination of its normal $\boldsymbol{n}$ with respect to the global $x$ axis is defined by angle $\alpha_{ED}^{e}$. Let $B^{e}$ contain the indices of those bulk integration points where the major principal stress equals (within numerical tolerance) $\sigma_{p}^{e}$:
	\begin{equation}
	\label{eq:SetOfBIP}
	B^{e} = \{ id \in \{1,2,3,4 \}: \sigma_{p1}^{e,id} = \sigma_{p}^{e} \}
	\end{equation}
	Point $\boldsymbol{x}_{ED}^{e}$ is the geometric centre of the integration points in the set $B^{e}$
	\begin{equation}
	\label{eq:EDPosition}
	\boldsymbol{x}_{ED}^{e} = \frac{1}{N_{B^{e}}} \sum_{i \in B^{e}}\boldsymbol{x}^{e,i} ,
	\end{equation}
	where $N_{B^{e}}$ is the number of elements in $B^{e}$, and $\boldsymbol{x}^{e,i}$ 
	is the position vector of the integration point $i$ in the global coordinate system. 

	The orientation of the discontinuity line is determined by the normal vector $\boldsymbol{n}$ (see Fig.~\ref{fig:CrackNuclCrit}), whose direction is defined by angle $\alpha_{ED}^{e}$. The latter is determined by computing the average stresses for points from set $B^{e}$ as
	\begin{equation}
	\label{eq:EDStress}
	\boldsymbol{\sigma}^{e}_{ED} = \frac{1}{N_{B^{e}}} \sum_{i \in B^{e}} \boldsymbol{\sigma}^{e,i}, 
	\end{equation}
	where $\boldsymbol{\sigma}^{e,i}$ are the stresses in the integration point $i$. From $\boldsymbol{\sigma}^{e}_{ED}$ one finds the principal angle that coincides with the direction of the normal vector $\boldsymbol{n}$ to the discontinuity line as:
	\begin{equation}
	\label{eq:EDAngle}
	\alpha_{ED}^{e}=\frac{1}{2} \arctan\left({\frac{2 \sigma_{ED,xy}^{e}}
	{\sigma_{ED,xx}^{e}-\sigma_{ED,yy}^{e}}}\right) .
	\end{equation}
	
	%\noindent 
	We reiterate that we do not use any crack tracking algorithm that enforces two adjacent	discontinuity lines to be geometrically connected on the common edge of two finite elements in order to make the crack path continuous on the mesh level. In the present formulation, the local cracks (embedded discontinuities) are independent of the discontinuity positions in the neighbouring finite elements. A crack in an element appears whenever the principal stress $\sigma_{p}^{e}$ in that element exceeds the material tensile strength $\sigma_{un}$. Therefore, the crack path propagates naturally according to the stress state in the solid domain. Once the discontinuity is embedded --- i.e.\ a crack is nucleated --- its position and orientation are fixed and do not change during the simulation. From the viewpoint of the finite element implementation, only one discontinuity can be embedded in one finite element. However, in a mesh, several new embedded discontinuities can be activated in one load increment (or time increment), one for each element that fulfils the crack nucleation criterion.

	\IncMargin{1em}
	\begin{algorithm}[!hbt]
		\SetKwData{Left}{left}\SetKwData{This}{this}\SetKwData{Up}{up}
		\SetKwFunction{Union}{Union}\SetKwFunction{FindCompress}{FindCompress}
		\SetKwInOut{Input}{Input}\SetKwInOut{Output}{Output}
		\Input{Stress state in finite elements and location of bulk integration points ($i$)}
		\Output{Newly activated embedded discontinuities}
		\BlankLine
		\BlankLine
		\For{$e = 1$ \KwTo $N_{el}$}{
			\eIf{finite element $e$ does not have discontinuity line}{
				\For{$i = 1$ \KwTo $4$}
				{
					Compute first principal stress $\sigma_{p}^{e,i}$ 
					(Eqs. \eqref{eq:PrincipStress} and \eqref{eq:PrincipStress2}).
				}
				Compute the maximum principal stress $\sigma_{p}^{e}$ 
				(Eq. \eqref{eq:MaxPrincipStress})\; 
				\eIf{$\sigma_{p}^{e} > \sigma_{un}$}
				{
					1. Determine position $\boldsymbol{x}_{ED}^{e}$ and orientation $\alpha_{ED}^{e}$ 
					of the discontinuity line\\ \quad (Eqs. \eqref{eq:SetOfBIP}--\eqref{eq:EDAngle})\;
					2. Activate the embedded discontinuity line in finite element $e$.
				}{
					Keep the standard finite element formulation Q6.
				}
			}{
				Keep the finite element formulation Q6 with the embedded discontinuity line.
			}
		}
		\caption{Crack nucleation algorithm.
		}
		\label{alg:EmbeddingDisc}
	\end{algorithm}\DecMargin{1em}

	Let us mention in the conclusion of this section two ``numerical tricks'' related to the use of the exponential function for the softening, applied in the code in order to get better numerical performance of the novel finite element:
	
	\noindent (i) For the purpose of continuation of the solution process immediately after embedding the discontinuity line, we use a coefficient $\overline{\overline{\kappa}}_{0}$ that has a small value. It is introduced because the compliance for rigid-damage softening law equals infinity when softening is triggered, see Fig. \ref{fig:ExpSoft}. To that end, variables $\overline{\overline{\xi}}_{n,n}$ and $\overline{\overline{Q}}_{n,n}$ (i.e. the input data for the updating algorithm described in Section~\ref{sec:UpdateDisc}) are redefined at the crack nucleation as (see Eq. \eqref{eq:StressLikeVar}):
	\begin{equation}
	\label{eq:xi_init}
	\overline{\overline{q}}_{n,n} = \sigma_{un} \overline{\overline{\kappa}}_{0} ,
	\qquad 
	\overline{\overline{\xi}}_{n,n} = -\frac{G_{fn}}{\sigma_{un}}   \text{ln}\left(1-\overline{\overline{\kappa}}_{0}\right) ,
	\qquad
	\overline{\overline{Q}}_{n,n} = \frac{\overline{\overline{\xi}}_{n,n}}{\sigma_{un} \left(1-\overline{\overline{\kappa}}_{0}\right)}
	\end{equation}
	In the same manner, we compute variables $\overline{\overline{\xi}}_{m,n}$ and $\overline{\overline{Q}}_{m,n}$ in the tangential direction. In this way,	$\overline{\overline{\kappa}}_{0}$ enables starting the update procedure for the softening variables  immediately after the discontinuity line is embedded.
	
	\noindent (ii) When the traction value comes very close to zero (i.e.\ the exponential softening function is almost zero, see Fig.~\ref{fig:ExpSoft}) this affects the solution convergence in a negative sense. For this reason do the following: when the traction $t_{n}$ drops below $t_{n,min}=k_{FullSoft} \, \sigma_{un}$, where $k_{FullSoft}$ is a small number, the value of the traction $t_{n}$ is set to $t_{n,min}$ for further inelastic loading. The corresponding consistent tangent operator (see \eqref{eq:SoftTangOperModeII}) is set to $-10^{-6}$ for inelastic loading and to $10^{-6}$ for elastic unloading. The shear traction $t_{m}$, with $\sigma_{um}$ used instead of $\sigma_{un}$, is treated in the same manner.

	\section{Solution procedure}
	\label{sec:EqSys}
	The solution %of the governing equations 
	is computed at discrete time points 
	${\tau}_1,\dots,{\tau}_n,{\tau}_{n+1},\dots,T$.
	 Let us illustrate the $\ell$-th iteration when searching for the solution at 
	 ${\tau}_{n+1}$ by the Newton-Raphson method or one of its relatives.  Note that the 
	 term ``time'' is used hereinafter for both statics and dynamics. In dynamics, it 
	 represents real time, while in statics it only refers to the level of the external 
	 loading, which is assumed to depend linearly on time.

	\subsection{Update at discontinuity}
	\label{sec:UpdateDisc}
	At the Gauss integration points with embedded discontinuity, the time update of the	tractions requires numerical solution of evolution differential equations for the	internal variables of the discrete damage model. A return-mapping strategy is applied 	to this end. For the increment $\left[{\tau}_n,{\tau}_{n+1}\right]$ and mode I opening, the given data are the crack opening ${\overline{\overline{u}}}_{n,n+1}$ and internal variables ${\overline{\overline{\xi}}}_{n,n}$, ${\overline{\overline{Q}}}_{n,n}$. The output of the update algorithm is traction $t_{n,n+1}$, internal variables ${\overline{\overline{\xi}}}_{n,n+1}$ and ${\overline{\overline{Q}}}_{n,n+1}$, and the consistent tangent modulus.

	The trial values of traction and softening force-like variable are computed as $t^{trial}_{n,n+1} = {\overline{\overline{u}}}_{n,n+1}/{\overline{\overline{Q}}}_{n,n}$ and ${\overline{\overline{q}}}^{trial}_{n,n+1} = {\overline{\overline{q}}}_n\left({\overline{\overline{\xi}}}_{n,n}\right)$, see \eqref{eq:TractionStress} and \eqref{eq:StressLikeVar}, respectively. Their admissibility is checked by the failure function \eqref{eq:FailureFunModeI} as
	\begin{equation}
		\label{eq:TrialFailureFunModeI}
		{\overline{\overline{\phi}}}_{n,n+1} = {\overline{\overline{t}}}^{trial}_{n,n+1} - 
		\left({\sigma}_{un} - {\overline{\overline{q}}}^{trial}_{n,n+1}\right)\le \ 0 .
	\end{equation}
	If the trial values are admissible,	then $t_{n,n+1} = t^{trial}_{n,n+1}$, ${\overline{\overline{\xi}}}_{n,n+1} = {\overline{\overline{\xi}}}_{n,n}$ and ${\overline{\overline{Q}}}_{n,n+1} = {\overline{\overline{Q}}}_{n,n}$, otherwise a correction is performed. The corrected values have to fulfil \eqref{eq:TrialFailureFunModeI}; thus, the nonlinear equation ${\overline{\overline{\phi}}}_{n,n+1}({\overline{\overline{\gamma}}}_{n,n+1})=0$ is formed and the Lagrange multiplier ${\overline{\overline{\gamma }}}_{n,n+1}$ is computed by iteration with Newton's method
	\begin{equation}
		\label{eq:SoftNewtonItGamma} 
		{\overline{\overline{\gamma}}}^{\left(k+1\right)}_{n,n+1} = 
		{\overline{\overline{\gamma}}}^{\left(k\right)}_{n,n+1} - 
		{\left(\diff{{\overline{\overline{\phi}}}^{(k)}_{n,n+1}}
		{{\overline{\overline{\gamma}}}^{(k)}_{n,n+1}}\right)}^{-1} 
				{\overline{\overline{\phi}}}^{\left(k\right)}_{n,n+1} , \qquad
		k = 1,\dots,\text{convergence}\quad\Rightarrow \quad{\overline{\overline{\gamma }}}_{n,n+1} ,
	\end{equation}
	thereby permitting an update of the softening variables, traction, and compliance as
	\begin{align}
		\label{eq:SoftUpdateCrackOp}
		{\overline{\overline{\xi}}}_{n,n+1} &= {\overline{\overline{\xi}}}_{n,n} + {\overline{\overline{\gamma}}}_{n,n+1} , \\
%	\end{equation}
%	
%	\begin{equation}
		\label{eq:SoftUpdateStressLikeVar}	
		{\overline{\overline{q}}}_{n,n+1} &= {\sigma}_{un}\left(1-
		\exp\left[-\frac{{\sigma}_{un}}{G_{fn}}{\overline{\overline{\xi}}}_{n,n+1}\right]\right) , \\
%	\end{equation}
%	
%	\begin{equation}
		\label{eq:SoftUpdateTraction}
		t_{n,n+1}&={\sigma}_{un}-{\overline{\overline{q}}}_{n,n+1}, \\
%	\end{equation}
%	
%	\begin{equation}
		\label{eq:SoftUpdateDamCompl}
		{\overline{\overline{Q}}}_{n,n+1} &= {\overline{\overline{Q}}}_{n,n} + 
		\frac{{\overline{\overline{\gamma}}}_{n,n+1}}{t_{n,n+1}} 
		\left(1+\frac{{\sigma}_{un}}{G_{fn}}{\overline{\overline{\xi}}}_{n,n+1}\right) .
	\end{align}
	
	%\noindent 
	For mode II opening, the procedure is the same, except that subscript $n$ is replaced by $m$ in \eqref{eq:SoftNewtonItGamma}--\eqref{eq:SoftUpdateDamCompl}, and the failure function \eqref{eq:FailureFunModeI} is replaced by \eqref{eq:FailureFunModeII}. The consistent tangent operators $\frac{\di t_{n,n+1}}{{\di \overline{\overline{u}}}_{n,n+1}}$ and $\frac{\di t_{m,n+1}}{{\di\overline{\overline{u}}}_{m,n+1}}$ are needed in computations. 
	The latter is 
	\begin{equation}
		\label{eq:SoftTangOperModeII}
		\diff{t_{m,n+1}}{{\overline{\overline{u}}}_{m,n+1}} = 
		\frac{{\left({\overline{\overline{Q}}}_{m,n+1}\right)}^{-1}}
		{1+s\, {\overline{\overline{u}}}_{m,n+1} 
		\diff{{\left({\overline{\overline{Q}}}_{m,n+1}\right)}^{-1}}
		{{\overline{\overline{\gamma}}}_{m,n+1}} 
		{\left(\diff{{\overline{\overline{q}}}_{m,n+1}}
		{{\overline{\overline{\gamma}}}_{m,n+1}}\right)}^{-1}} ,
	\end{equation}
	where $s=\mrm{sign}(t_{m,n+1})$. The former is obtained if subscript $m$ is replaced by $n$ in \eqref{eq:SoftTangOperModeII} and one sets $s=1$.

	\subsubsection{Computation of unknowns}
	When there is no crack in the mesh of $N_e$ elements, solving the equations reduces to
	the linear problem
	\begin{equation}
		\label{eq:AugmentedSystEq}
		\Aop_{e=1}^{N_e} {\left[\begin{array}{cc}
			{\boldsymbol{K}}^e_{dd} & {\boldsymbol{K}}^e_{d\rho } \\ 
			{\boldsymbol{K}}^e_{\rho d} & {\boldsymbol{K}}^e_{\rho \rho } \end{array}
			\right]}_{n+1}^\ell {\left\lbrace \begin{array}{c}
			{\mathrm{\Delta} \boldsymbol{d}}^{e,\ell}_{n+1} \\ 
			{\mathrm{\Delta} \boldsymbol{\rho}}^{e,\ell}_{n+1} \end{array}
			\right\rbrace} = \Aop_{e=1}^{N_e} {\left\lbrace \begin{array}{c}
			-{\boldsymbol{R}}^{e,\ell}_{d,n+1} \\ 
			-{\boldsymbol{R}}^{e,\ell}_{\rho,n+1} \end{array}
			\right\rbrace} ,
	\end{equation}
	where ${\boldsymbol{\rho}}^e = {\left[{\rho}_{11},{\rho}_{12},{\rho}_{21},{\rho}_{22}\right]}^T$, 
	${\boldsymbol{K}}_{dd,n+1}^{e,\ell} = {\left.\frac{\partial {\boldsymbol{R}}_{d}^{e}}
	{{\partial \boldsymbol{d}}^{e}}\right|}^\ell_{n+1}$, ${\boldsymbol{K}}_{d\rho,n+1}^{e,\ell}=
	{\left.\frac{\partial {\boldsymbol{R}}^e_d}{\partial {\boldsymbol{\rho }}^e}\right|}^\ell_{n+1}$, 
	${\boldsymbol{K}}^{e,\ell}_{\rho d,n+1}={\left.\frac{\partial {\boldsymbol{R}}^e_{\rho }}
	{\partial {\boldsymbol{d}}^e}\right|}^\ell_{n+1}$, and 
	${\boldsymbol{K}}^{e,\ell}_{\rho \rho ,n+1}={\left.\frac{\partial {\boldsymbol{R}}^e_{\rho }}
	{\partial {\boldsymbol{\rho }}^e}\right|}^\ell_{n+1}$. The ${\boldsymbol{R}}^e_d$ and 
	${\boldsymbol{R}}^e_{\rho}$ are those from \eqref{eq:GlobalEq} and \eqref{eq:LocalEq3},
	 respectively, with ${\alpha}_i=0$. The system \eqref{eq:AugmentedSystEq} is simplified 
	 by the condensation of the parameters of the incompatible modes:
	\begin{align}
		\label{eq:CondAugmentedSys}
		\upDelta\boldsymbol{\rho}^{e,\ell}_{n+1} &= -{\left[
		{\boldsymbol{K}}^{e,\ell}_{\rho \rho,n+1}\right]}^{-1} 
		\left(-{\boldsymbol{R}}^{e,\ell}_{\rho ,n+1} - 
		{\boldsymbol{K}}^{e,\ell}_{\rho d,n+1} {\upDelta\boldsymbol{d}}^{e,\ell}_{n+1}\right) ,\\
%	\end{equation}
%	\begin{equation}
		\label{eq:CondK}
%		\begin{split}
		{\boldsymbol{K}}^{e,\ell}_{c,n+1} &= {\boldsymbol{K}}^{e,\ell}_{dd,n+1} - 
		{\boldsymbol{K}}^{e,\ell}_{d\rho ,n+1} 
{\left[{\boldsymbol{K}}^{e,\ell}_{\rho \rho ,n+1}\right]}^{-1}{\boldsymbol{K}}^{e,\ell}_{\rho d,n+1},
		\\ \label{eq:CondK-2}
		{\boldsymbol{R}}^{e,\ell}_{c,n+1} &= {\boldsymbol{R}}^{e,\ell}_{d,n+1} - 
		{\boldsymbol{K}}^{e,\ell}_{d\rho ,n+1} 
		{\left[{\boldsymbol{K}}^{e,\ell}_{\rho \rho ,n+1}\right]}^{-1}
		{\boldsymbol{R}}^{e,\ell}_{\rho ,n+1}, 
%		\end{split}  
%	\end{equation}
	\end{align}
	leading to the linear problem 
 	\begin{equation}
		\label{eq:CondSys}
		\Aop_{e=1}^{N_e}\left({\boldsymbol{K}}^{e,\ell}_{c,n+1}
		{\upDelta\boldsymbol{d}}^{e,\ell}_{n+1}\right) = 
		\Aop_{e=1}^{N_e}\left(-{\boldsymbol{R}}^{e,\ell}_{c,n+1}\right) .
	\end{equation}
	
	%\noindent 
	When there is at least one element with a crack in the mesh, the problem becomes nonlinear, and
	the following system of equations needs to be solved iteratively for the nodal displacement
	increment  ${\upDelta\boldsymbol{d}}^{e,\ell}_{n+1}$:
	\begin{equation} 
		\label{eq:CondensedSys}
		\Aop_{e=1}^{N_e}\left({\boldsymbol{K}}^{e,\ell}_{c,n+1}
		{\upDelta\boldsymbol{d}}^{e,\ell}_{n+1}\right) = 
		\Aop_{e=1}^{N_e}\left(-{\boldsymbol{R}}^{e,\ell}_{c,n+1}\right) .
	\end{equation}
	For an element with no crack, the ${\boldsymbol{K}}_{c,n+1}^{e,l}$ and 
	${\boldsymbol{R}}_{c,n+1}^{e,l}$ contributions are those from \eqref{eq:CondK}
	and \eqref{eq:CondK-2}, respectively. 
	For an element with a crack, they are obtained by performing static condensation of the incompatible mode and the crack opening parameters. The following element equations are considered
	\begin{equation}
		\label{eq:AugmentedSystem2}
		{\left[ \begin{array}{ccc}
			{\boldsymbol{K}}^e_{dd} & {\boldsymbol{K}}^e_{d\rho } & {\boldsymbol{K}}^e_{d\alpha } \\ 
			{\boldsymbol{K}}^e_{\rho d} & {\boldsymbol{K}}^e_{\rho \rho } & 
			{\boldsymbol{K}}^e_{\rho \alpha } \\ 
			{\boldsymbol{K}}^e_{\alpha d} & {\boldsymbol{K}}^e_{\alpha \rho } & 
			{\boldsymbol{K}}^e_{\alpha \alpha } \end{array}
			\right]}^\ell_{n+1} \left[ \begin{array}{c}
		{\upDelta\boldsymbol{d}}^{e,\ell}_{n+1} \\ 
		{\upDelta\boldsymbol{\rho }}^{e,\ell}_{n+1} \\ 
		{\upDelta\boldsymbol{\alpha }}^{e,\ell}_{n+1} \end{array}
		\right]  = \left[  \begin{array}{c}
		{-\boldsymbol{R}}^{e,\ell}_{d,n+1} \\ 
		{-\boldsymbol{R}}^{e,\ell}_{\rho ,n+1} \\ 
		{-\boldsymbol{h}}^{e,\ell}_{n+1} \end{array} 
		\right] ,
	\end{equation}
	where some of the terms were already defined above, and
	\begin{equation}
		\label{eq:AugmentedSystem2Eq}
		\begin{split}
		& {\boldsymbol{K}}^{e,\ell}_{d\alpha ,n+1} = {\left.\frac{\partial {\boldsymbol{R}}^e_d}
		{{\partial \boldsymbol{\alpha }}^e}\right|}_{n+1}^{\ell},\qquad \qquad
		{\boldsymbol{K}}^{e,\ell}_{\rho \alpha ,n+1} = {\left.\frac{\partial {\boldsymbol{R}}^e_{\rho}}
		{\partial {\boldsymbol{\alpha }}^e}\right|}_{n+1}^{\ell} 
		\\
		& {\boldsymbol{K}}^{e,\ell}_{\alpha d,n+1} = {\left.\frac{\partial {\boldsymbol{h}}^e}
		{{\partial \boldsymbol{d}}^e}\right|}_{n+1}^\ell,\quad
		{\boldsymbol{K}}^{e,\ell}_{\alpha \rho ,n+1} = {\left.\frac{\partial {\boldsymbol{h}}^e}
		{\partial {\boldsymbol{\rho }}^e}\right|}_{n+1}^\ell,\quad
		{\boldsymbol{K}}^{e,\ell}_{\alpha \alpha ,n+1} = {\left.\frac{\partial {\boldsymbol{h}}^e}
		{{\partial \boldsymbol{\alpha }}^e}\right|}_{n+1}^\ell .
		\end{split}
	\end{equation}
	The $h_i,\ i=1,2$, in ${\boldsymbol{h}}^{e,l}_{n+1}={\left[h_1,h_2\right]}^T|^l_{n+1}$ 
	were given in \eqref{eq:LocalEqh1}--\eqref{eq:LocalEqh2}. For the sake of brevity, 
	let us rewrite \eqref{eq:AugmentedSystem2} as
	\begin{equation}
		\label{eq:AugmentedSystem2Eq2}
		{\left[ \begin{array}{cc}
			{\boldsymbol{K}}^e_{dd} & {\boldsymbol{K}}^e_{d\beta } \\ 
			{\boldsymbol{K}}^e_{\beta d} & {\boldsymbol{K}}^e_{\beta \beta } \end{array}
			\right]}^\ell_{n+1} \left[  \begin{array}{c}
		{\upDelta\boldsymbol{d}}^{e,\ell}_{n+1} \\ 
		{\upDelta\boldsymbol{\beta }}^{e,\ell}_{n+1} \end{array}
		\right] =-\left[  \begin{array}{c}
		{\boldsymbol{R}}^{e,\ell}_{d,n+1} \\ 
		{\boldsymbol{R}}^{e,\ell}_{\beta ,n+1} \end{array}
		\right] ,
	\end{equation}
	where
	\begin{equation}
		\label{eq:AugmentedSystem2Eq2Kbb}
		\begin{split}
		&{\boldsymbol{K}}_{\beta \beta ,n+1}^{e,\ell} = {\left[ \begin{array}{cc}
			{\boldsymbol{K}}_{\rho \rho }^e & {\boldsymbol{K}}_{\rho \alpha }^e \\ 
			{\boldsymbol{K}}_{\alpha \rho }^e & {\boldsymbol{K}}_{\alpha \alpha }^e \end{array}
			\right]}_{n+1}^\ell \\ 
		&{\boldsymbol{K}}_{d\beta ,n+1}^{e,\ell} = {\left[ \begin{array}{cc}
			{\boldsymbol{K}}_{d\rho }^e & {\boldsymbol{K}}_{d\alpha }^e \end{array}
			\right]}^\ell_{n+1}, \qquad 
		{\boldsymbol{K}}_{\beta d,n+1}^{e,\ell}={\left[ \begin{array}{c}
			{\boldsymbol{K}}_{\rho d}^e \\ 
			{\boldsymbol{K}}_{\alpha d}^e \end{array}
			\right]}_{n+1}^\ell ,
		\end{split}
	\end{equation}
		and
	\begin{equation}
		\label{eq:AugmentedSystem2Eq2Rbb}
		{\upDelta\boldsymbol{\beta }}_{n+1}^{e,\ell} = \left[ \begin{array}{c}
		{\upDelta\boldsymbol{\rho }}_{n+1}^{e,\ell} \\ 
		{\upDelta\boldsymbol{\alpha }}_{n+1}^{e,\ell} \end{array}
		\right]  , \qquad
		{\boldsymbol{R}}_{\beta ,n+1}^{e,\ell} = \left[  \begin{array}{c}
		{\boldsymbol{R}}^{e,\ell}_{\rho ,n+1} \\ 
		{\boldsymbol{h}}^{e,\ell}_{n+1} \end{array}
		\right] .
	\end{equation}
	Static condensation of \eqref{eq:AugmentedSystem2Eq2} yields the following equations
	\begin{equation}
		\label{eq:CondSys2}
		{\upDelta\boldsymbol{\beta }}^{e,\ell}_{n+1}= - 
		{\left[{\boldsymbol{K}}^{e,\ell}_{\beta \beta ,n+1}\right]}^{-1}
		\left(-{\boldsymbol{R}}^{e,\ell}_{\beta ,n+1}-{\boldsymbol{K}}^{e,\ell}_{\beta d,n+1}
		{\upDelta\boldsymbol{d}}^{e,\ell}_{n+1}\right)  ,
	\end{equation}
		and
	\begin{equation}
		\label{eq:CondensedSys2}
		\begin{split}
		&{\boldsymbol{K}}^{e,\ell}_{c,n+1} = {\boldsymbol{K}}^{e,\ell}_{dd,n+1} - 
		{\boldsymbol{K}}^{e,\ell}_{d\beta ,n+1}
		{\left[{\boldsymbol{K}}^{e,\ell}_{\beta \beta ,n+1}\right]}^{-1} 
		{\boldsymbol{K}}^{e,\ell}_{\beta d,n+1} \\ 
		&{\boldsymbol{R}}^{e,\ell}_{c,n+1} = {\boldsymbol{R}}^{e,\ell}_{d,n+1} - 
		{\boldsymbol{K}}^{e,\ell}_{d\beta ,n+1}
		{\left[{\boldsymbol{K}}^{e,\ell}_{\beta \beta ,n+1}\right]}^{-1}
		{\boldsymbol{R}}^{e,\ell}_{\beta ,n+1} .
		\end{split}
	\end{equation}
	For an element with a crack, the quantities from \eqref{eq:CondensedSys2} are used in \eqref{eq:CondensedSys}. By solving \eqref{eq:CondensedSys} for ${\upDelta\boldsymbol{d}}_{n+1}^{e,\ell}$ and by performing post-computations \eqref{eq:CondSys2} to get ${\upDelta\boldsymbol{\rho}}_{n+1}^{e,\ell}$ and 	${\upDelta\boldsymbol{\alpha}}_{n+1}^{e,\ell}$, new iterative values of element total displacements, incompatible mode parameters, and opening parameters are obtained
	\begin{equation}
		\label{eq:NewtItUpdate}
		{\boldsymbol{d}}_{n+1}^{e,\ell} = {\boldsymbol{d}}^{e,\ell-1}_{n+1} + 
		{\upDelta\boldsymbol{d}}^{e,\ell}_{n+1},\quad {\boldsymbol{\rho}}_{n+1}^{e,\ell} = 
		{\boldsymbol{\rho}}^{e,\ell-1}_{n+1}+{\upDelta}\boldsymbol{\rho}^{e,\ell}_{n+1}\ ,\quad
	 {\boldsymbol{\alpha }}^{e,\ell}_{n+1}={\boldsymbol{\alpha }}^{e,\ell-1}_{n+1}+
	 {\upDelta\boldsymbol{\alpha }}^{e,\ell}_{n+1} . 
	\end{equation}
	Once the convergence tolerance is reached, the values from \eqref{eq:NewtItUpdate} are accepted as solution at ${\tau}_{n+1}$, and the computation can advance to the solution at ${\tau}_{n+2}$. Otherwise, another iteration will be performed by incrementing $\ell$.

	\subsection{Extension to dynamics}
	\label{sec:Dyn}
	Extension of the above described static formulation to dynamics requires the formulation of the inertial forces and the choice of the time-stepping scheme. As for the latter, we choose the trapezoidal rule, also known as Newmark's method \cite{Ibrahimbegovic-book}. We note that the above introduced points $\tau_n$ and $\tau_{n+1}$ are in dynamics real time points in the time interval of interest $\left[0,T\right]$. 
	
	%\noindent 
	The element force vector $\boldsymbol{R}^e_d$ in \eqref{eq:GlobalEq} now becomes
	\begin{equation}
		\label{eq:DynGlobalEq}
		\boldsymbol{R}^e_d \rightarrow
		\boldsymbol{R}^{e,dyn}_d=\boldsymbol{f}^{ine,e}+\boldsymbol{f}^{int,e}-\boldsymbol{f}^{ext,e},
	\end{equation}
	where $\vek{f}^{int,e}$ and $\vek{f}^{ext,e}$ are as before, 
	and the vector of  element inertial forces is defined as
	\begin{equation}
		\label{eq:DynIntWork}
		\boldsymbol{f}^{ine,e} = 
		\underbrace{
			\left[ t^e \int_{\Omega^e} \sum_{a=1}^{4} \boldsymbol{N}_a^T	\,	\bar{\rho} \,
			\sum_{c=1}^{4}\vek{N}_c 
			\,	\di\Omega \right]}_{\boldsymbol{M}^e}  \,
		\ddot{\boldsymbol{d}} ;
	\end{equation}
	where $\bar{\rho}$ is the mass density, $\ddot{\boldsymbol{{d}}}$ is the acceleration, and $\boldsymbol{M}^e$ is the element mass matrix. After the application of the time discretisation, the velocity $\dot{\vek{d}}_{n+1}$ and acceleration $\ddot{\vek{d}}_{n+1}$ at time $\tau_{n+1}$ can expressed from the trapezoidal rule as
	\begin{equation}
		\label{eq:DynCondensedSys2}
		\begin{split}
		&  \dot{\boldsymbol{{d}}}_{n+1}=\frac{\gamma}{\beta \upDelta \tau} 
		\left (  \boldsymbol{d}_{n+1} - \boldsymbol{d}_{n} \right) -
		\frac{\gamma-\beta}{\beta} \, \dot{\boldsymbol{d}}_{n} - \frac{\gamma-2 \beta}{2 \beta} 
		\upDelta \tau \, \ddot{\boldsymbol{d}}_{n} ,  \quad \text{ and}
		\\ 
		&  
		\ddot{\boldsymbol{d}}_{n+1}=
		\frac{1}{\beta \upDelta \tau^2} 
		\left (  \boldsymbol{d}_{n+1} - \boldsymbol{d}_{n} \right)
		- \frac{1}{\beta \upDelta \tau} \,
		\dot{\boldsymbol{d}}_{n} - 
		\frac{1-2\beta}{2 \beta} \, \ddot{\boldsymbol{d}}_{n} ,
		\end{split}
	\end{equation}
	where $\gamma=1/2$, $\beta=1/4$, and $\upDelta \tau=\tau_{n+1}-\tau_n$. 
	Using \eqref{eq:DynCondensedSys2} in \eqref{eq:DynIntWork}, one obtains 
	the dynamic element force vector $\boldsymbol{R}^{e,dyn}_{d, n+1}$ that replaces $\boldsymbol{R}^{e}_{d, n+1}$ in the above iteration equations. Accordingly, the part of the stiffness matrix ${\boldsymbol{K}}_{dd,n+1}^{e}$ in those iteration equations is replaced by the dynamic stiffness ${\boldsymbol{K}}_{dd,n+1}^{e,dyn}={\boldsymbol{K}}_{dd,n+1}^{e}+	\frac{1}{\beta \upDelta \tau^2} \boldsymbol{M}^e$. With these changes, one obtains the dynamic formulation for the embedded discontinuity quadrilateral considered in this work. Note that there is no inertial force associated neither with the incompatible modes or with the crack opening parameters, and thus the static condensation carried out earlier remains unchanged.

\section{Numerical examples}
\label{sec:NumExamples}
A number of numerical examples are presented here in order to illustrate and test the derived formulation. We start with a simple uniaxial cyclic test to illustrate the functioning of the softening law described in Subsection~\ref{sec:Cons_rel}. The proposed procedure for the cracking mechanism (see Subsection~\ref{sec:ED_Criterion}) is examined by a tension test example, where the specimen is discretised both by a regular and a distorted mesh. Further, quasi-static benchmark examples, namely thermal shock, a four point bending test, and a double edge notched specimen test are shown. Finally, it is demonstrated that the proposed formulation can also be used for the numerical simulation of the \correction{dynamic fracture problems: Kalthoff's experimental test, where kinking and branching of crack paths appear, and a single edge notched specimen test for pure crack branching.} We note that the element used (sometimes referred to hereinafter as ``Q6 with ED'' or Q6ED) allows for crack opening in mode I and crack sliding in mode II. However, activation of both modes (a combined mode) is possible, with a predominant role of mode I. The additional finite element constants discussed in Subsection~\ref{sec:ED_Criterion} were chosen for the numerical examples presented below as: $\overline{\overline{\kappa}}_{0} = \{10^{-5}, 10^{-5}, 10^{-5}, 10^{-5}, 10^{-6}, 10^{-6}, 10^{-8},\}$ and $k_{FullSoft} = \{10^{-5}, 10^{-5}, 10^{-4},5 \times 10^{-6}, 10^{-5}, 5 \times 10^{-5}, 10^{-3}\}$; the values are arranged by the order of numerical examples.

Recall, that two uncoupled damage traction-separation laws to describe the tractions in the discontinuity are used. The amount of dissipated fracture energy is computed separately for mode I and mode II as: 
	\begin{align}
		\label{eq:DissFracModeI}
		W_{D, I} &= \sum_{\ell=1}^{n} \sum_{e=1}^{N_{el}} t^e l_{\Gamma}^{e} 
		(\overline{\overline{\xi}}_{n,\ell}-\overline{\overline{\xi}}_{n,\ell-1}) t_{n,\ell}, \\ 
%\quad \text{for all} \quad \overline{\overline{\xi}}_{n,l}-\overline{\overline{\xi}}_{n,l-1} \ge 0
%	\end{equation}
%	
%	\begin{equation}
		\label{eq:DissFracModeII}
		W_{D, II} &= \sum_{\ell=1}^{n} \sum_{e=1}^{N_{el}} t^e l_{\Gamma}^{e} 
		(\overline{\overline{\xi}}_{m,\ell}-\overline{\overline{\xi}}_{m,\ell-1}) 
		\left|t_{m,\ell}\right| , 
%\quad \text{for all} \quad \overline{\overline{\xi}}_{m,l}-\overline{\overline{\xi}}_{m,l-1} \ge 0
	\end{align}
where $n$ is the number of completed time steps, $N_{el}$ is the number of finite elements, $t^e$ is the element thickness, $l_{\Gamma}^{e}$ is length of the embedded discontinuity in the finite element, $t_{n}$ is the tensile traction, $t_{m}$ is the shear traction, and $\overline{\overline{\xi}}_{n}$ and $\overline{\overline{\xi}}_{m}$ are the softening variables in the normal and tangent direction, respectively.

\subsection{Uniaxial cyclic test}
	\begin{figure}[!hbt]
		\centering
		\includegraphics[width=14cm]{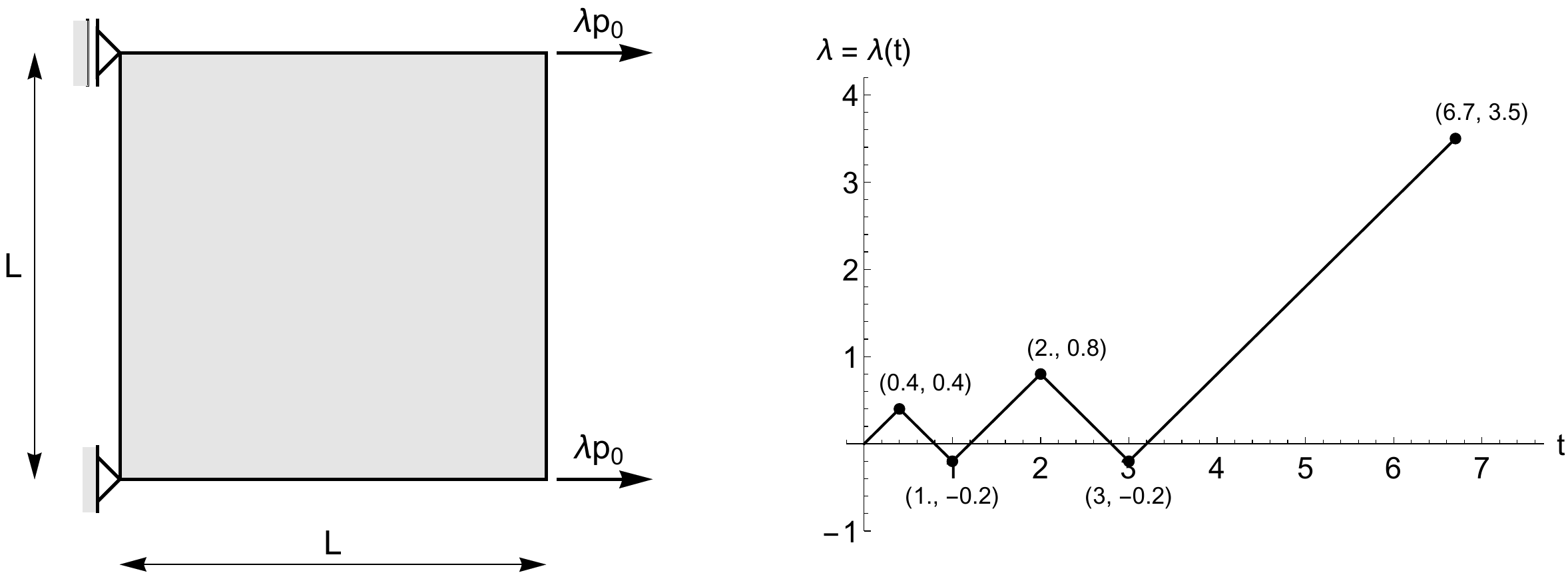}
		\caption{Uniaxial cyclic test. \text{Left}: Problem illustration. \text{Right}: 
		Loading protocol with the load factor $\lambda = \lambda \left( t \right)$ as a function 
		of the monotonically increasing quasi-time parameter $t \ge 0$.
		%  The loading protocol
		%is described by lines that connects the points $\left(t,\lambda\left(t\right)\right)$ . 
		}
		\label{fig:TensIgeometry}
	\end{figure}

Consider a square element of edge length $L = 200$ mm and thickness $10$ mm that is in a plane stress state. A static analysis is performed where the element is restrained on one side and subjected to an imposed horizontal displacement $p = \lambda p_{0}$ on the opposite side, see Fig.~\ref{fig:TensIgeometry} (left), where $p_{0} = 0.1$ mm and $\lambda$ is the load multiplier. The material parameters are Young's modulus $E=30000$ N/mm$^2$, Poisson's ratio $\nu=0.2$, the tensile strength $\sigma_{un} = 3$ N/mm$^2$, and shear the strength $\sigma_{um}=1$ N/mm$^2$. The mode I and mode II fracture energies are $G_{fn}=0.2$ N/mm and $G_{fm}=0.2$ N/mm, respectively. The loading procedure (see Fig.~\ref{fig:TensIgeometry} (right)) is composed of two tension-compression cycles that are completed by tension until the element separates into two parts. The element is in a homogeneous stress state. Therefore, the crack embedding criterion (see Subsection~\ref{sec:ED_Criterion}) is met in all four bulk integration points at the same time and the discontinuity line is positioned through the centre of the element.

	\begin{figure}[!hbt]
		\centering
		\includegraphics[width=15cm]{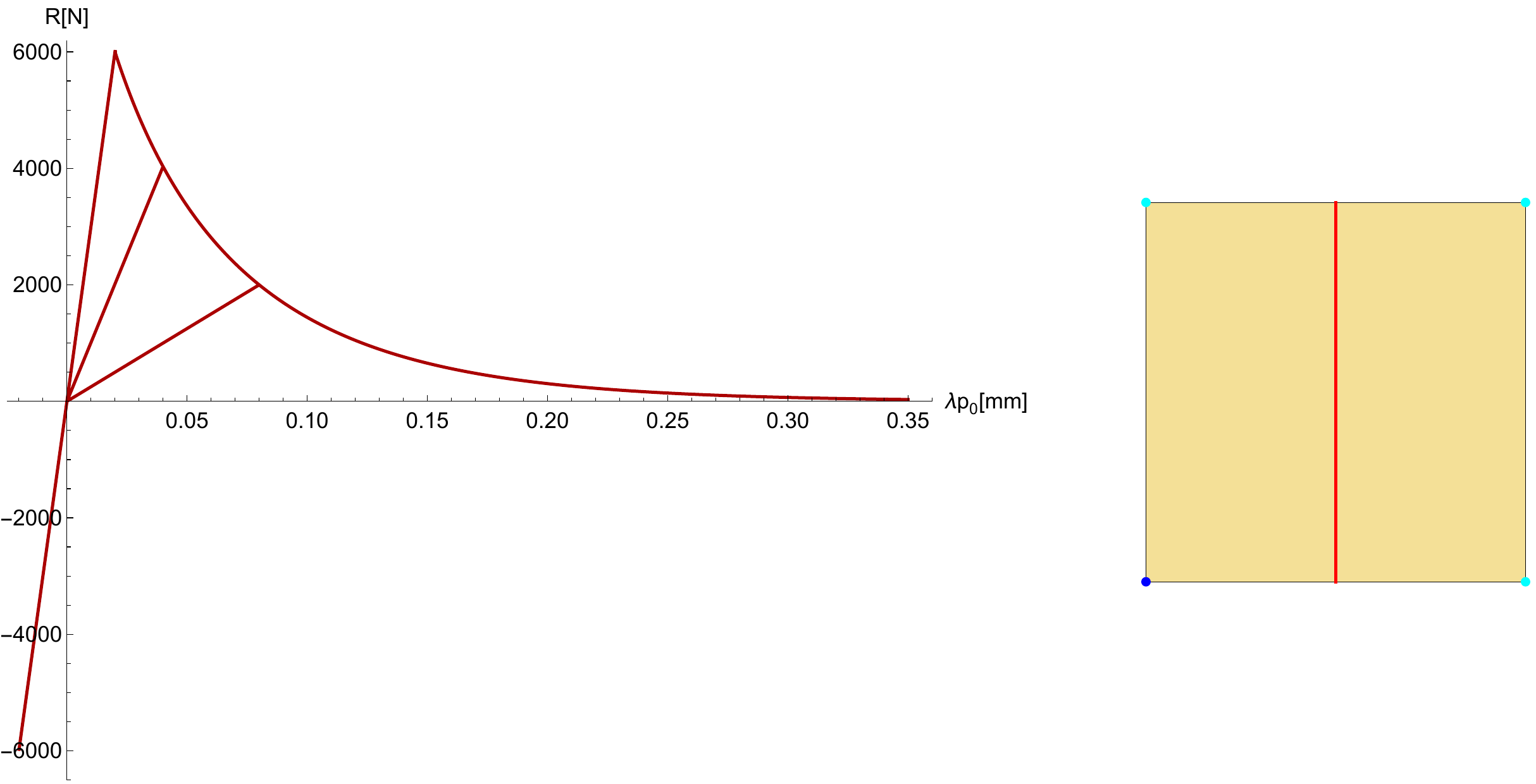}
		\caption{Uniaxial cyclic test. \text{Left}: Reaction versus imposed displacement. 
		\text{Right}: Position of embedded discontinuity.} 
		\label{fig:TensIresponse}
	\end{figure}

\begin{figure}[!hbt]
		\centering
		\includegraphics[width=16cm]{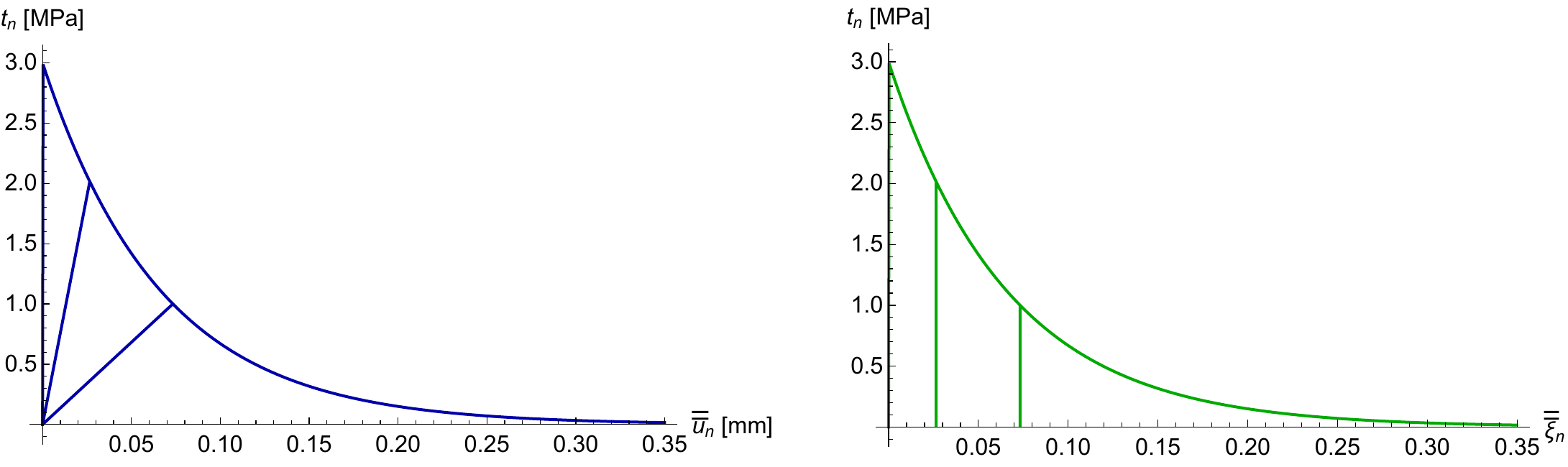}
		\caption{Uniaxial cyclic test. \text{Left}: Traction versus jump in displacement. 
			%$\overline{\overline{u}}_{n}$. 
			\text{Right}: Traction versus softening variable. %$\overline{\overline{\xi}}_{n}$.
		} 
		\label{fig:TensIinternal}
	\end{figure}

%	\noindent 
The reaction force versus imposed displacement is shown in the graph of Fig.~\ref{fig:TensIresponse}. Once the tensile strength $\sigma_{un}$ is exceeded, softening is triggered and the element resistance decreases with the crack opening in tension. In the unloading phase, the crack is closing and it stays closed in compression. Fig.~\ref{fig:TensIinternal} shows two graphs of internal variables at the discontinuity line: traction versus jump-in-displacement and traction versus softening variable $\overline{\overline{\xi}}_{n}$.

	\subsection{Tension test} 
	\label{sec:TensionTest}
\begin{figure}[!hbt]
		\centering
		\includegraphics[width=10cm]{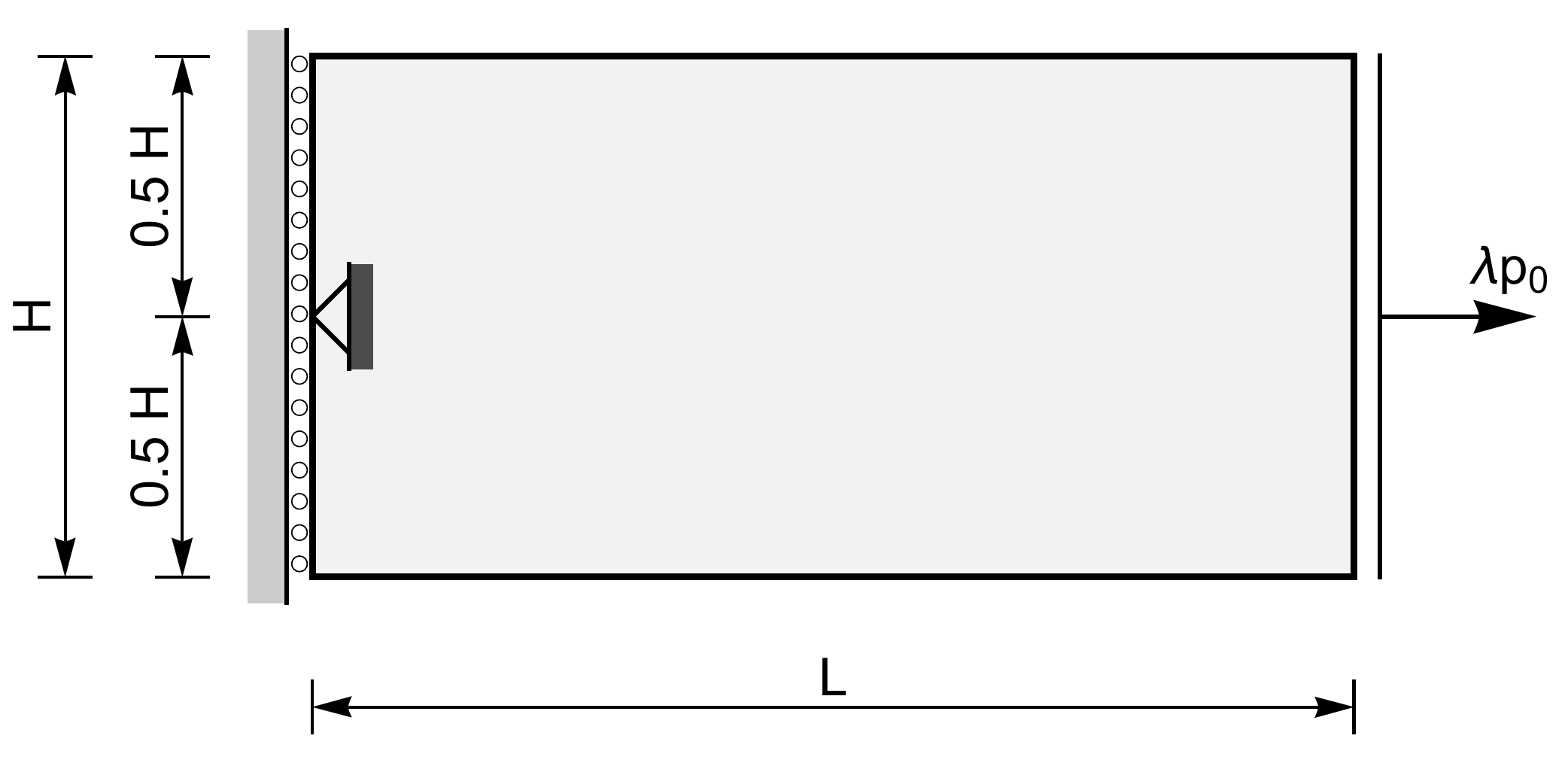}
		\caption{Tension test: Problem illustration.} 
		\label{fig:TensIIgeometry}
	\end{figure}

In order to assess the influence of the nucleation criterion from Subsection~\ref{sec:ED_Criterion} on the global cracking mechanism, a tension test is performed. Consider a plane-stress solid of rectangular shape with length $L = 400$ mm, height $H = 200$ mm, and thickness $t^e=10$ mm, see Fig.~\ref{fig:TensIIgeometry}. On the supported edge of the solid, the horizontal displacement is prevented, and the mid-point is fixed. The opposite edge is subjected to an imposed horizontal displacement $p = \lambda p_{0}$, where $p_{0}=0.1$ mm and $\lambda$ is a load multiplier. The material parameters are Young's modulus $E=30000$ N/mm$^2$, Poisson's ratio $\nu=0.2$, the tensile strength $\sigma_{un} = 3$ N/mm$^2$, the shear strength $\sigma_{um}=1$ N/mm$^2$, and fracture energies for mode I and mode II given by $G_{fn}=0.2$ N/mm and $G_{fm}=0.2$ N/mm, respectively. 

\begin{figure}[H]
		\centering
		\includegraphics[width=16cm]{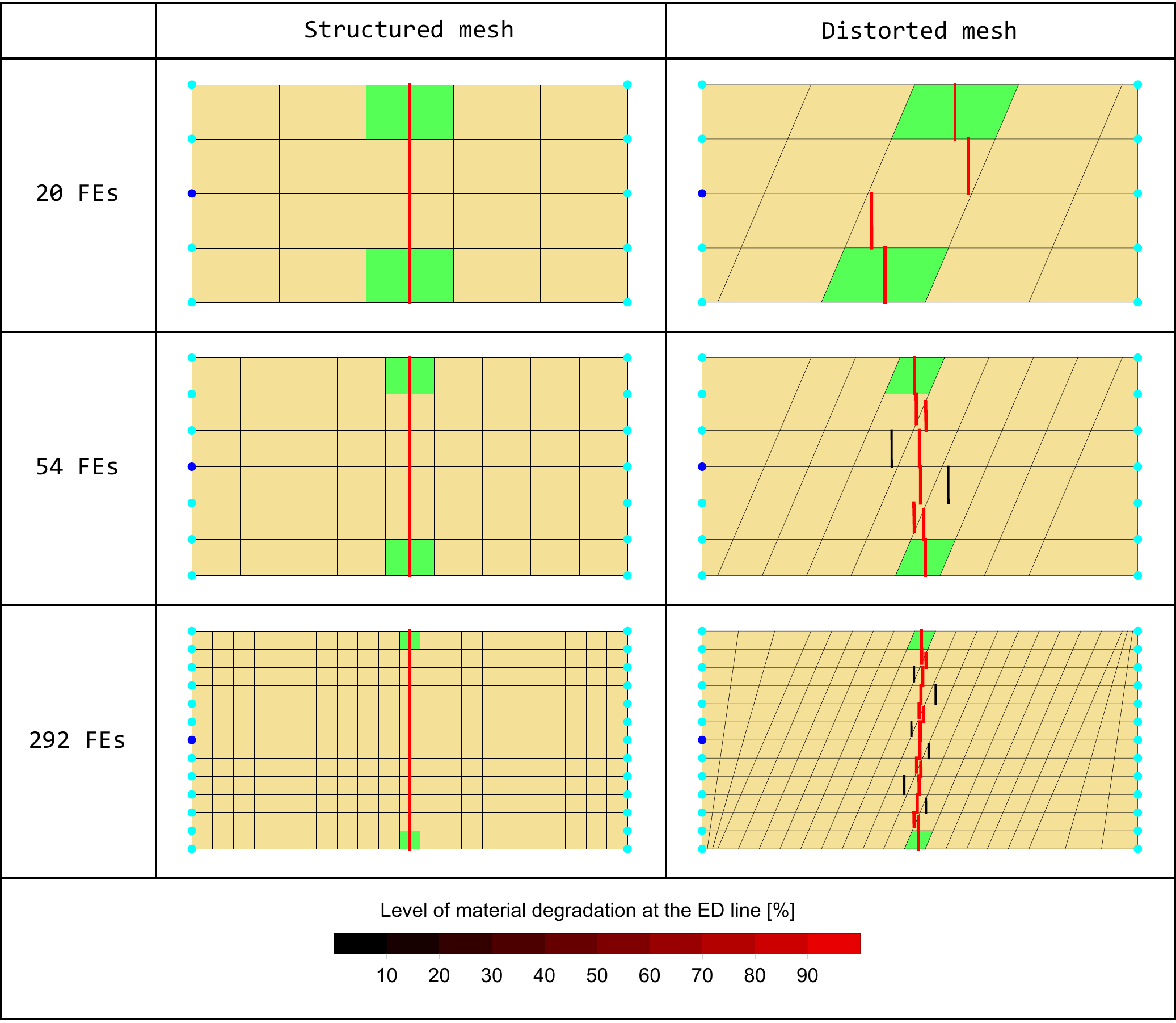}
		\caption{Tension test: Meshes with nucleated cracks, green elements have weakened strength.} 
		\label{fig:TensIIcracksAll}
	\end{figure}

\begin{figure}[H]
		\centering
		\includegraphics[width=15cm]{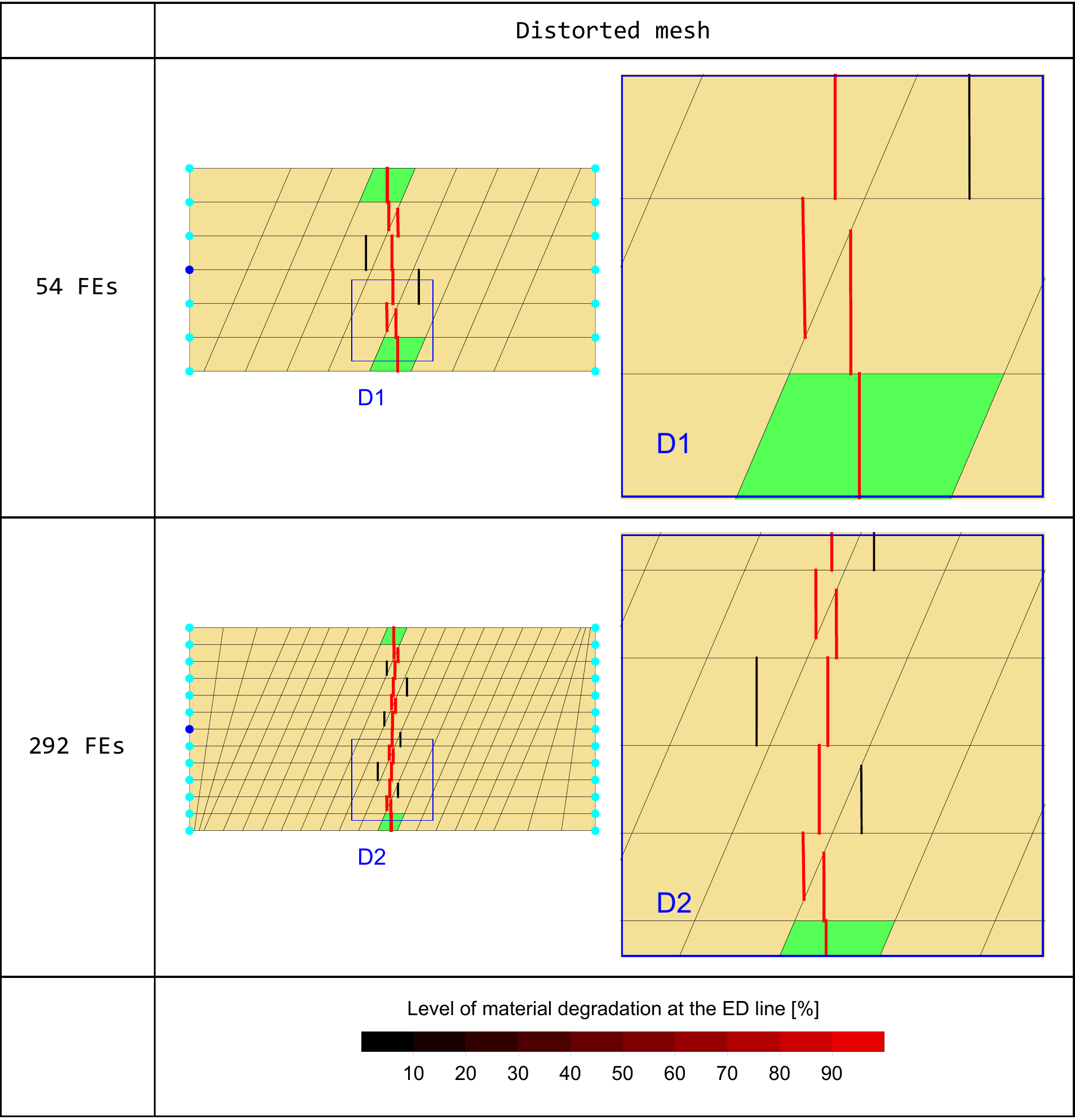}
		\caption{Tension test: Details of distorted meshes with nucleated cracks.} 
		\label{fig:TensIIcracksPart}
	\end{figure}

\begin{figure}[H]
		\centering
		\includegraphics[width=16cm]{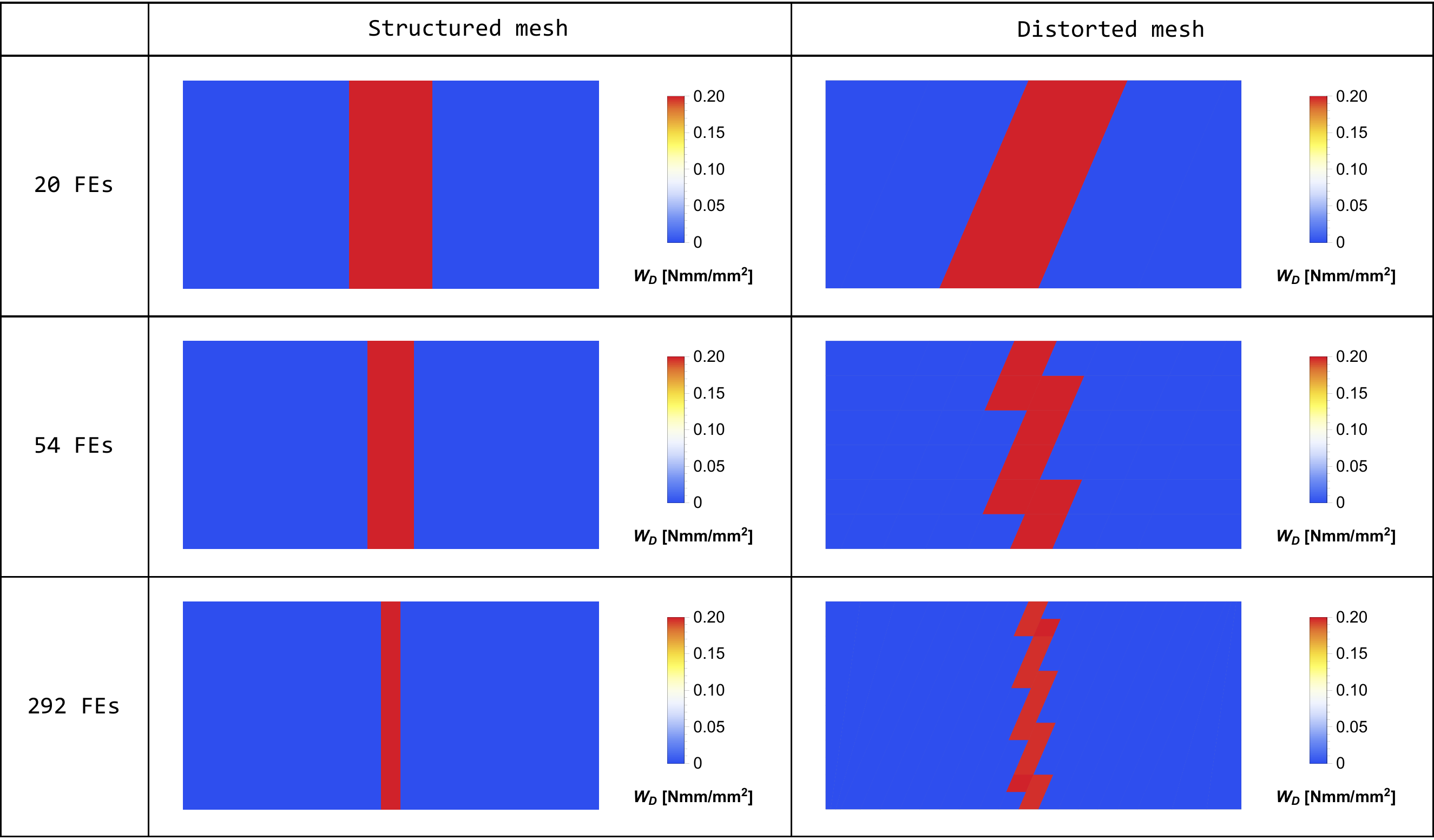}
		\caption{Tension test: Element-wise presentation of dissipated fracture energy at $\lambda=0.05$.} 
		\label{fig:TensIIdissAll2}
	\end{figure}

\begin{figure}[H]
		\centering
		\includegraphics[width=\linewidth]{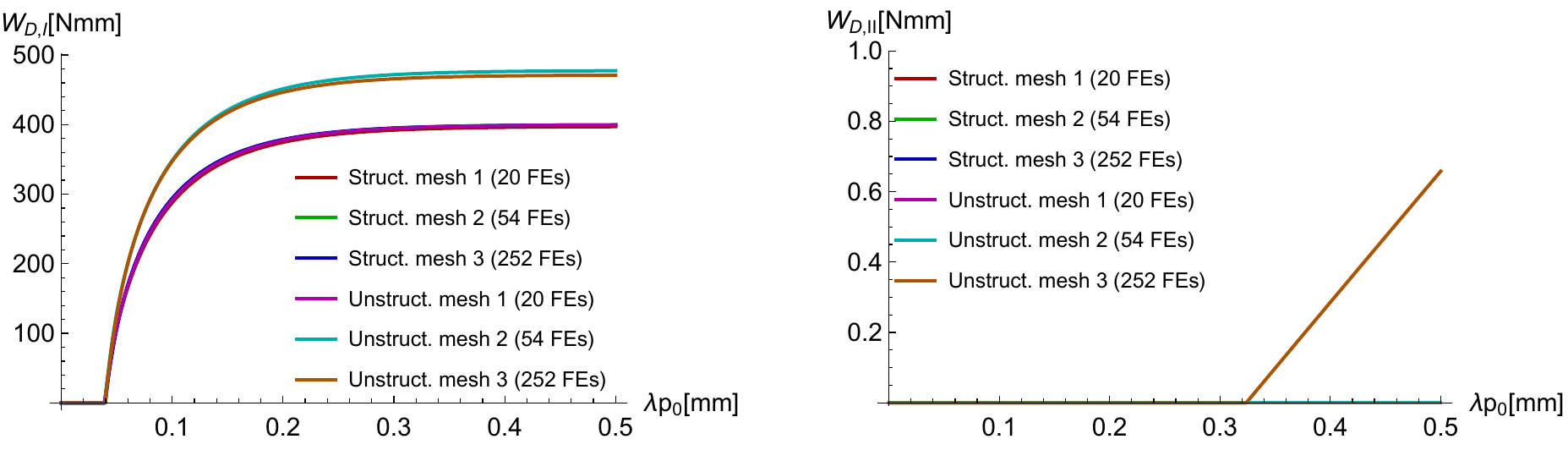}
		\caption{Tension test: Dissipated fracture energy versus imposed displacement for mode I (left) and mode II (right).} 
		\label{fig:TensIIGraphDissp}
	\end{figure}

\begin{figure}[H]
		\centering
		\includegraphics[width=\linewidth]{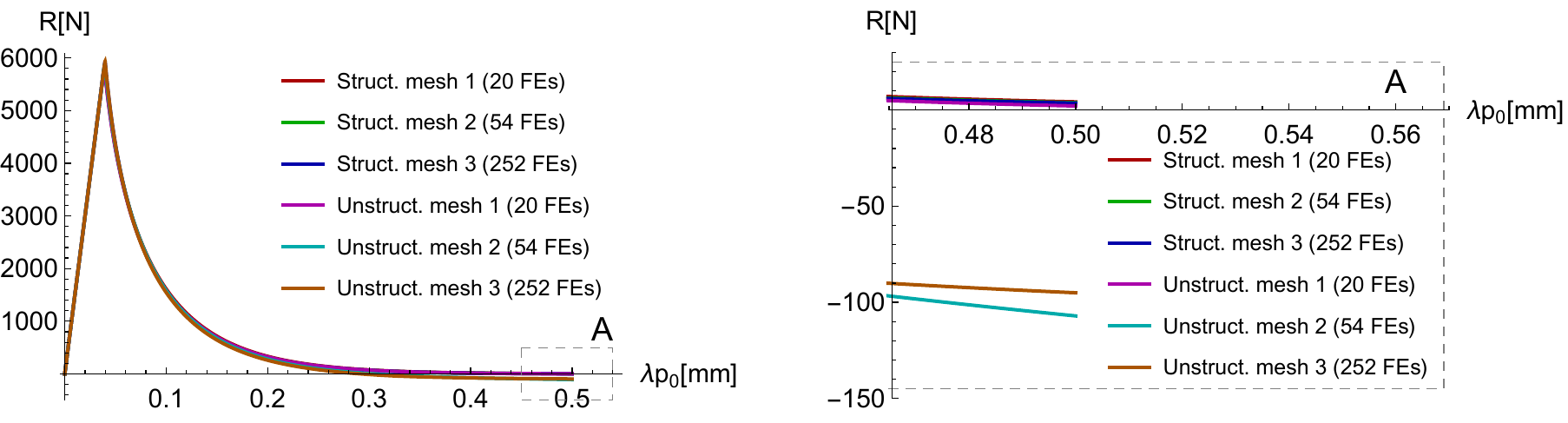}
		\caption{Tension test: Reaction versus imposed displacement. Right: Detail A.} 
		\label{fig:TensIIGraphRp}
	\end{figure}

For the tension test three structured and three distorted meshes are considered. In order to initiate a specific crack location and avoid computational bifurcation that stems from the initially  homogeneous stress state due to the uniformly imposed displacement along the width, the tensile strength is lowered to $\sigma_{un,weak}=0.95 \sigma_{un}$ for two elements, positioned on the upper and lower side of the solid at its half length. Fig.~\ref{fig:TensIIcracksAll} shows the meshes, with the elements having weaker strength in green, and it also displays the state of the cracks at $\lambda=5$ for each mesh. The cracks are coloured according to the level of softening, i.e.\ from black ($0\%$) to red ($100\%$), where $100\%$ means fully softened crack without any cohesion left. Due to the uniaxial stress state, all cracks are vertical. From Fig.~\ref{fig:TensIIcracksAll} it may be seen that the cracks in the structured meshes form an one-line crack pattern, whereas this is not the case for the distorted meshes, although the red cracks tend to converge towards a one-line crack pattern with increasing density of the mesh --- recall that the black coloured cracks are not active and do not dissipate fracture energy. Fig.~\ref{fig:TensIIcracksPart} displays zoomed portions of the distorted meshes, showing closely the crack patterns. From Fig.~\ref{fig:TensIIcracksPart} it may be observed that the total length of the traction-less cracks $L_{crack}$ is larger for the distorted than for the structured meshes, where the crack-length is $L_{crack}=H=200$ mm. For the distorted mesh with 54 elements, the crack-length amounts to $L_{crack} = 240.5$ mm, and for the distorted mesh with 292 elements the crack-length is $L_{crack} = 235.8$ mm. Fig.~\ref{fig:TensIIdissAll2} illustrates, in an element-wise manner, the dissipated fracture energy per unit crack surface (in Nmm/mm$^{2}$). The curves in Fig.~\ref{fig:TensIIGraphDissp} (left) show that a larger $L_{crack}$ leads to larger dissipated fracture energy for mode I, because the dissipated fracture energy depends on the crack length, see \eqref{eq:DissFracModeI}. The final value for the distorted meshes with 54 and 292 elements is approximately $20.3\%$ higher than that for the other meshes. The curves in Fig.~\ref{fig:TensIIGraphDissp} (right) show that the distorted mesh with 292 elements dissipate also a negligible amount of mode II fracture energy. The meshes influence also the reaction-displacement curves shown in Fig.~\ref{fig:TensIIGraphRp}, but only slightly, when the full collapse of the specimen happens. When the total horizontal reaction force $R$ is around zero, the response curves for the distorted meshes with 54 and 292 elements slightly deviate from the curves of structured finite element meshes ($R$ even gets small negative values). This example shows that the crack pattern can be expected to converge to the ``correct one'' when the distorted mesh is reasonably dense, that the mesh distortion has a small influence on the load-displacement curve, and that the mesh distortion influences the amount of the dissipated fracture energy. These are the insights to keep in mind when interpreting the results of more complex examples.

\subsection{Thermal shock test} \label{sec:therm-shock}
The next problem to be analysed is the thermal shock test from \cite{Farrell et al 2017-art}, where the authors used the phase-field method for the computation of the crack pattern. Originally, this kind of test was performed as an experiment on brittle ceramics in order to investigate crack patterns with periodic array of cracks (see \cite{Bahr et al 1986-art, Geyer Nemat 1982-art, Shao et al 2011-art}). In the experiment, small rectangular solids of elastic homogeneous material were stacked together into a thick slab, which was uniformly heated to a temperature of $T_{0}$, and then thrown into a cold bath causing a temperature drop of $\upDelta T$ on the lateral surfaces. The sudden change of temperature induces an in-homogeneous stress field that causes the nucleation of many cracks forming a quasi periodic pattern.
	
	\begin{figure}[!hbt]
		\centering
		\includegraphics[height=5.5cm]{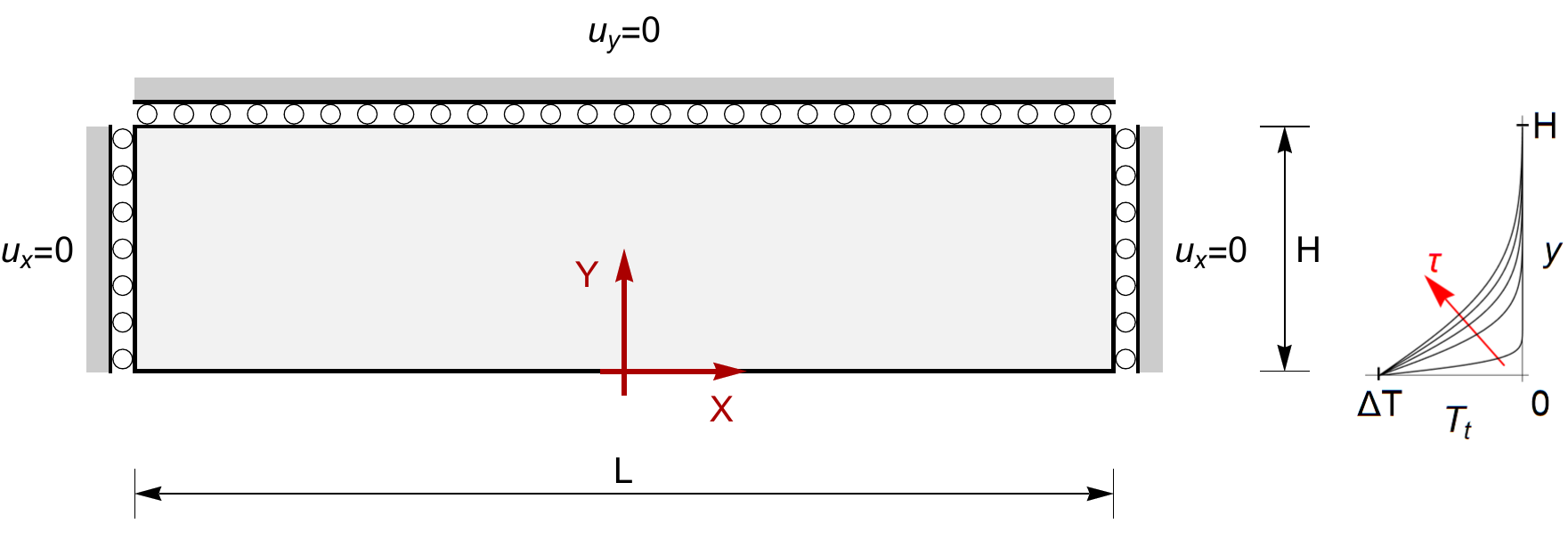}
		\caption{Thermal shock. Left: Geometry and boundary conditions.
			% of numerical model. 
			Right: Error complementary functions (Erfc) distribution along the 
			%model
			height.} 
		\label{fig:ThermalShockGeometry}
	\end{figure}
	
%	\noindent 
As in \cite{Farrell et al 2017-art}, we prepare a simplified finite element model of the experimental test in dimensionless form. The considered portion of the brittle slab (hereinafter called specimen) has length $L = 40$, height $H = 10$ and thickness $t^e=1$ in plane stress condition. Fig.~\ref{fig:ThermalShockGeometry} shows the boundary conditions applied on the top and side surfaces, where the displacement perpendicular to the edge is prevented, while the free bottom side (at $y=0$) is exposed to a temperature change. The material properties are: Young's modulus $E = 1$, Poisson's ratio $\nu = 0.3$, tensile strength $\sigma_{un} = 0.61$, shear strength $\sigma_{um} = 0.61$, and the fracture energies for modes I and II are $G_{fn} = 1$ and $G_{fm} = 1$. The thermal shock effect is modelled via a prescribed time-dependent thermal strain (in vector form):
	\begin{equation}
		\label{eq:ThermalShock_TempStrain}
		\boldsymbol{\epsilon}_{t} \left(\tau, y\right) = \beta_{th} \, 
		T_{t}\left(\tau, y\right) \, \boldsymbol{b}, \quad 
		\boldsymbol{b} = {\left[1,1,0\right]}^T ,
	\end{equation}
where $\beta_{th} = 1$ is the thermal expansion coefficient, $T_{t}\left(\tau, y \right)$ is the temperature change field in the specimen, and $\tau$ is the time. The prescribed strains $\boldsymbol{\epsilon}_{t}$ are subtracted from the total strains to obtain the effective strains $\boldsymbol{\epsilon}_{\text{eff}}$:
	\begin{equation}
		\label{eq:ThermalShock_EffStrain}
  \boldsymbol{\epsilon}_{\text{eff}} = \boldsymbol{\epsilon}  
		 \left(\boldsymbol{\xi},{\Gamma}^{e} \right) - \boldsymbol{\epsilon}_{t} \left(\tau, y\right),
	\end{equation}
which are used in the equilibrium equations (see \eqref{eq:GlobalEq} and 
\eqref{eq:LocalEqh1}--\eqref{eq:LocalEq3}) instead of $\boldsymbol{\epsilon}$ from \eqref{eq:strains}. The temperature change field $T_{t}$ is (see Fig.~\ref{fig:ThermalShockGeometry})
	\begin{equation}
		\label{eq:ThermalShock_TempField}
		T_{t}\left(\tau, y \right) = -\upDelta T \, 
		\text{erfc}\left(\frac{y}{2 \sqrt{k_c \tau}} \right), \quad  
		\upDelta T = a \, \upDelta T_{c} ,
	\end{equation}
where $\text{erfc}$ is the complementary error function, $k_{c} = 1$ is the thermal conductivity, $\upDelta T$ is the applied temperature drop, $\upDelta T_{c} = \sqrt{\frac{3}{8}}$ is the critical temperature drop, and $a$ is an intensity factor. When $a < 1$ ($\upDelta T < \upDelta T_{c}$), the solution is purely elastic with no cracks, and for $a > 1$ ($\upDelta T > \upDelta T_{c}$) the cracks appear. According to \cite{Farrell et al 2017-art}, $T_{t}$ is the analytical solution of an approximate thermal diffusion problem with Dirichlet boundary conditions on the temperature for a semi-infinite homogeneous slab of thermal conductivity $k_{c}$ (the influence of the cracks on the thermal conductivity is neglected).

	\begin{figure}[!hbt]
		\centering
		\includegraphics[width=17cm]{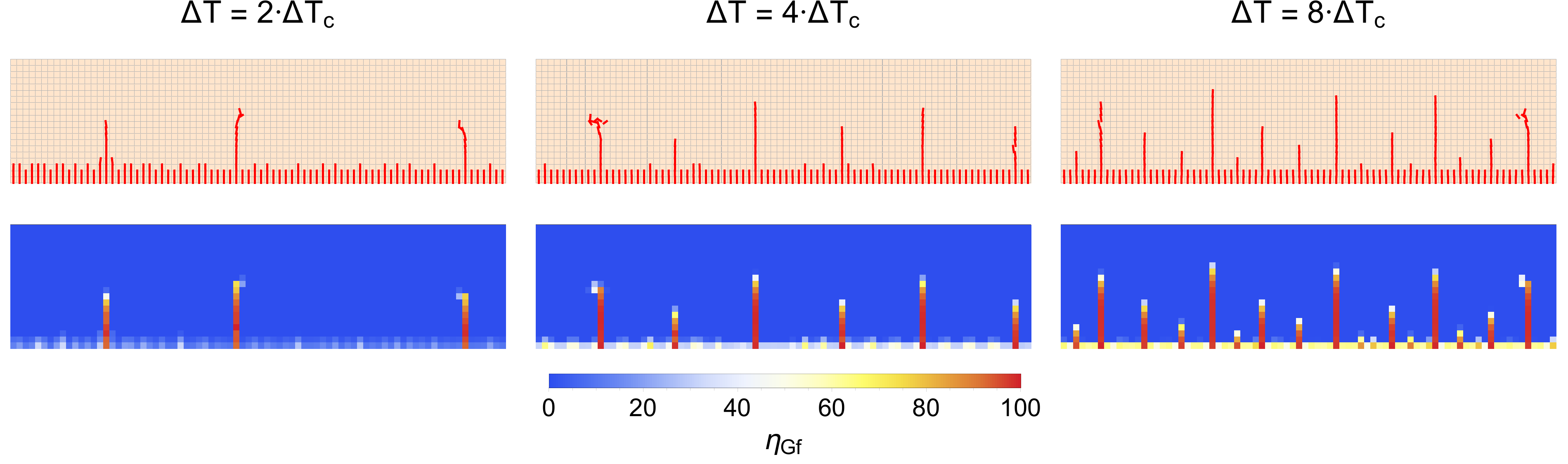}
		\caption{Thermal shock: Nucleated cracks (first row) and specific dissipated energy (second row) 
		at $\tau=5$ for three different intensities of the temperature drop (2, 4 and 8). 
		The dissipated fracture energy is presented as a percentage of $G_{fn}$, denoted 
		as $\eta_{Gf}$. The mesh size is $80 \times 20$.} 
		\label{fig:ThermalShockCracks}
	\end{figure}

As the temperature change field $T_{t}$ is constant along the horizontal $x$-direction, the bottom side is exposed to uniform strains, which induces homogeneous tension along the $x$-direction. Again, in order to initiate crack nucleation and avoid computational bifurcation, we introduce a small material heterogeneity in the tensile strength $\sigma_{un}$ for the finite elements that lie along the bottom edge by adding a perturbation $\upDelta \sigma_{un}$ as  
	\begin{equation}
		\label{eq:ThermalShock_Pert}
		\sigma_{un}^{pert}=\sigma_{un}+\upDelta \, \sigma_{un} ,
	\end{equation}
where $\upDelta \, \sigma_{un}$ is uniform white noise, i.e.\ randomly and independently generated for each element from the uniform distribution with bounds $\pm 0.02 \, \sigma_{un}$. 

	\begin{figure}[!hbt]
		\centering
		\includegraphics[width=10.5cm]{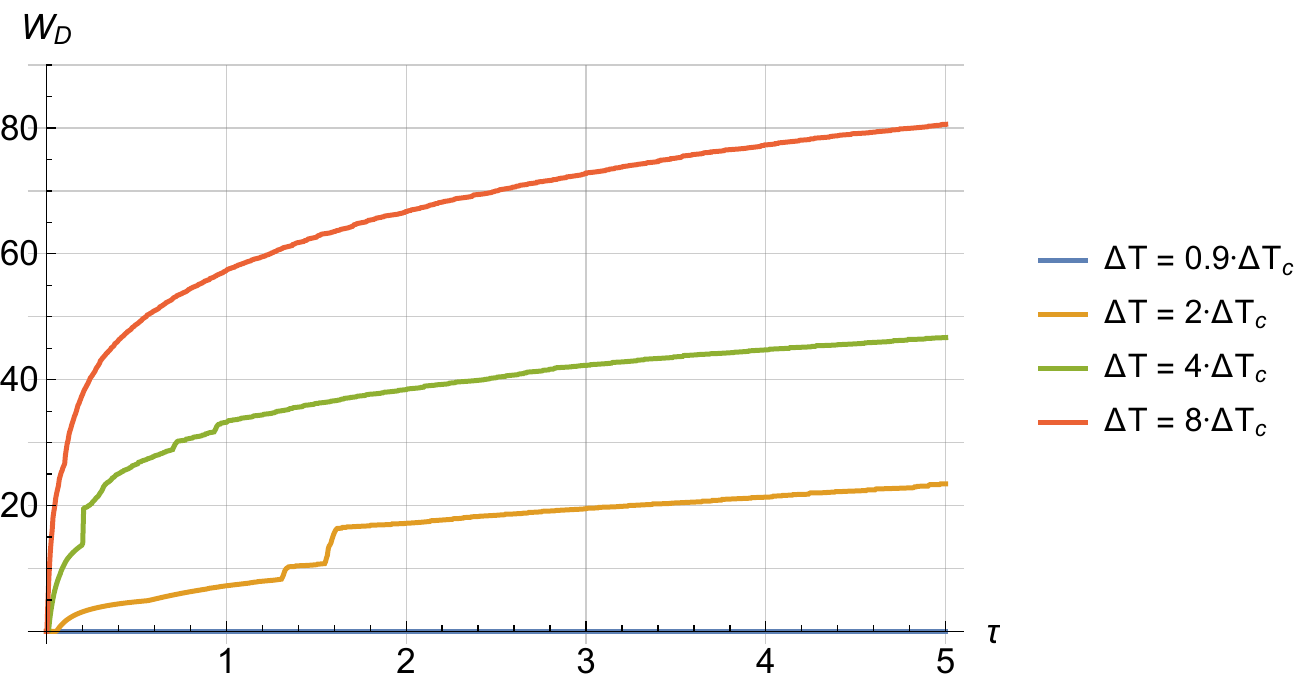}
		\caption{Thermal shock: Dissipated fracture energy for different intensities of
		the temperature drop ($a=$0.9, 2, 4 and 8).} 
		\label{fig:ThermalShockDiss}
	\end{figure}

%	\noindent 
Simulations were performed for four values of the intensity factor, i.e.\ for $a \in \{0.9, 2, 4, 8 \}$. The results are presented in Figs.~\ref{fig:ThermalShockCracks} and~\ref{fig:ThermalShockDiss}. The cracks open in mode I, as presented in Fig.~\ref{fig:ThermalShockCracks}, which shows finite elements with cracks and dissipated fracture energy. One can observe essentially periodic crack patterns with significant distances between the two larger cracks of the same length. For the higher values of the intensity factors $a$, there is a larger number of longer cracks with shorter cracks appearing in-between. The evolution of the total dissipated fracture energy is shown in Fig.~\ref{fig:ThermalShockDiss} (note that for $a=0.9$, the specimen stays elastic with no cracks). The shapes of the curves are similar to those from \cite{Farrell et al 2017-art}: the curve is smooth for $a=8$, while for $a=2$ the curve has a jump almost at the same time as in \cite{Farrell et al 2017-art}.

	\begin{figure}[!hbt]
		\centering
		\includegraphics[width=12.5cm]{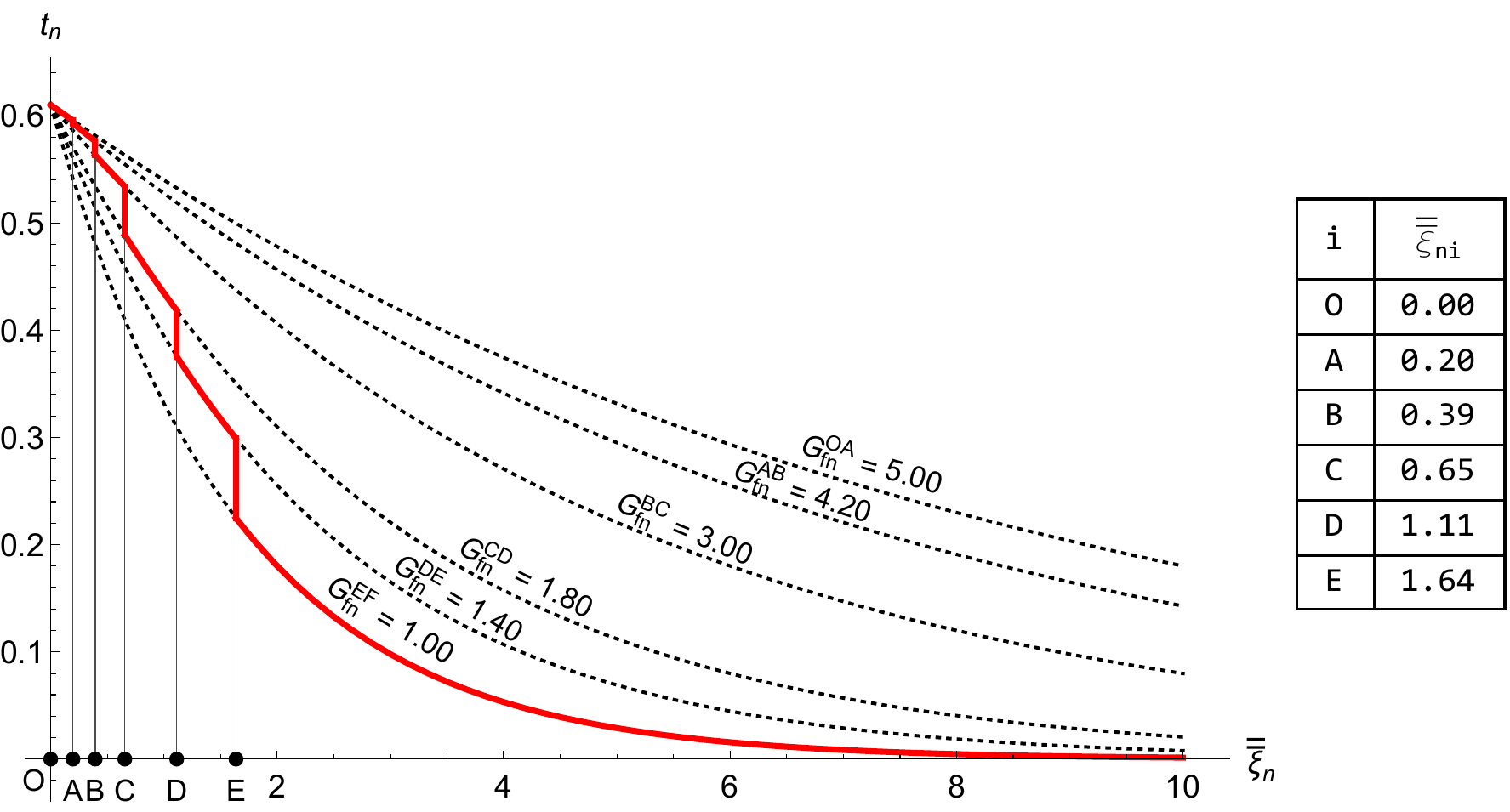}
		\caption{Thermal shock: The ``zig-zag'' softening function.} 
		\label{fig:ThermalShockGradualFunction}
	\end{figure}

%	\noindent 
As described in Section~\ref{sec:EqSys}, we solve the global and local equilibrium equations simultaneously. This way of solving the augmented system of equations (see \eqref{eq:AugmentedSystem2}) is sensitive to low material fracture energies $G_{fn}$ and $G_{fm}$, as is the case of the thermal shock test example. Namely, when softening is triggered, the initial tangents to the exponential softening function are very steep, causing corresponding large terms in the augmented stiffness matrix (see \eqref{eq:AugmentedSystem2Eq}) %are large (in absolute values) 
in comparison to other terms. Consequently, the convergence may be problematic. We regulated this problem with a ``zig-zag'' softening function illustrated in Fig.~\ref{fig:ThermalShockGradualFunction}, which is used instead of the smooth softening function from Fig.~\ref{fig:ExpSoft}. The initial part of the exponential function in the interval between $\xi_{O}$ and $\xi_{E}$ was modified to drop in five jumps from the softening function with $G_{fn}^{OA} = 5.00$ to the desired function with $G_{fn}^{EF}=G_{fn}=1.00$. The same zig-zag function is also used for softening in mode II. Note that the choice of initial fracture energy $G_{fn}^{OA}$ is based on the results of thermal shock simulations for different input values for $G_{fn}$ and $G_{fm}$. For a higher fracture energy $G_{f}$, the convergence was good, but the crack patterns as presented in Fig.~\ref{fig:ThermalShockCracks} did not appear. The conclusion is that for the thermal shock test example the fracture energy determines the velocity of crack propagation, and it plays a crucial role in the numerical simulations in terms of robustness.

\subsection{Four point bending test} 
\label{sec:FourPB}

	\begin{figure}[!htb]
		\centering
		\includegraphics[width=10cm]{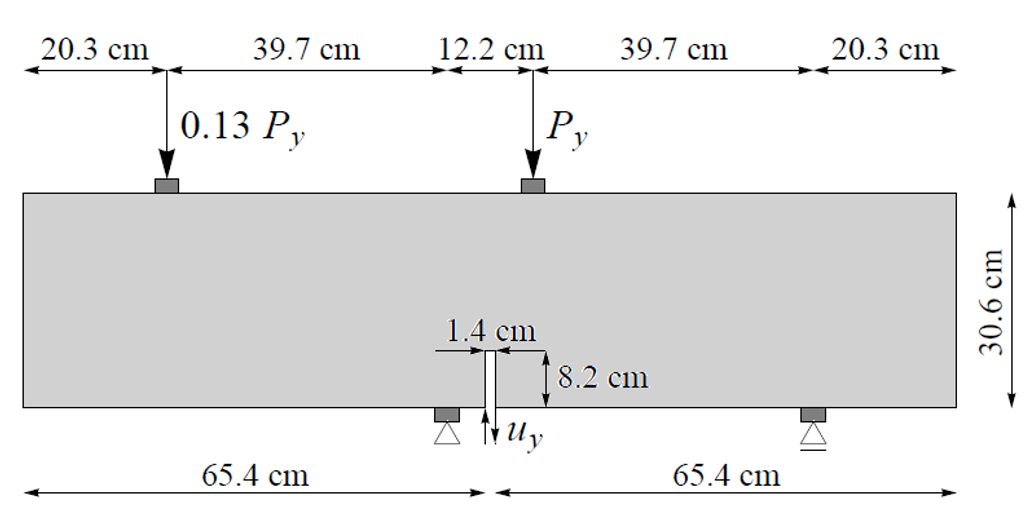}
		\caption{Four point bending test: Geometry, boundary and loading conditions.} 
		\label{fig:FourPBgeometry}
	\end{figure}

\begin{figure}[!htb]
		\centering
		\includegraphics[width=10cm]{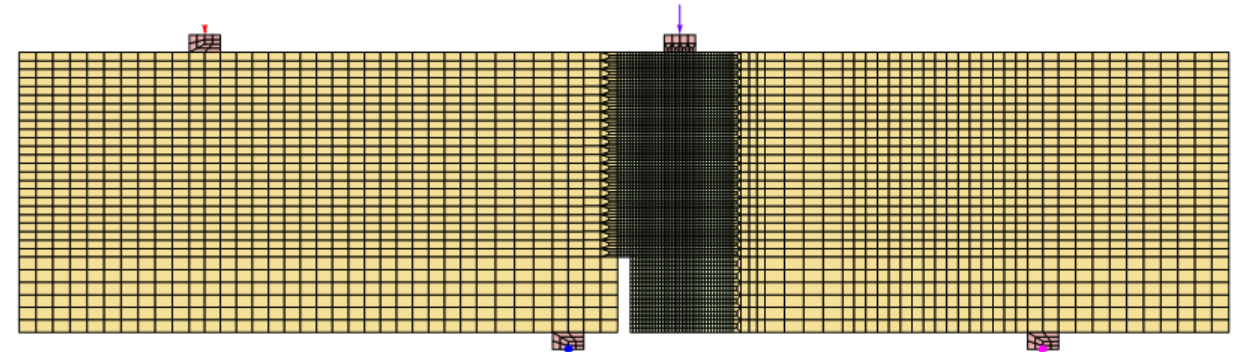}
		\caption{Four point bending test: Finite element mesh.} 
		\label{fig:FourPBmesh}
	\end{figure}

\begin{figure}[!htb]
		\centering
		\includegraphics[width=11cm]{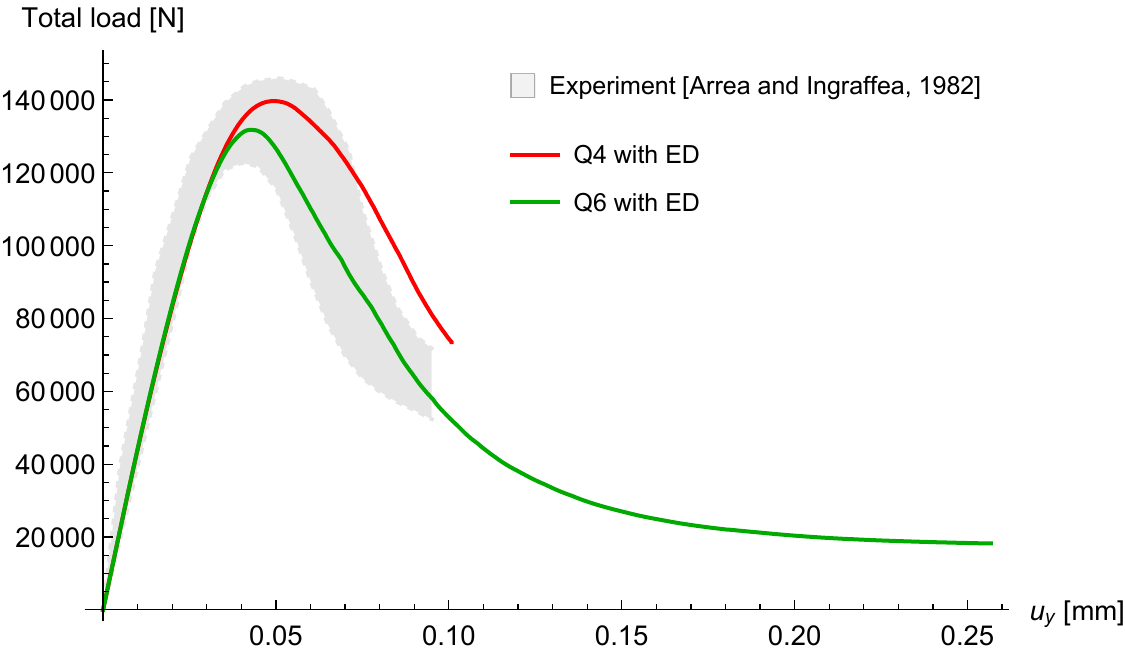}
		\caption{Four point bending test: Total load versus relative displacement at the notch mouth.} 
		\label{fig:FourPBResponse}
	\end{figure}

Another challenging example is the four point bending test from \cite{Arrea Ingraffea 1982-art} with geometry as well as boundary- and loading-conditions illustrated in Fig.~\ref{fig:FourPBgeometry}. It is a concrete beam with a predefined notch at its half-length subjected to two vertical forces $0.13P_{y}$ and $P_{y}$, where $P_{y}=\lambda P_{0}$, $P_{0}=1$~N, and $\lambda$ is the load multiplier. The considered plane stress model has a  thickness of 156 mm. The material properties for
the concrete are: Young's modulus $E=32000$ N/mm$^2$, Poisson's ratio $\nu=0.2$, tensile strength $\sigma_{un}=2.8$ N/mm$^2$, and shear strength $\sigma_{um}=1$ N/mm$^2$. The mode I and mode II fracture energies are $G_{fn}=0.09$ N/mm and $G_{fm}=0.09$ N/mm. The stiff caps (coloured dark gray in Fig.~\ref{fig:FourPBgeometry}) are modelled as elastic with Young's modulus $E=3.2\cdot10^{9}$ N/mm$^2$ and a Poisson's ratio of $\nu=0.2$. Fig.~\ref{fig:FourPBmesh} shows the finite element mesh refined in the region where crack path is expected to propagate. For the computation of the example we used a path-following method that controls a single constantly increasing/decreasing degree of freedom \cite{Stanic-Brank-Korelc 2016-art}, chosen here to be the vertical displacement at the notch mouth. Fig.~\ref{fig:FourPBResponse} shows the total load versus relative displacement at the notch mouth, which is in good agreement with the range of experimental data.

%	\noindent 
	\begin{figure}[!htb]
		\centering
		\includegraphics[width=10cm]{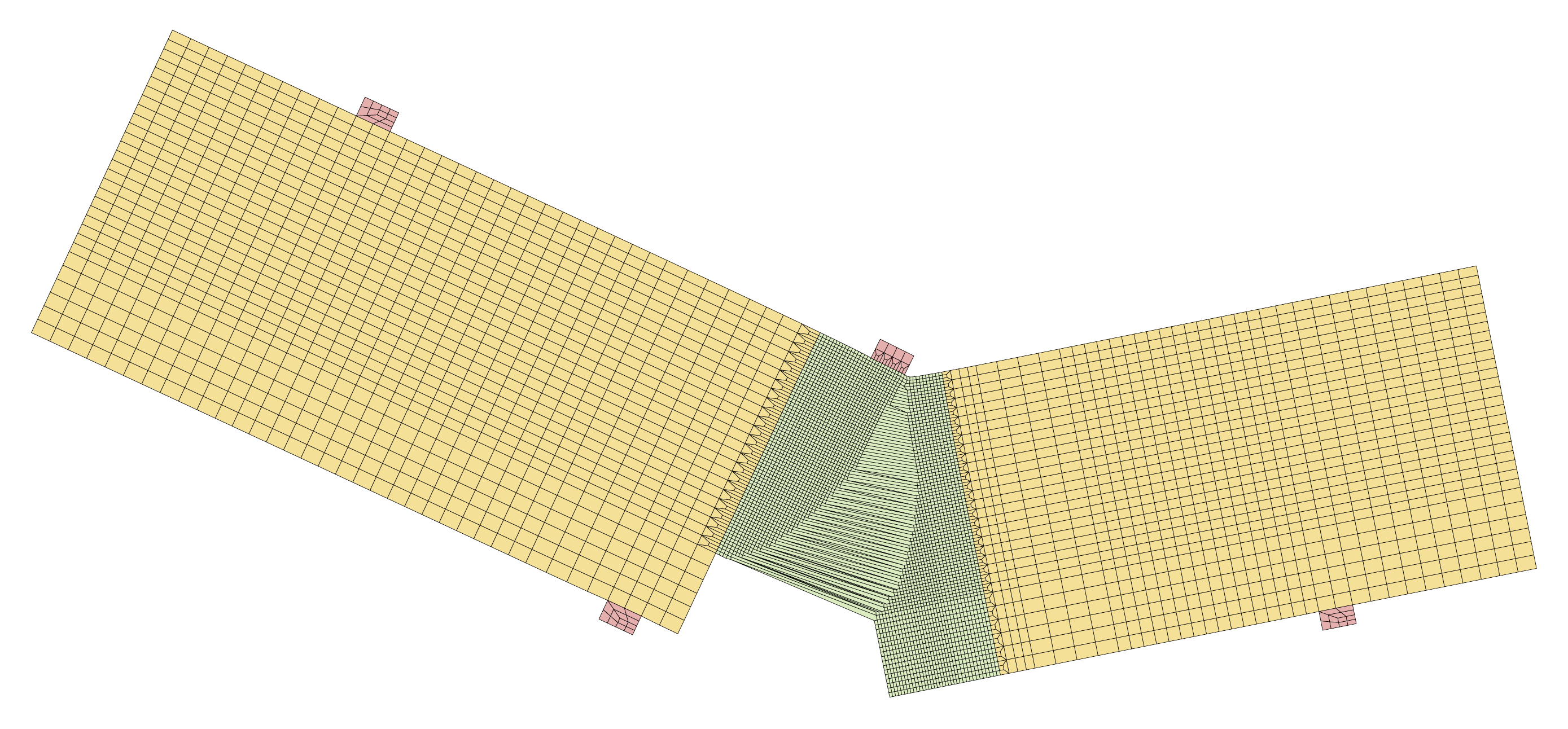}
		\caption{Four point bending test: Deformed finite element mesh 
		(displacements magnified by a factor of 250).} 
		\label{fig:FourPBDeformedConf}
	\end{figure}
	\begin{figure}[!htb]
		\centering
		\includegraphics[width=17cm]{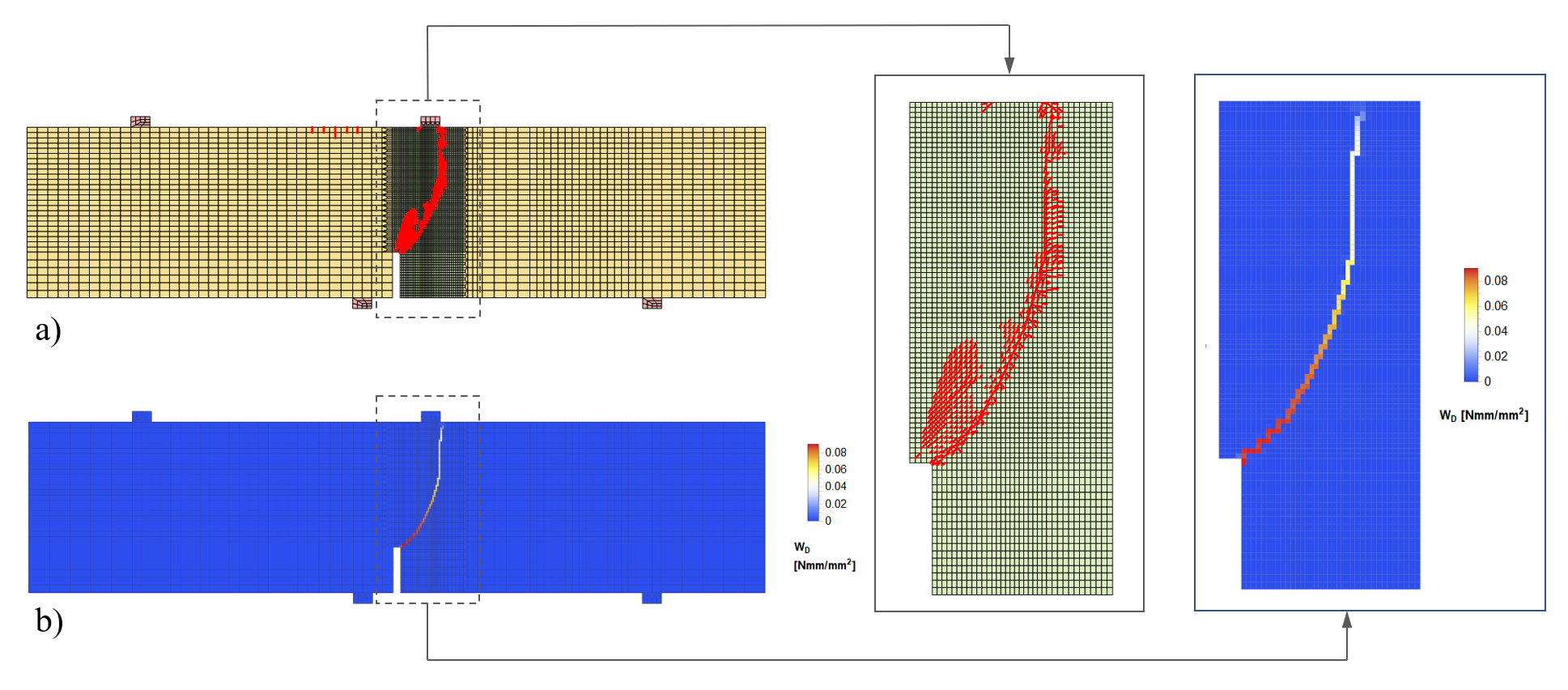}
		\caption{Four point bending test: a) Nucleated cracks in the mesh and 
		b) specific dissipated fracture energy at the end of the analysis. 
		% (at $\lambda=16114$).
		}
		\label{fig:FourPBcracks}
	\end{figure}

	\begin{figure}[!htb]
	\centering
	\includegraphics[width=17cm]{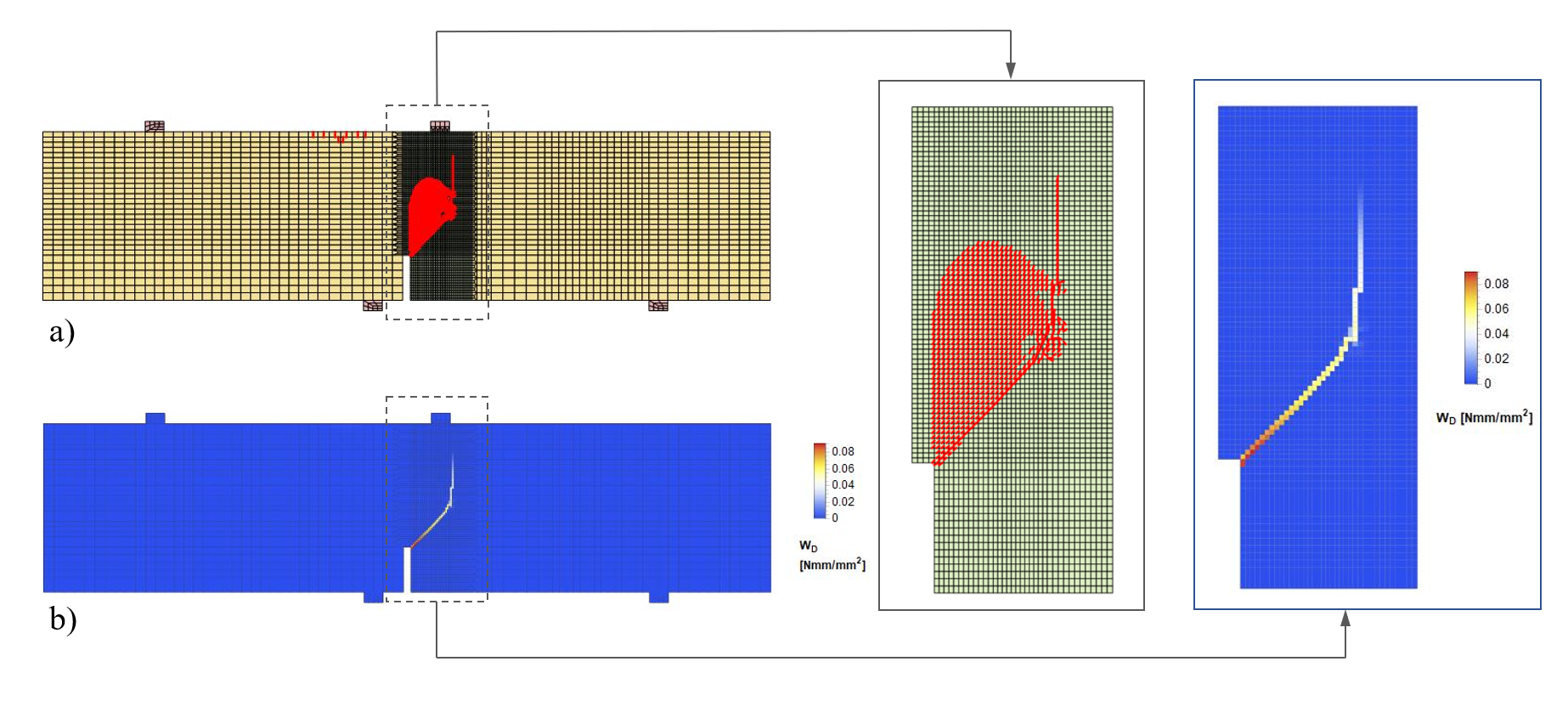}
	\caption{\correction{Four point bending test calculated without incompatible modes (Q4 with ED): a)~Nucleated cracks in the mesh and b)~specific dissipated fracture energy at the end of the analysis.}.
	}
	\label{fig:FourPB_Q4_cracks}
	\end{figure}

The criterion for crack nucleation and hence the embedded discontinuity activation (see Subsection~\ref{sec:ED_Criterion}) is met in many elements of the mesh (see Fig.~\ref{fig:FourPBcracks}a)), especially in the region around the notch. However, the fracture energy dissipates only in some of them (see Fig.~\ref{fig:FourPBcracks}b)), namely those that form the crack path line. The deformed finite element mesh in Fig.~\ref{fig:FourPBDeformedConf} confirms that the global crack path runs as indicated by the dissipated fracture energy in Fig.~\ref{fig:FourPBcracks}b). 

\correction{To illustrate the influence of the incompatible modes, the four point bending test is recomputed by suppressing the incompatible modes (and thus using the Q4ED element). The results are presented in Figs. \ref{fig:FourPBResponse} and \ref{fig:FourPB_Q4_cracks}. Fig. \ref{fig:FourPBResponse} compares the load-displacement curves: Q4ED computes higher limit load and stops earlier due to convergence loss. The discrepancy in limit loads stems from different crack paths. Fig. \ref{fig:FourPB_Q4_cracks} presents nucleated cracks in the mesh and the dissipated energy at the end of the Q4ED analysis. Fig. \ref{fig:FourPB_Q4_cracks}b shows that the shape of the Q4ED crack path is almost bi-linear. In contrast, Q6ED crack path bends gradually, see Fig. \ref{fig:FourPBcracks}b, which is an important detail. We can conclude that the incompatible modes considerably improve the simulation of crack path propagation.}

\subsection{Double edge notched specimen (DENS) test} 
	\begin{figure}[!htb]
	  \centering
		\begin{subfigure}{.49\textwidth}
			\centering
			% include first image
			\includegraphics[width=8.0 cm]{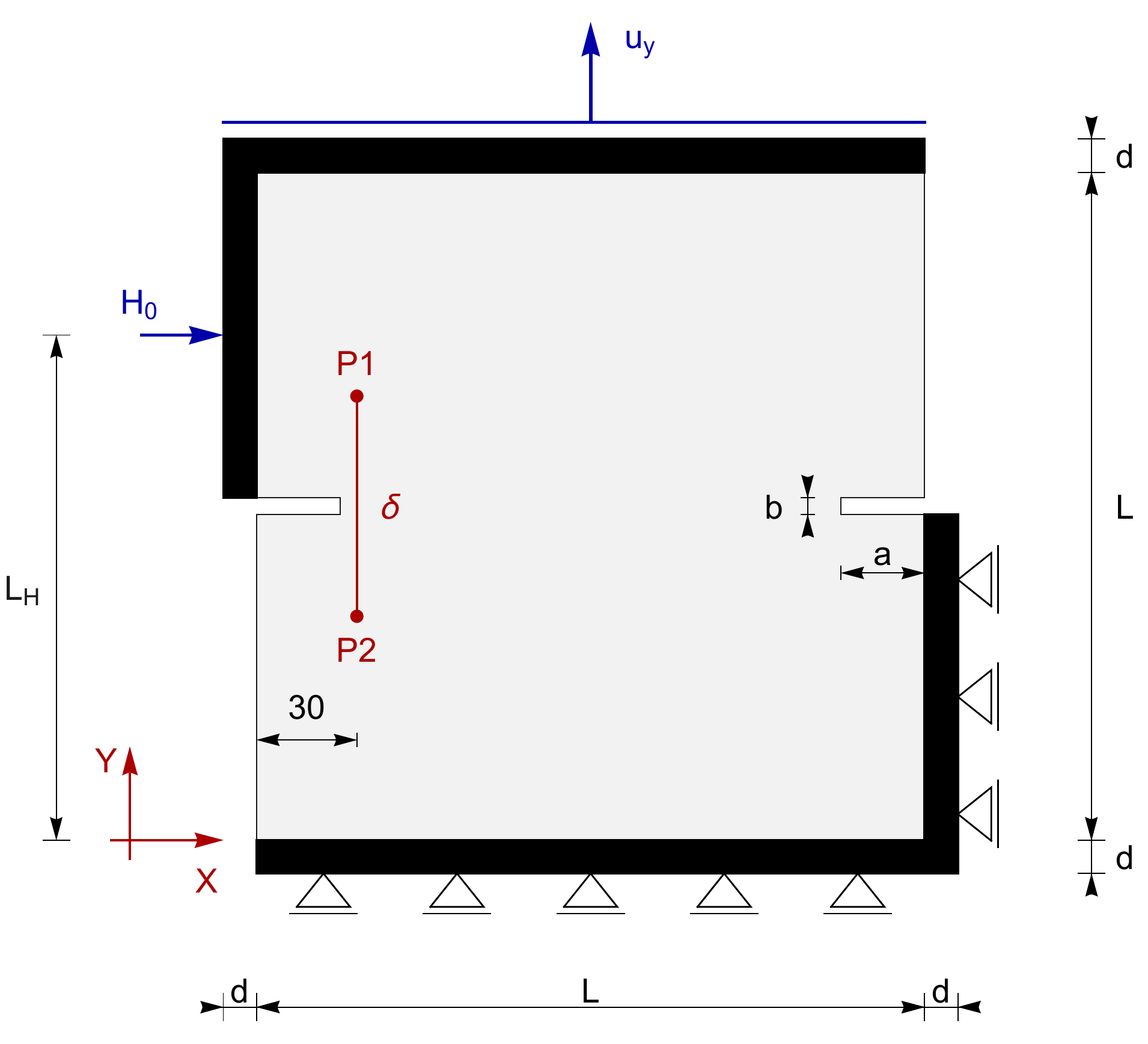}  
			\caption{Model A}
			\label{fig:DENSgeometryA}
		\end{subfigure}
		\hfil
		\begin{subfigure}{.49\textwidth}
			\centering
			% include second image
			\includegraphics[width=8.0 cm]{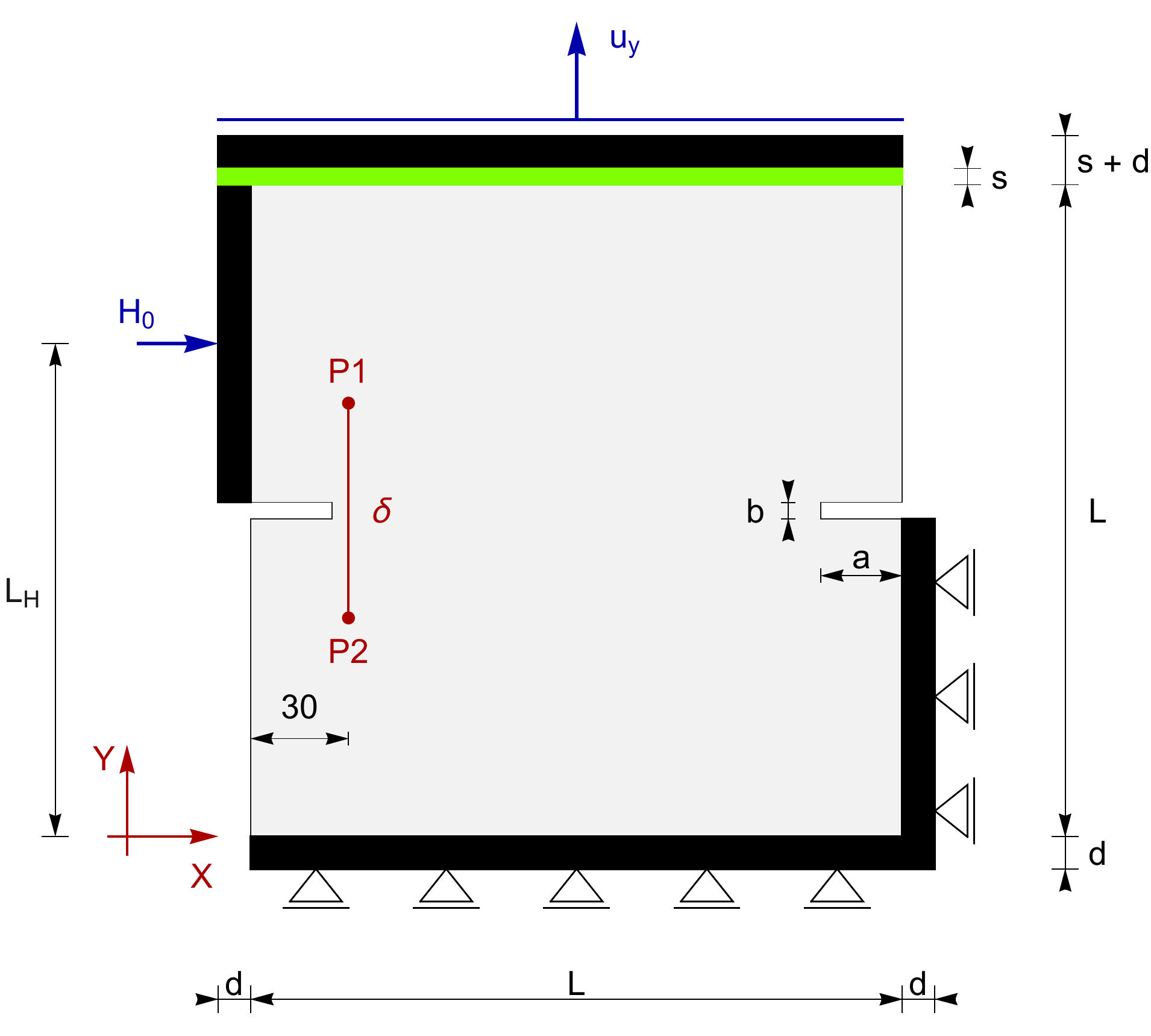}  
			\caption{Model B}
			\label{fig:DENSgeometryB}
		\end{subfigure}
		\caption{DENS: Geometry, boundary and loading conditions.}
		\label{fig:DENSgeometry}
	\end{figure}
%We reproduce 
A tensile-shear test, taken from \cite{NooruMohamed-phd}, is the last quasi-static example to be considered; there it is denoted as ``load path 4b (46-05)''. Fig.~\ref{fig:DENSgeometryA} shows the geometry as well as the loading and boundary conditions for the numerical model of the test performed on a concrete block with two notches ($a=25$ mm, $b=5$ mm) on opposite sides; it has a thickness of $t=50$ mm and a side length of $L=200$ mm. The black lines around the specimen in Fig.~\ref{fig:DENSgeometry} represent a steel loading frame with a width of $d=10$ mm and a thickness equal to that of the specimen. The concrete block and the frame are assumed to be in plane stress state.

	\begin{figure}[!hbt]
	\centering
	\begin{subfigure}{.49\textwidth}
		\centering
		% include first image
		\includegraphics[width=7 cm]{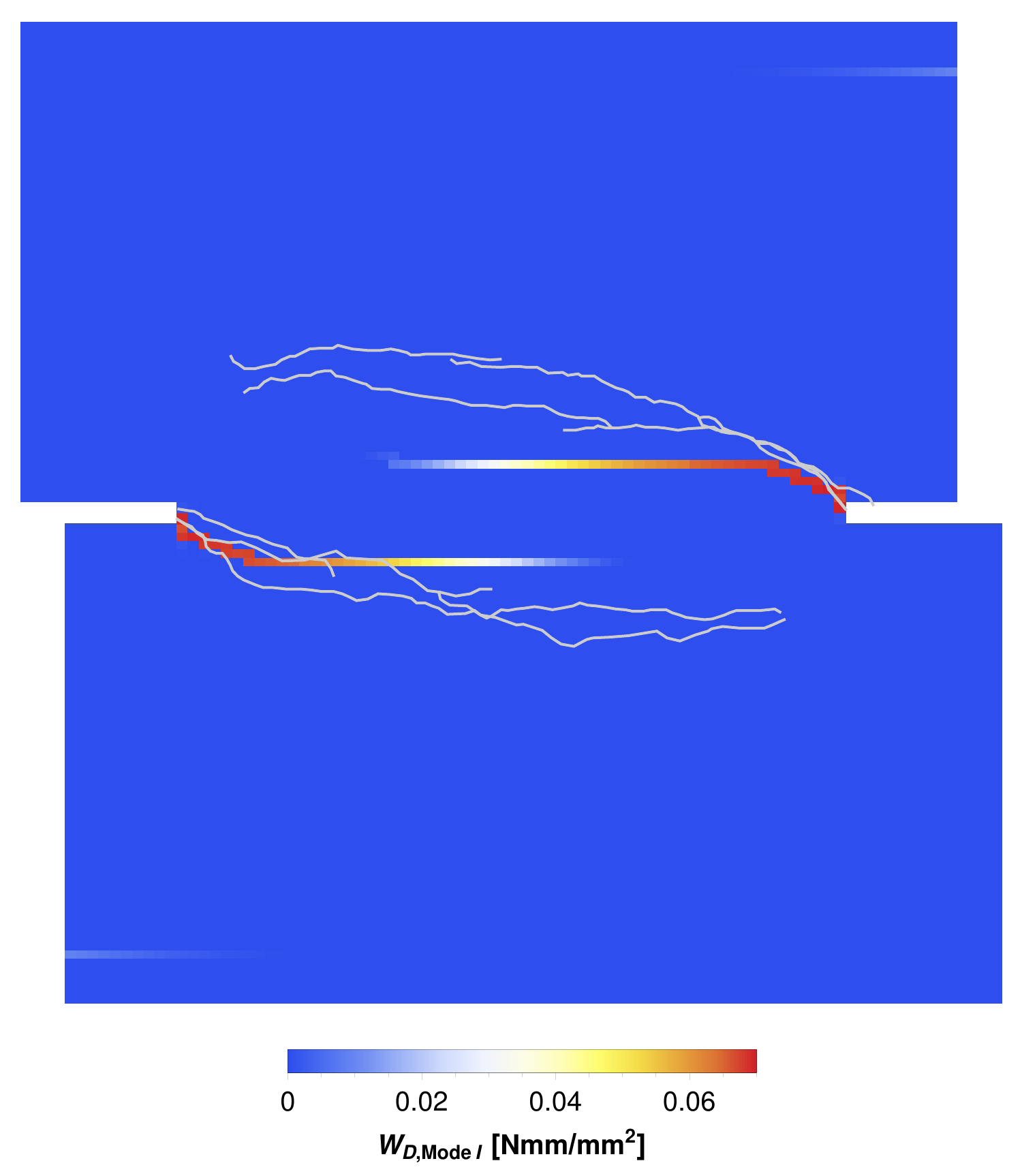}  
		%		\caption{Model A}
		%		\label{fig:DENSModeI_dissA}
	\end{subfigure}
	\hfil
	\begin{subfigure}{.49\textwidth}
		\centering
		% include second image
		\includegraphics[width=7 cm]{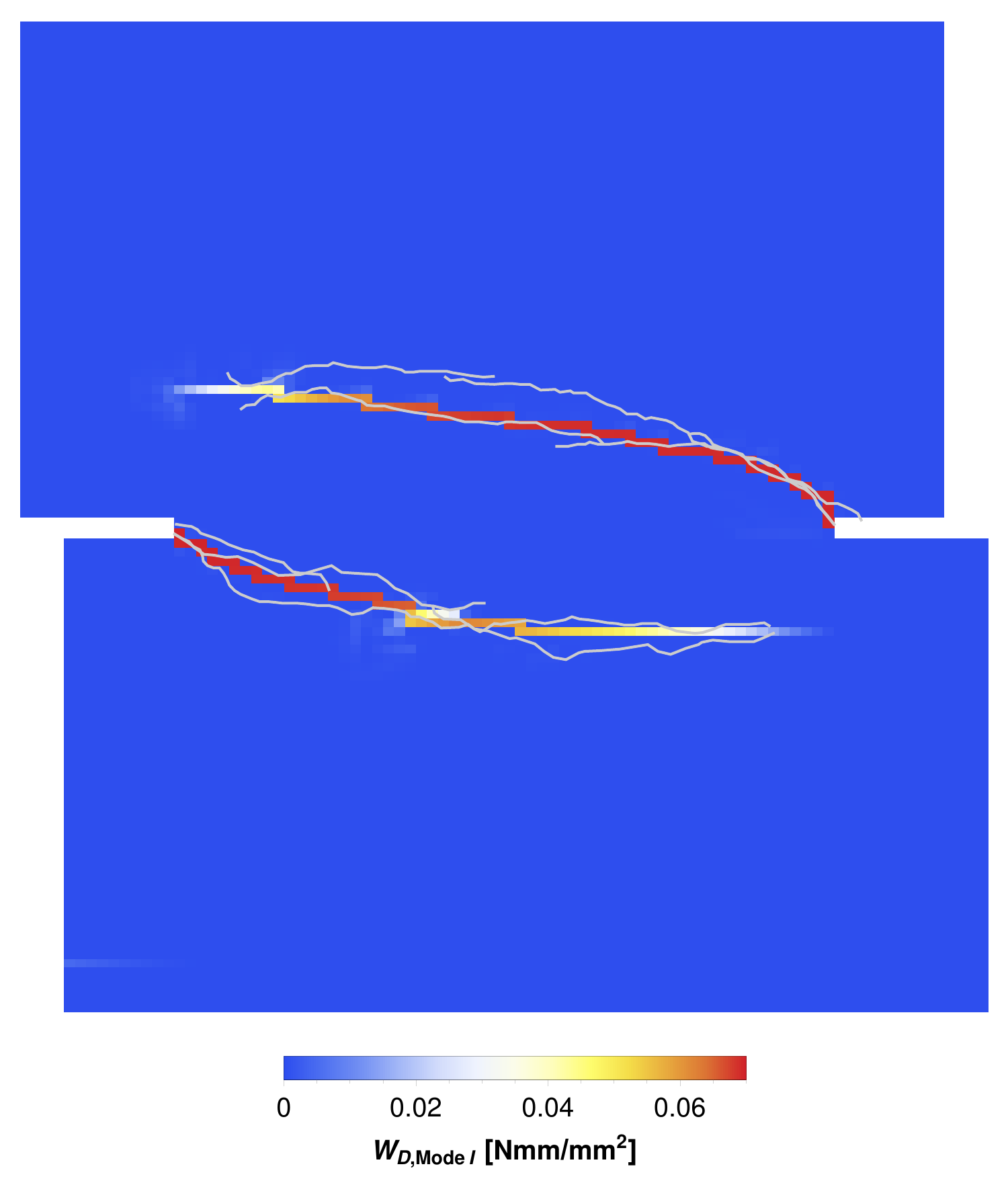}  
		%		\caption{Model B}
		%		\label{fig:DENSModeI_dissB}
	\end{subfigure}
	\caption{DENS: Specific dissipated mode I fracture energy for
		%due to separation in fashion of Mode I in 
		models A (left) and B (right).}
	\label{fig:DENSModeI_diss}
	\end{figure}

The loading is composed of constant horizontal force $H_{0}=10$ kN (the point of application is at $L_{H}=151.25$ mm) and a vertical imposed displacement $u_y=\lambda \, u_{y0}$ (where $u_{y0}=1$ mm and $\lambda$ is a load multiplier) along the top edge of the loading frame. The material properties for the concrete are: Young's modulus $E=30000$ N/mm$^2$, Poisson's ratio $\nu=0.2$, tensile strength $\sigma_{un}=3.0$ N/mm$^2$, and shear strength $\sigma_{um}=2.4$ N/mm$^2$. The mode I and mode II fracture energy are $G_{fn}=0.07$ Nmm/mm$^{2}$ and $G_{fm}=0.10$ Nmm/mm$^{2}$, respectively. For these material constants, the fracture energy $G_{fn}=0.07$ Nmm/mm$^{2}$ is too low to achieve reasonable numerical convergence, a similar problem already discussed in Subsection~\ref{sec:therm-shock}. Therefore, we apply the zig-zag form of the softening function, as in the thermal shock test example (see Fig. \ref{fig:ThermalShockGradualFunction} in Subsection~\ref{sec:therm-shock}), with the following data: $\xi_{O} = 0$, $\xi_{A} = 0.003$, $\xi_{B} = 0.006$, $\xi_{C} = 0.009$, $\xi_{D} = 0.016$ and $\xi_{E} = 0.023$;  $G_{fn}^{OA} = 0.100$, $G_{fn}^{AB} = 0.094$, $G_{fn}^{BC} = 0.085$, $G_{fn}^{CD} = 0.076$, $G_{fn}^{DE}=0.073$ and $G_{fn}^{EF}=G_{fn}=0.070$. The numerical simulation of the experimental test is performed in two phases. Firstly, with vertical displacements on the top edge fixed, the horizontal force is incrementally increased until it reaches its final value $H_{0}$. In the second phase, the horizontal force $H_{0}$ is fixed, and the vertical displacements $u_{y}$ of the top edge are incrementally increased.

%	\noindent
Initial numerical tests showed a high sensitivity of the crack paths to the exact boundary conditions. The report on the experiment in \cite{NooruMohamed-phd} states that there were two independent parallel loading frames, one frame for applying the horizontal force, and another one for the vertical load. For this reason, it may be that small rotations of the (non-active) frame occured during the loading. To model and investigate this possible influence, we prepared an additional numerical model with a thin soft layer of width $s=5$ mm below the top loading frame (see the green layer in Fig.~\ref{fig:DENSgeometryB}). The soft layer is modelled as an elastic material with Young's modulus $E=200$ N/mm$^2$ and Poisson's coefficient $\nu=0.3$, and it allows the top half of the concrete block to slightly rotate during the loading procedure. Thus, there are two models: model A considers the frame setting as stiff and rigidly connected to the specimen, whereas model B allows for a relative rotation between the vertical frame (where the force is applied) and the horizontal frame (where the imposed displacement is applied).

	\begin{figure}[!hbt]
	\centering
	\begin{subfigure}{.49\textwidth}
		\centering
		% include first image
		\includegraphics[width=7 cm]{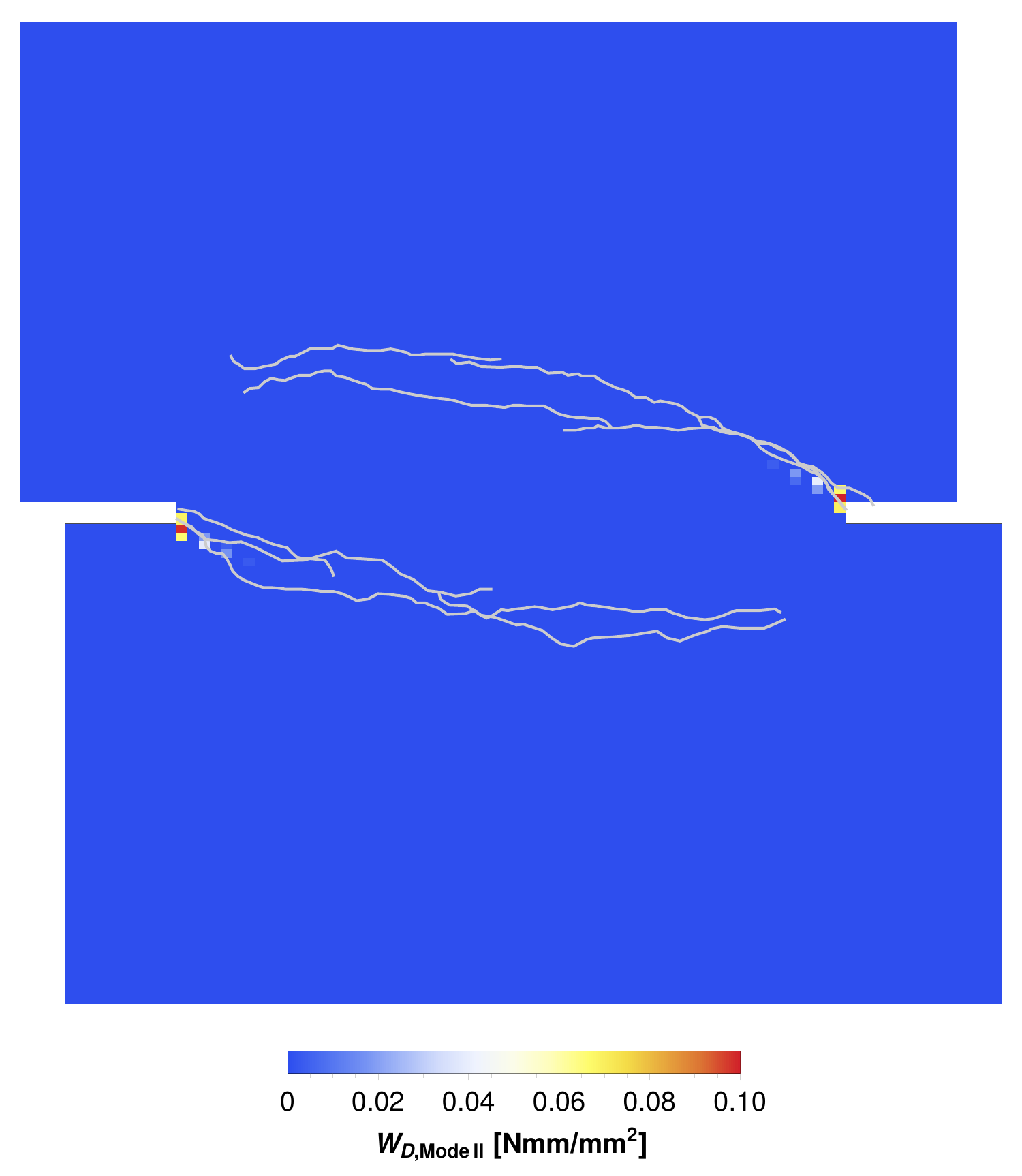}  
		%		\caption{Model A}
		%		\label{fig:DENSModeII_dissA}
	\end{subfigure}
	\hfil
	\begin{subfigure}{.49\textwidth}
		\centering
		% include second image
		\includegraphics[width=7 cm]{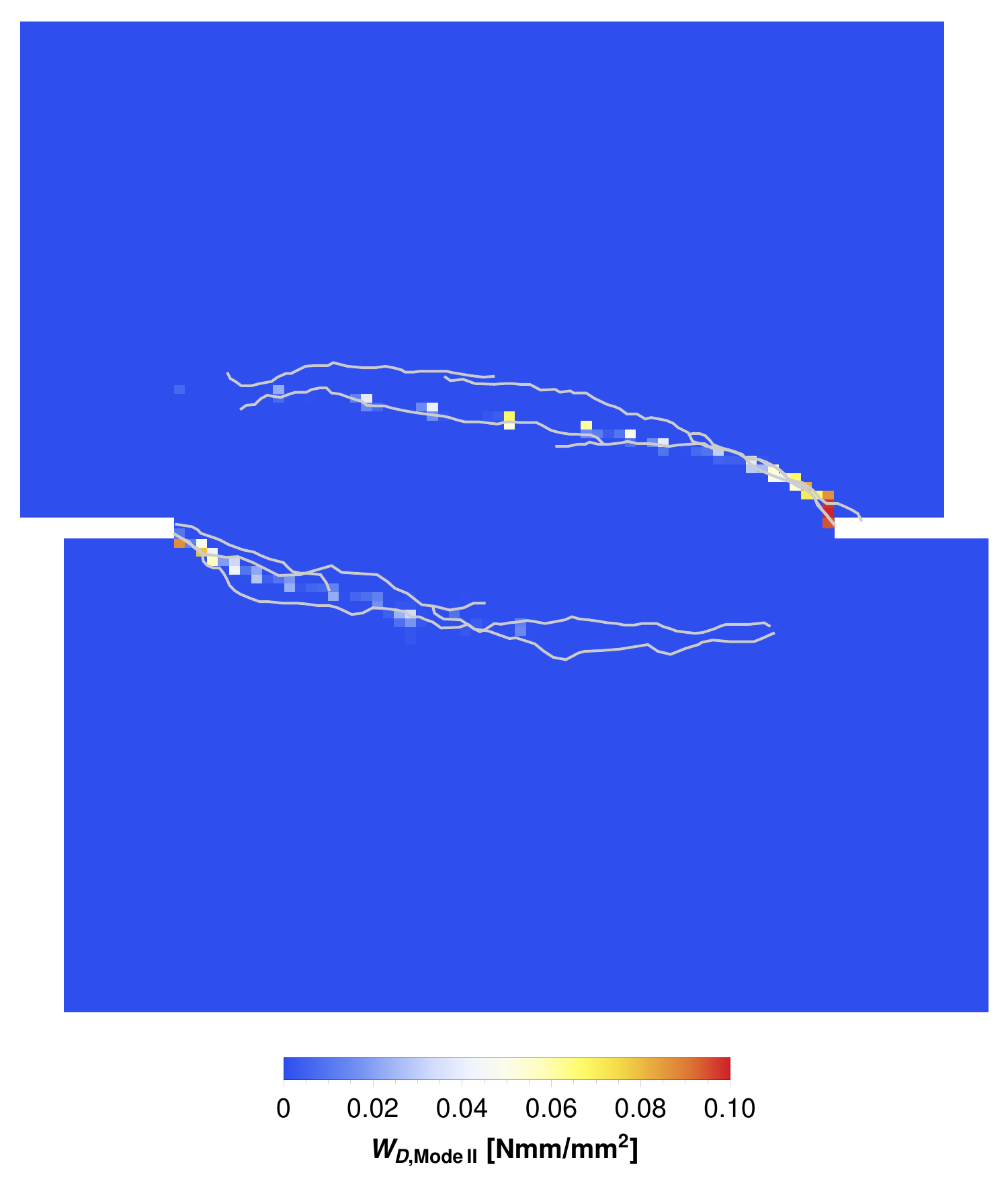}  
		%		\caption{Model B}
		%		\label{fig:DENSModeII_dissB}
	\end{subfigure}
	\caption{DENS: Specific dissipated mode II fracture energy for 
	models A (left) and B (right).}
	\label{fig:DENSModeII_diss}
	\end{figure}

\begin{figure}[!hbt]
	\centering
	\begin{subfigure}{.49\textwidth}
		\centering
		% include first image
		\includegraphics[width=7 cm]{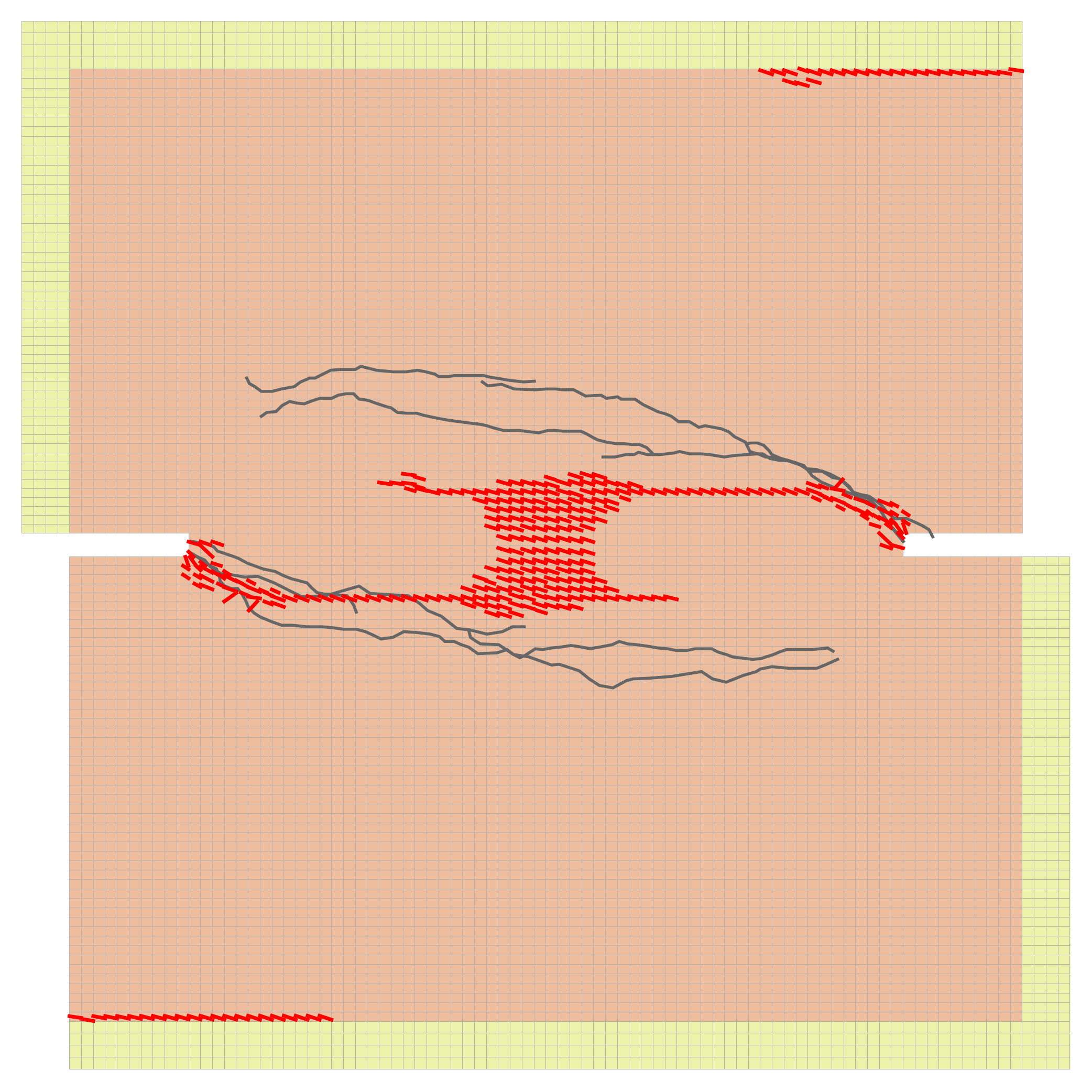}  
		%			\caption{Model A}
		%			\label{fig:DENSModel_crackA}
	\end{subfigure}
	\hfil
	\begin{subfigure}{.49\textwidth}
		\centering
		% include second image
		\includegraphics[width=7 cm]{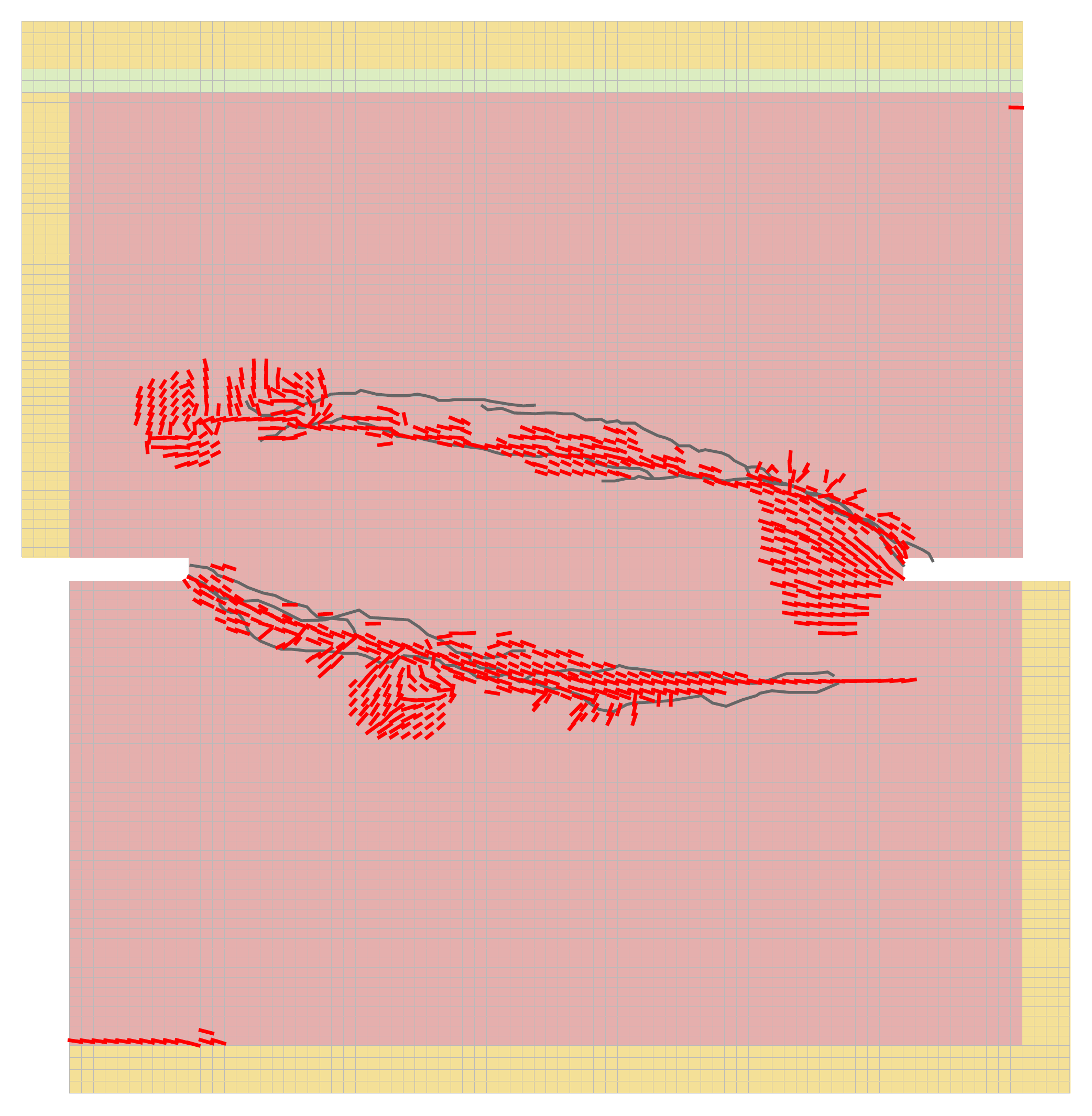}  
		%			\caption{Model B}
		%			\label{fig:DENSModel_crackB}
	\end{subfigure}
	\caption{DENS: Nucleated cracks in model A (left) and model B (right) in comparison 
	with experimentally observed crack paths (grey continuous lines) from \cite{NooruMohamed-phd}.}
	\label{fig:DENSModel_crack}
	\end{figure}

\begin{figure}[!hbt]
	\centering
		\begin{subfigure}{.49\textwidth}
			\centering
			% include first image
			\includegraphics[width=7 cm]{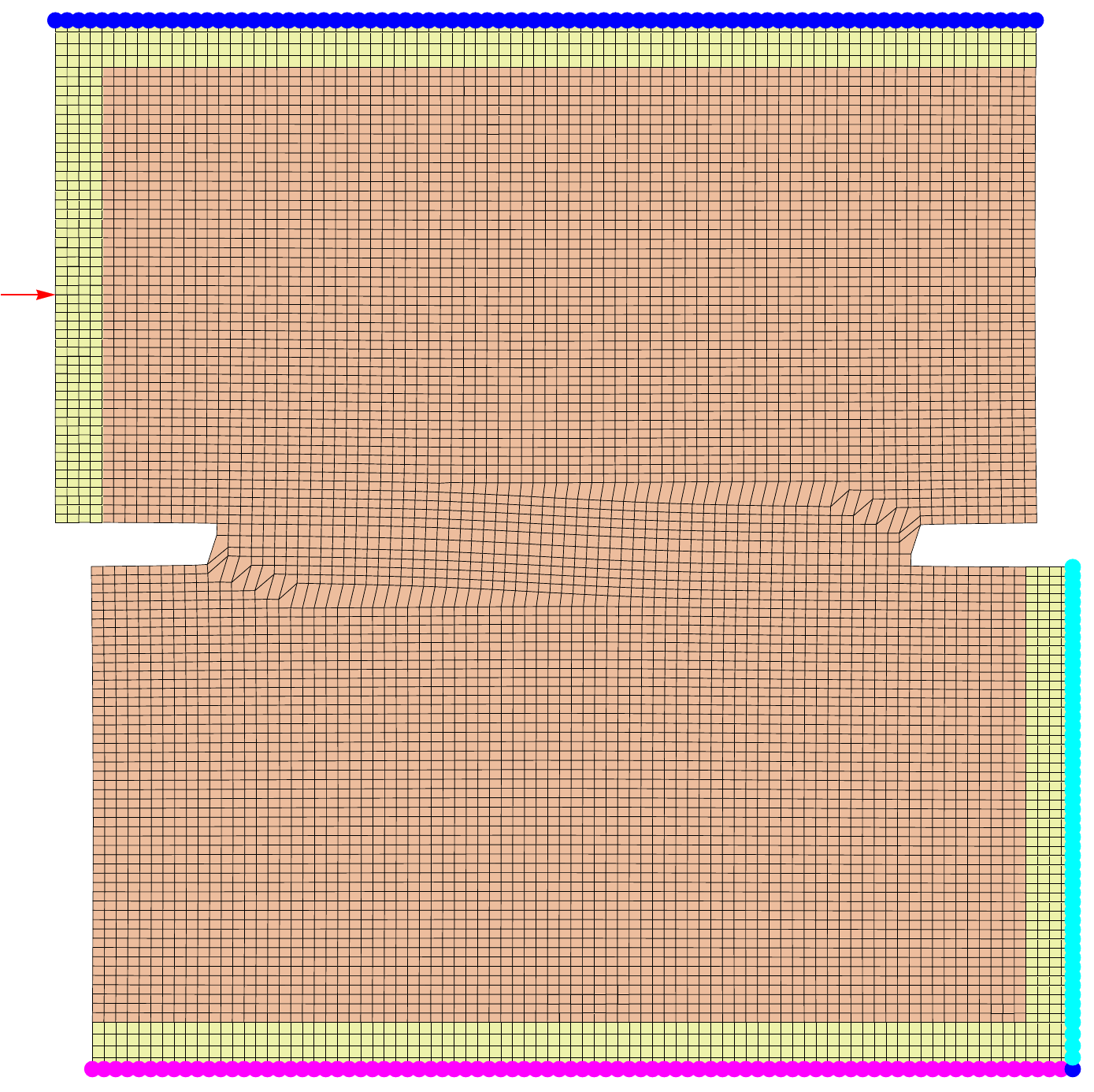}  
			%		\caption{Model A}
			%		\label{fig:DENSModel_defConfA}
		\end{subfigure}
		\hfil
		\begin{subfigure}{.49\textwidth}
			\centering
			% include second image
			\includegraphics[width=7 cm]{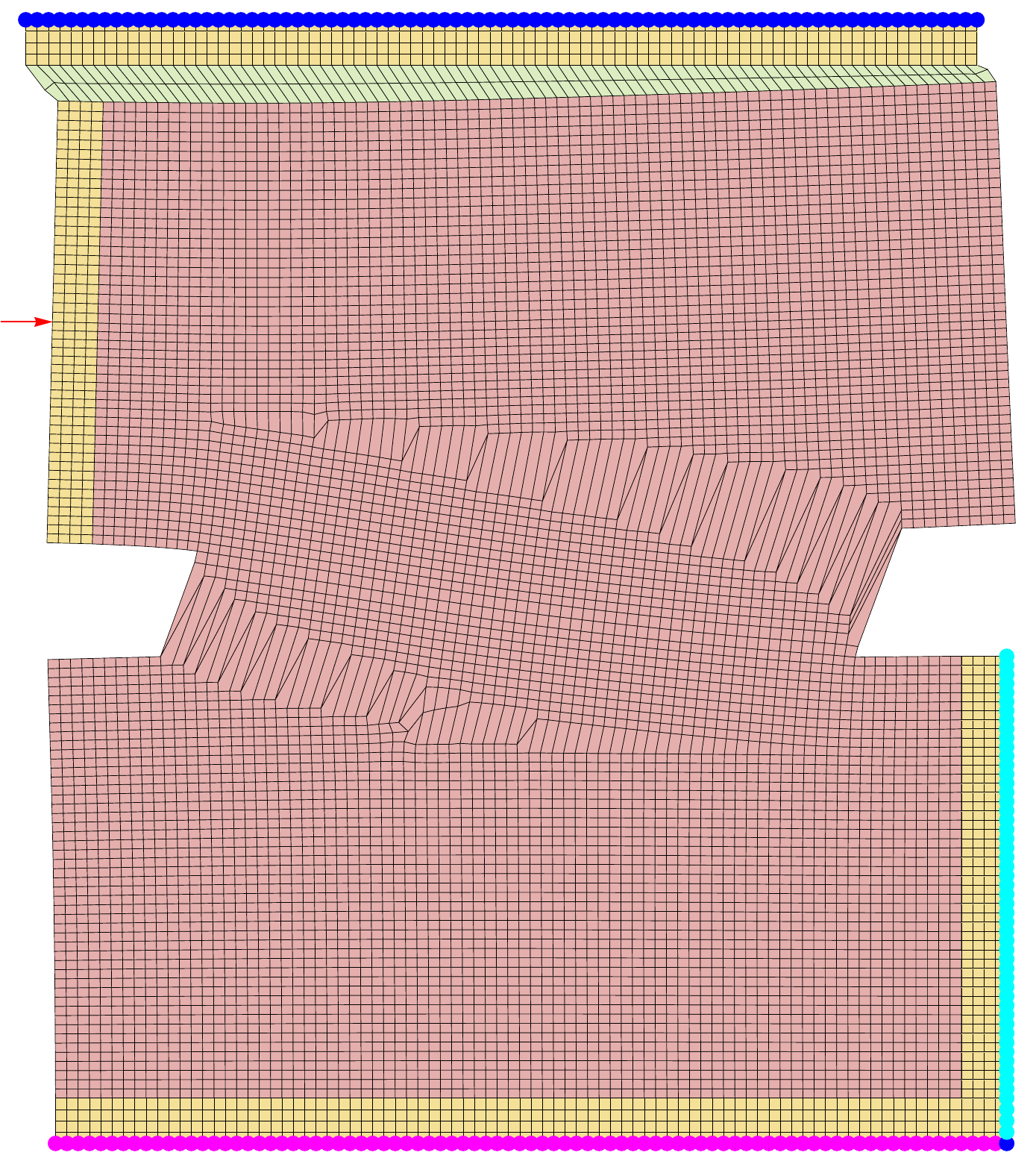}  
			%		\caption{Model B}
			%		\label{fig:DENSModel_defConfB}
		\end{subfigure}
		\caption{DENS: Deformed meshes (150$\times$ magnified) at the end of simulation 
		(at $\lambda = 0.03 $ and $\lambda = 0.14 $) for model A (left) and model B (right). 
		%\fixAS{Popravljeno.}
		}
		\label{fig:DENSModel_defConf}
	\end{figure}

% \noindent 
Fig.~\ref{fig:DENSModeI_diss} displays the cracks and the dissipated fracture energy for models A and B, and compares the so-obtained patterns with the experimentally observed crack paths. In the experiment, one crack path formed from each notch. Also for models A and B, one crack path propagates from each notch. However, for model A, both crack paths are flat and do not match the experimentally observed ones. The crack paths for model B match the experimentally observed crack lines much better, which implies that in this example the boundary conditions play an important role. The dissipated fracture energies in Figs.~\ref{fig:DENSModeI_diss} and~\ref{fig:DENSModeII_diss} show that the cracks along the crack path open and slide in mode I and in mode II manners. Mode I is the dominant separation mode, while mode II is active at the beginning of the crack lines around the notch area. Fig.~\ref{fig:DENSModel_crack} displays all nucleated cracks in both numerical models. The ultimate tensile stress is exceeded in many elements, but the \emph{active} embedded discontinuities can be easily determined by checking for the elements which dissipate fracture energy (see Fig.~\ref{fig:DENSModeI_diss}). The crack paths can also be observed from the deformed configurations of the finite element mesh in Fig.~\ref{fig:DENSModel_defConf}.

\begin{figure}[H]
	\centering
	\begin{subfigure}{.49\textwidth}
		\centering
		% include first image
		\includegraphics[width=8 cm]{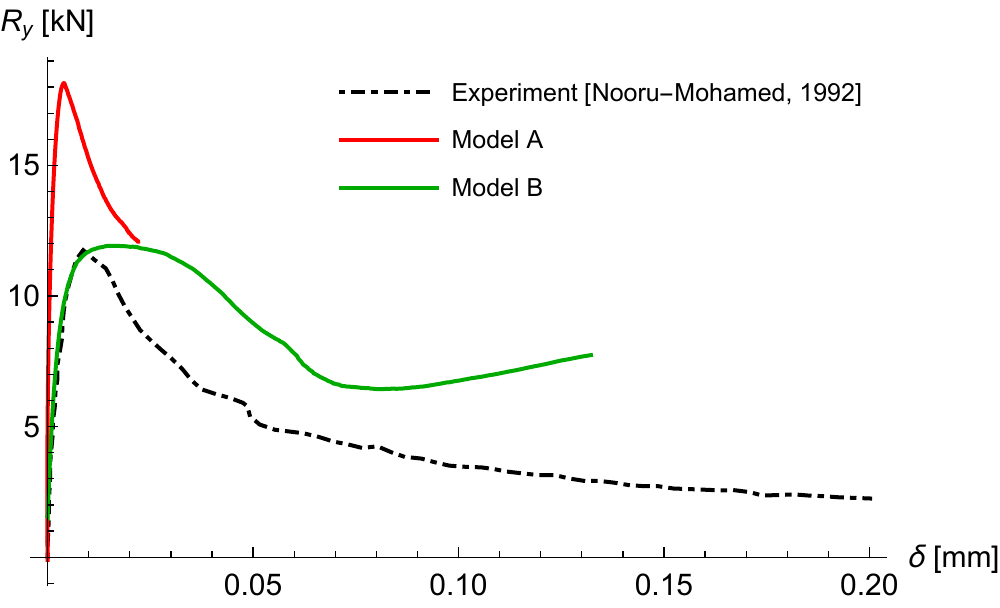}  
		\caption{Comparison with experiment.}
		\label{fig:DENS_Response01}
	\end{subfigure}
	\hfil
	\begin{subfigure}{.49\textwidth}
		\centering
		% include second image
		\includegraphics[width=8 cm]{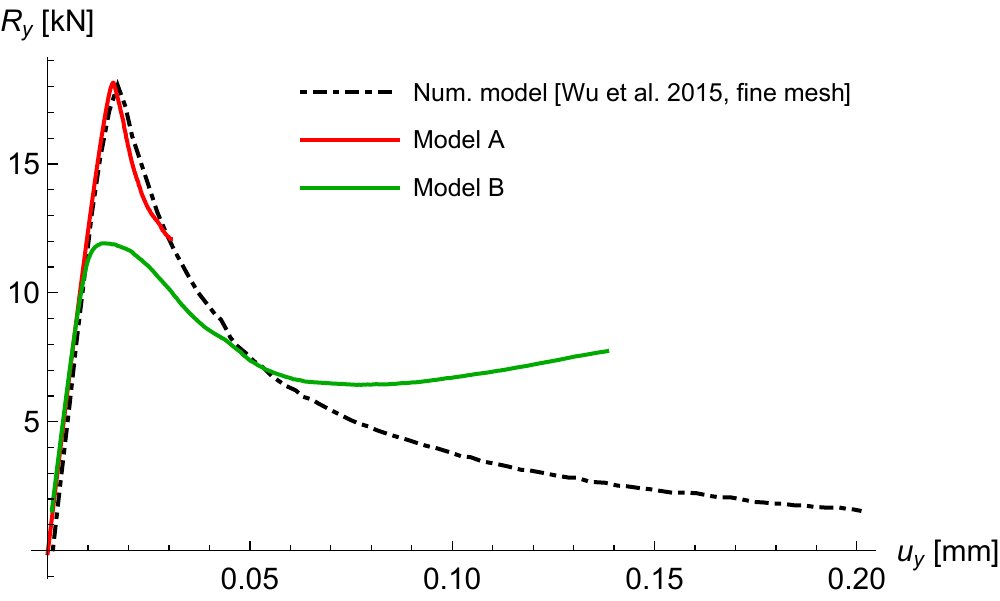}  
		\caption{Comparison with another numerical simulation.}
		\label{fig:DENS_Response02}
	\end{subfigure}
	\caption{DENS: Reaction force versus vertical displacements $\delta$ and $u_y$}
	\label{fig:DENS_Response}
	\end{figure}

%	\noindent 
Graphs of the computed results are presented in Fig.~\ref{fig:DENS_Response}. Reaction force versus relative displacement $\delta$ for models A and B is shown in Fig.~\ref{fig:DENS_Response01} in comparison with measurements from the experiment. In \cite{NooruMohamed-phd}, the author measured the difference in distance $\delta$ between two points: P1 with coordinates $(x,y)=(30, 132.5)$ mm, and P2 with coordinates $(x,y)=(30, 67.5)$ mm (see Fig.~\ref{fig:DENSgeometry}). In comparison with the experiment, model A has the same initial slope, but it achieves a considerably higher limit load. While model B exhibits the same initial stiffness, it also reaches similar limit resistance in comparison with the experiment. This suggests that model B indeed simulates the experimental conditions better than model A. Model B also produces a longer force-displacement curve than model A, but at some point it starts increasing. The reason for this is probably some kind of locking, because only constant separation modes are used in our formulation (see Fig.~\ref{fig:Q4_SepModes}). In Fig.~\ref{fig:DENS_Response02}, the computed response of the present models is compared with the numerical simulation from~\cite{Wu et al-art}, where the authors used a crack tracking algorithm and a formulation with constant and linear separation modes, where opening and sliding parameters were treated as global unknowns. The curves show reaction force versus applied vertical displacement $u_{y}$. Model A has quite similar response as the model in \cite{Wu et al-art}, but the response of the experiment is better matched with model B. The curve of model A stops early because of the loss of convergence in numerical analysis.

	\subsection{Kalthoff's test}
	\label{sec:Kalthoff} 
	
\begin{figure}[H]
		\centering
		\includegraphics[width=12cm]{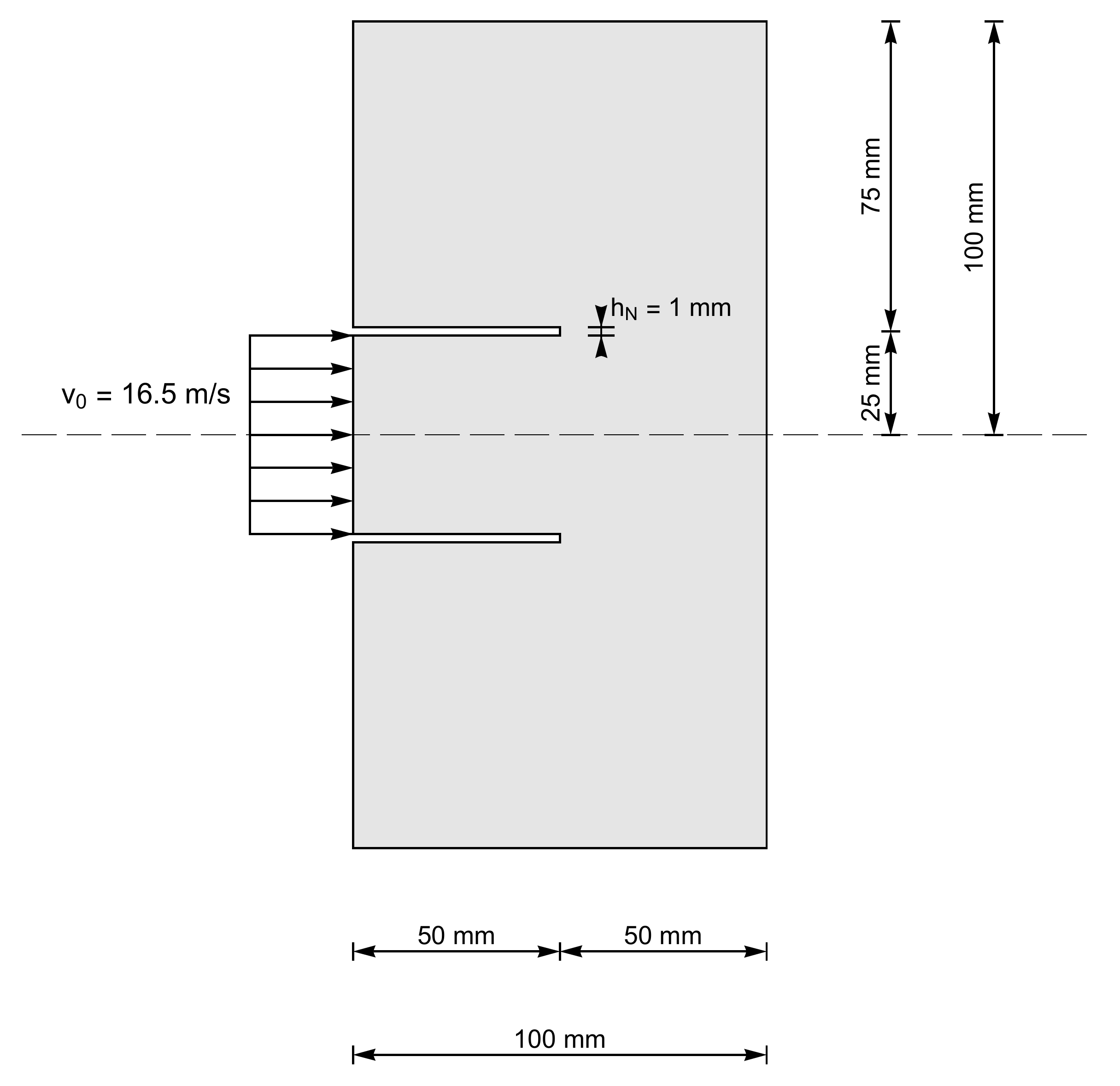}
		\caption{Kalthoff's test: Geometry and boundary conditions.} 
		\label{fig:KalthoffGeometry}
	\end{figure}

	The experiment from \cite{Kalthoff 2000-art}, known as Kalthoff's test, was simulated as \correction{an example of dynamic fracture}. A thin, high strength steel plate of thickness $t^e=9$ mm, with two parallel notches (see Fig.~\ref{fig:KalthoffGeometry}), is exposed to lateral impact load by a projectile that was launched from a gun muzzle with velocity of $33$ m/s. The experimental results showed that from each notch a tensile crack propagated at an angle of about $70^{\circ}$ with respect to the notch line. When preparing the numerical model, we took into account the assumption from \cite{Lee Freund 1989-art}, saying that the specimen and the projectile have the same elastic impedance, and thereby the impact load of the projectile can be modelled by an imposed constant velocity of $v_{0}=16.5$ m/s --- which is one-half of the projectile's velocity. The imposed velocity $v_{0}$, which is applied on the part of the edge between the two notches for $100 \, \upmu$s, is the only boundary condition in the finite element model. The material properties of the high strength steel are taken from \cite{Lloberas-Valls et al 2016-art}: Young's modulus $E=190 \cdot 10^3$ N/mm$^2$, Poisson's ratio $\nu=0.3$, tensile strength $\sigma_{un}=844$ N/mm$^2$, mode I fracture energy $G_{fn}=22.17$ Nmm/mm$^{2}$, mass density $\bar{\rho} = 8000 $ kg/m$^{3}$, and a Rayleigh wave speed of $v_R = 2.799 \cdot 10^6$ mm/s. The shear strength and mode II fracture energy are chosen as $\sigma_{um}=400$ N/mm$^2$ and $G_{fm}=12$ Nmm/mm$^{2}$, respectively. We prepared three meshes (with 3280, 7320 and 19800 elements) for plane strain as well as for plane stress conditions. Dynamic computations were performed by Newmark's standard implicit time-stepping scheme, the trapezoidal rule (see Section~\ref{sec:Dyn}).
	
\begin{figure}[H]
		\centering
		\includegraphics[width=\linewidth]{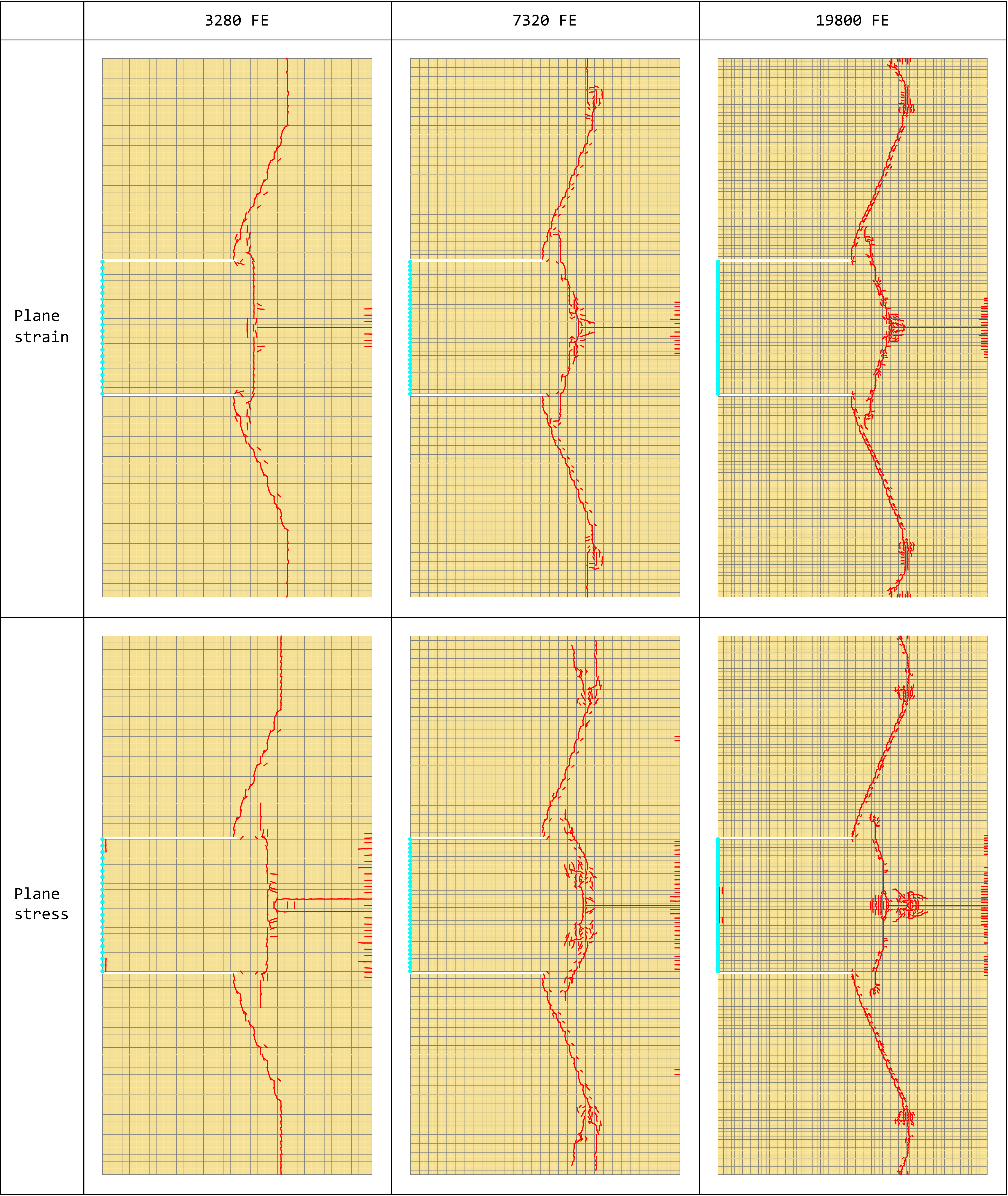}
		\caption{Kalthoff's test: Considered meshes with nucleated cracks.} 
		\label{fig:Kalthoff_Cracks}
	\end{figure}

\begin{figure}[!hbt]
		\centering
		\subfloat[Plane strain.]{{
				\includegraphics[width=6.5cm]{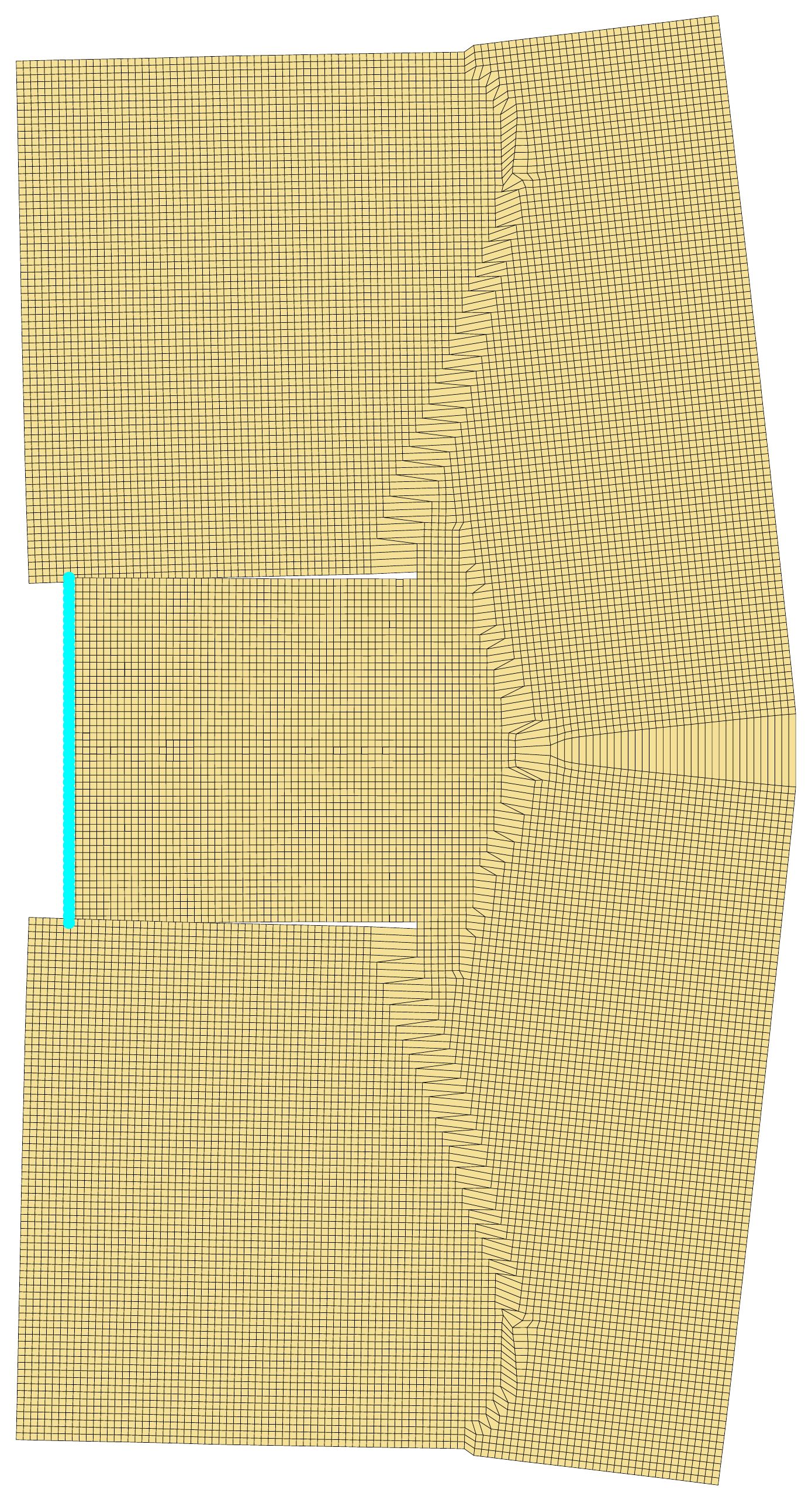} }}%
		\qquad
		\subfloat[Plane stress.]{{
				\includegraphics[width=6.5cm]{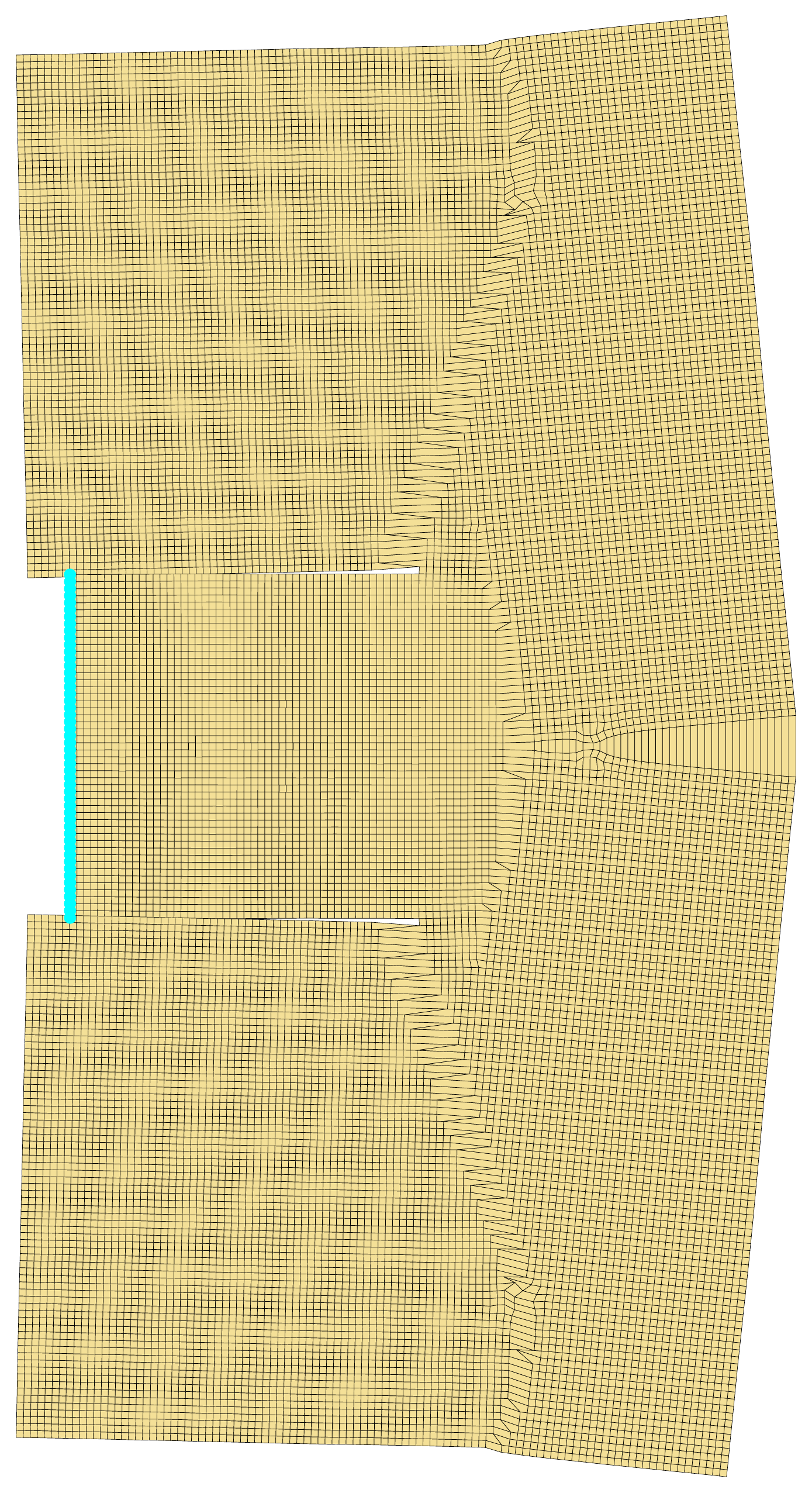} }}%
		\caption{Kalthoff's test: Deformed mesh (deformation magnified 5$\times$) at the end of   
		the simulation.}%
		\label{fig:Kalthoff_DefMesh}%
	\end{figure}

\begin{figure}[H]
		\centering
		\includegraphics[width=\linewidth]{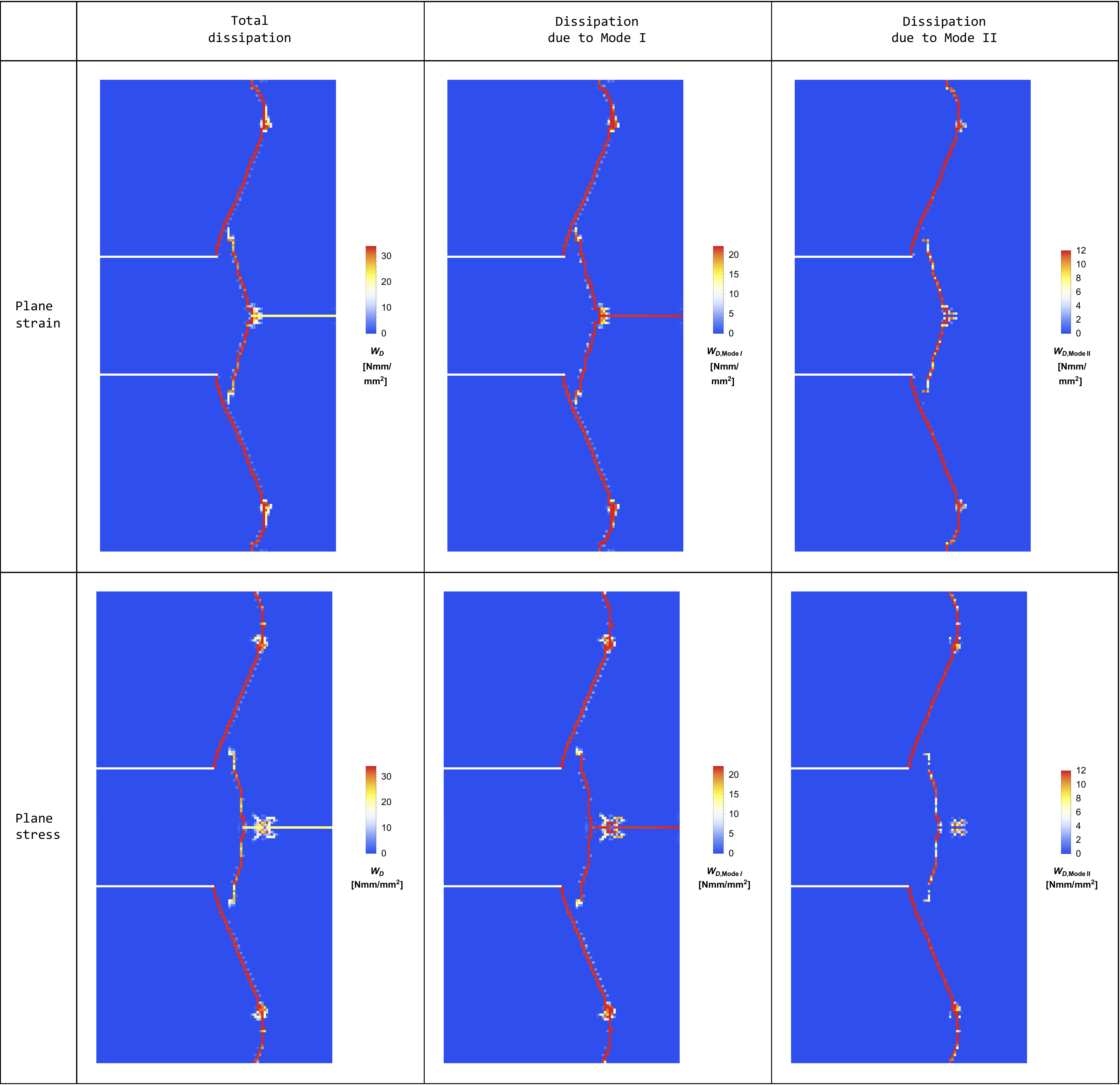}
		\caption{Kalthoff's test: Specific dissipated fracture energy (mesh with 19800 elements).} 
		\label{fig:Kalthoff_Diss}
	\end{figure}

\begin{figure}[H]
		\centering
		\subfloat[Plane strain.]{\label{fig:Kalthoff_PlaneStrain_DissGraph}			
			\includegraphics[width=0.96\linewidth]{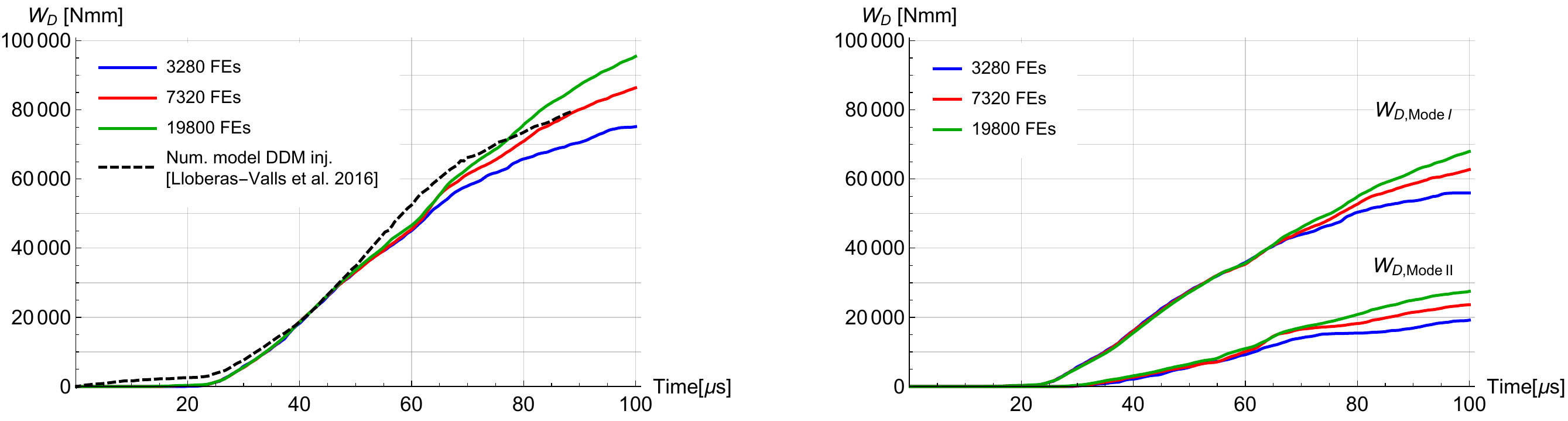}}
		\par\medskip
		\centering
		\subfloat[Plane stress.]{\label{fig:Kalthoff_PlaneStress_DissGraph}	
			\includegraphics[width=0.96\linewidth]{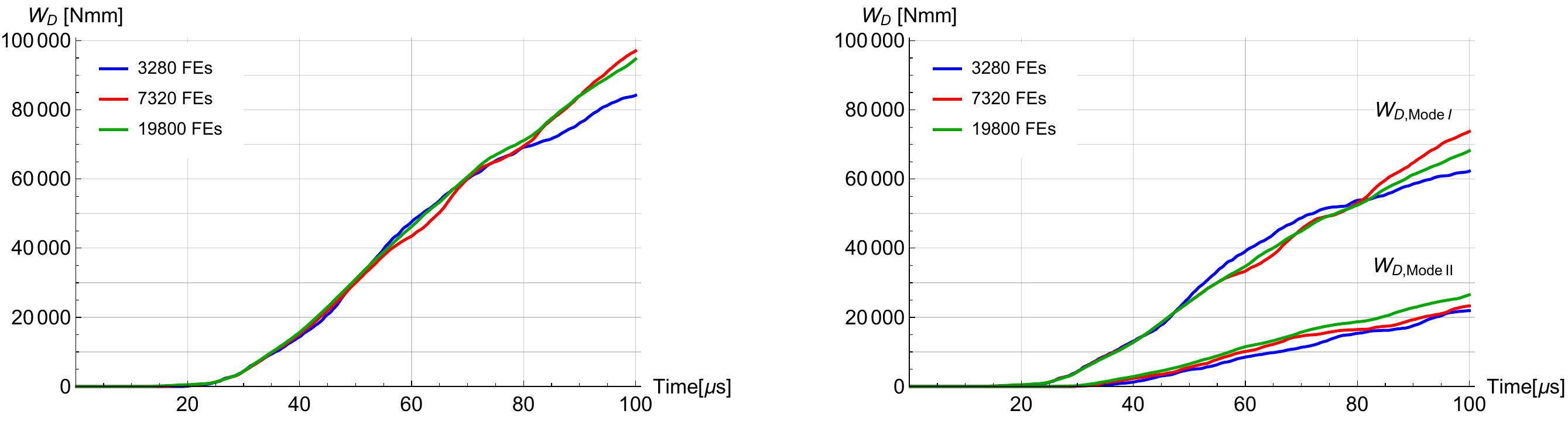}}
		\caption{Kalthoff's test: Dissipated energy. 
			%(a) Plane strain. (b) Plane stress.
		}
		\label{fig:Kalthoff_DissGraph}
	\end{figure}

\begin{figure}[H]
		\centering
		\includegraphics[width=16cm]{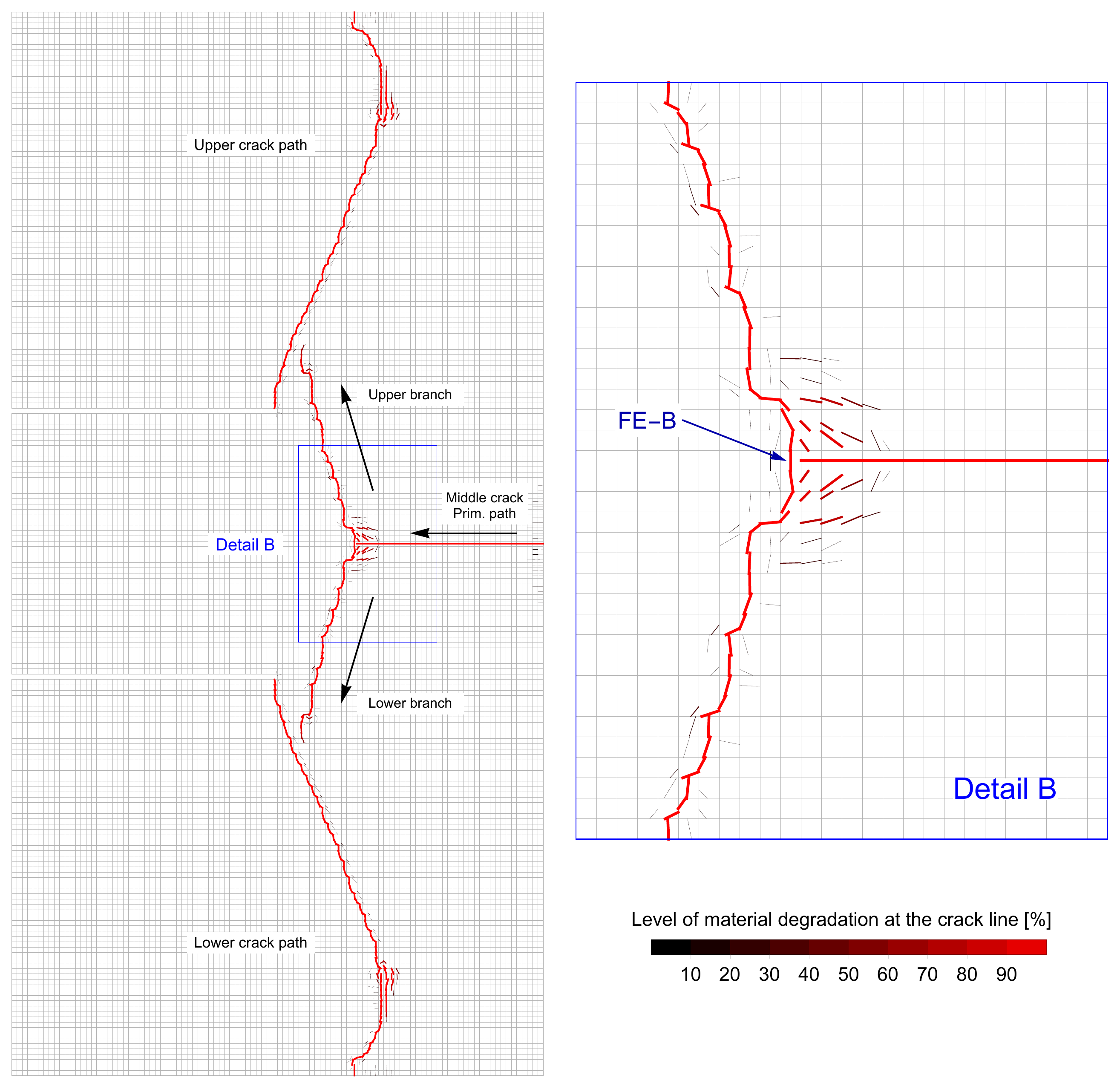}
		\caption{Kalthoff's test: Details of crack paths in the plain strain model 
		with 19800 finite elements.}
		\label{fig:Kalthoff_Branching}
	\end{figure}

\begin{figure}[H]
		\centering
		\subfloat[Middle crack]{\label{fig:Kalthoff_MiddleCrack_DissGraph}
			\includegraphics[width=10cm]{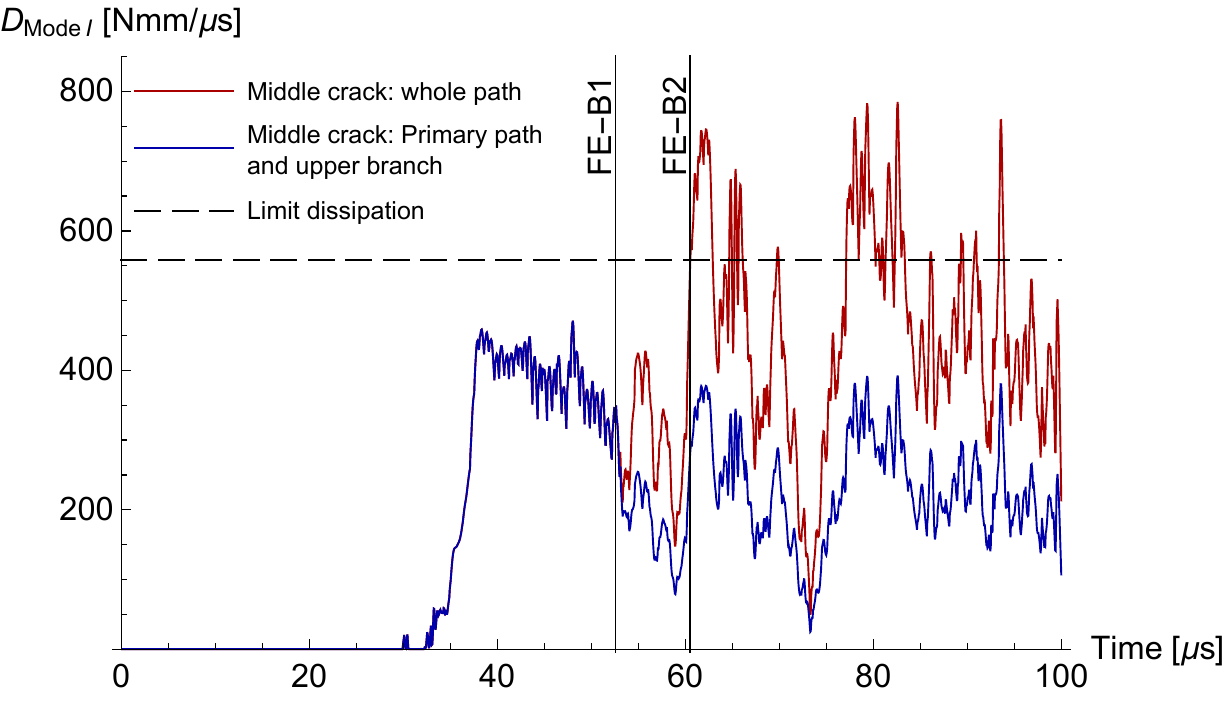}}
		\par
		\medskip
		\centering
		\subfloat[Upper crack]{\label{fig:Kalthoff_UpperCrack_DissGraph}
			\includegraphics[width=10cm]{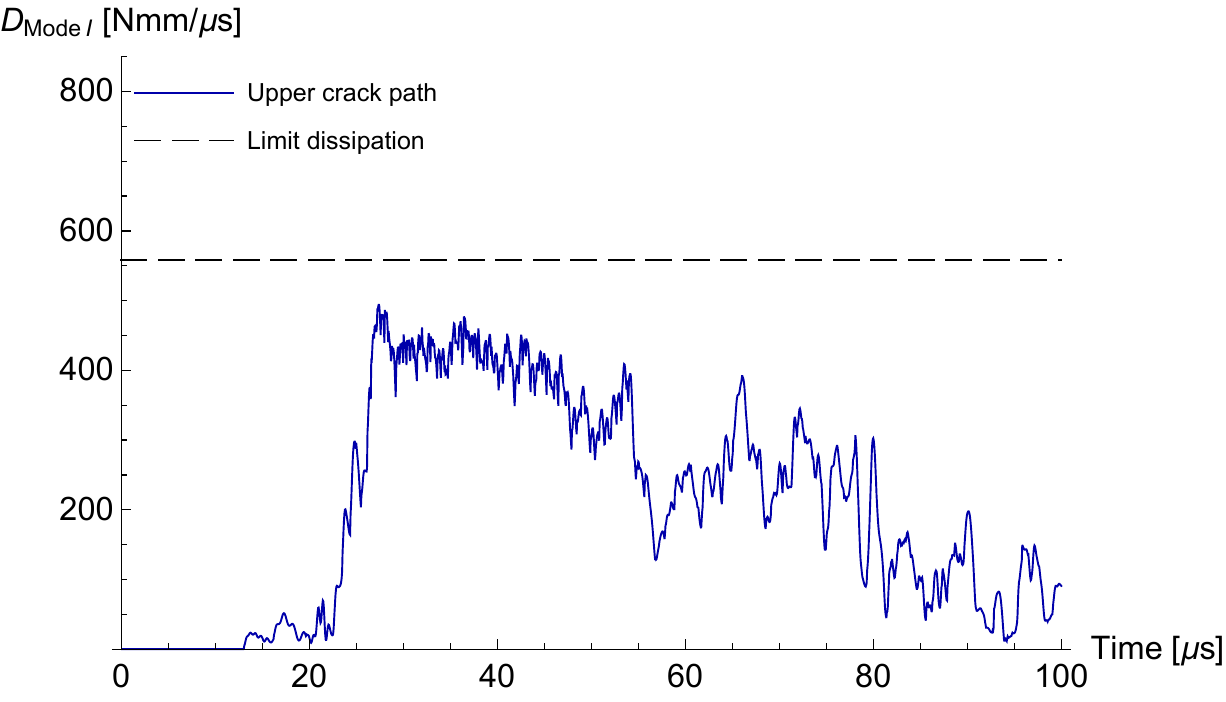}}
		\caption{Kalthoff's test. Dissipation due to mode I separation in plane strain model with 19800 finite element.}
		\label{fig:Kalthoff_DissGraph_Crack}
	\end{figure}

	%	\noindent 
	Results are presented in Figs.~\ref{fig:Kalthoff_Cracks}--\ref{fig:Kalthoff_DissGraph_Crack}. Patterns of nucleated cracks at the end of the analysis may be gleaned from Fig.~\ref{fig:Kalthoff_Cracks}, which considers the first $100 \, \upmu$s of the impact. Both the plane stress and the plane strain models capture all parts of the complex crack pattern well for all meshes. Two straight lines of cracks propagate from each notch at an angle of about $70^{\circ}$, and towards the end of the propagation they kink towards the edges parallel to the impact direction. Moreover, another straight line --- hereinafter called the middle crack --- propagates from the edge opposite to the impacted edge towards the middle area of the plate, and, at some point, it splits into an upper and bottom branch (see also Fig.~\ref{fig:Kalthoff_Branching}). All three plane strain meshes have very similar crack patterns. The main influence of the mesh density is that for the denser meshes in some parts the crack lines are slightly more curved. On the other hand, the crack patterns for the plane stress meshes show some differences when compared to each other.
	
	For the coarsest mesh (with 3280 finite elements), the first part of the middle crack exhibits two parallel straight lines close to each other. The reason for this is a bifurcation due to the stress state in the part of the edge where the middle crack starts. When the stress wave in the solid caused by the impact load is reflected on the opposite side of the specimen, there is a group of finite elements with stresses that exceed the material strength. Therefore, the simulation may predict one crack line --- as is the case for the two denser meshes --- or two close parallel crack lines for the primary part of the middle crack, as is the case for the coarsest mesh. Furthermore, the denser two meshes (with 7320 and 19800 finite elements) have a vague branching point for the plane stress model. It seems that there are two branching points, and that the upper and lower branches choose to propagate from one of them. Additionally, for the plane stress mesh with 7320 finite elements, the crack that goes from the notches forms a secondary crack at the kink. Fig.~\ref{fig:Kalthoff_DefMesh} illustrates the deformed state for the plane strain and the plane stress meshes, where all parts of the complex crack paths are clearly visible as well.

To allow easier comparison of dissipated energy, both plane stress and plain strain models were computed with the thickness of 9 mm and this is shown in the following graphs. Fig.~\ref{fig:Kalthoff_Diss} depicts dissipated energy per unit crack surface for the finest mesh. It shows that both mode I and mode II are active; however, the dominant separation mode is mode I. Only the primary path of the middle crack propagates solely in mode I fashion. Graphs of the total dissipated energy are given in Fig.~\ref{fig:Kalthoff_DissGraph}. Both plane stress and plane strain models dissipate a similar amount of fracture energy (for the same mesh). In Fig.~\ref{fig:Kalthoff_PlaneStrain_DissGraph} (left), we compare our results with those from \cite{Lloberas-Valls et al 2016-art}, where the authors use a much more complicated numerical fracture model with strain injection techniques with discontinuous displacement fields (DDM inj.) and a crack tracking algorithm. The numerical model from \cite{Lloberas-Valls et al 2016-art} consists of approximately 15000 finite elements. One may observe that our curves for the total dissipated fracture energy are similar to those from \cite{Lloberas-Valls et al 2016-art}.

%	\noindent 
Crack branching is usually a great modelling challenge in fracture analysis, especially when a reliable criterion to determine the branching point has to be applied. Note that no branching criteria need to be used here to compute crack patterns with branching, as shown in the above figures. Let us elaborate a bit more on the crack branching in the case of the plane strain model with the finest mesh, see Fig.~\ref{fig:Kalthoff_Branching}, where the detail (called Detail B) of the crack branching area is also shown. Detail B shows that the primary path of the middle crack extends from the edge of the mesh to the finite element denoted as ``FE-B'', which has the discontinuity line perpendicular to the orientation of the primary path of the middle crack. The upper and the lower crack branches propagate from ``FE-B''. The cracked finite elements, which lie in the area of the middle crack path, may be divided into three groups: primary path, upper and lower branch. Observe the velocity of dissipated fracture energy in these finite elements. Fig.~\ref{fig:Kalthoff_MiddleCrack_DissGraph} displays the dissipation of fracture energy in mode I for each time step, where two curves are shown: the red curve is dissipation in all elements from the middle crack path, and the blue line represents the contribution of finite elements from the primary path and the upper branch. We compute the rate of dissipated fracture energy (or dissipation) for mode I as:
	\begin{equation}
		\label{eq:DissModeI}
		D_{Mode I} = \frac{\upDelta W_{D,Mode I}}{\upDelta \tau}  ,
	\end{equation}
where $\upDelta W_{D,Mode I}$ is the incremental growth of the dissipated energy, and $\upDelta \tau$ is the time increment (in our simulation $\upDelta \tau = 0.05 \upmu$s). One can see that the red and blue curves are identical until $t_{FE-B1} = 52.60 \,\upmu$s, when the crack in ``FE-B'' is nucleated, and after that the red curve is above the blue one. The experimental results (see \cite{Lloberas-Valls et al 2016-art} and references therein) report that the critical crack tip velocity is below the Rayleigh wave velocity $v_R$ (which defines the limit speed of elastic waves on the surface of a solid domain). When the critical crack tip velocity is exceeded, crack path branches and new crack surfaces are formed in a short time interval, which enables larger dissipation of fracture energy. The dissipation rate depends on the crack tip velocity, so we can express the limit dissipation 
rate $D_{Mode I,lim}$ as:
	\begin{equation}
		\label{eq:LimDiss}
		D_{Mode I,lim} = v_R G_{fn} b ,
	\end{equation}
and compare it with the results in Fig.~\ref{fig:Kalthoff_DissGraph_Crack}. The red curve exceeds $D_{Mode I,lim}$ at $t_{FE-B2} = 60.5 \mu$s, when the dissipated energy in ``FE-B'' reaches $93 \%$ of $G_{fn}$. If only one branch is considered (blue line), then the dissipation is below $D_{Mode I,lim}$. The branching criterion is usually based on the maximal crack tip velocity, which is not useful in our case, because the crack branching appears when the dissipation rate is not maximal (it is not at peak value). We can use $D_{Mode I,lim}$ to check for possible branching events. Fig.~\ref{fig:Kalthoff_UpperCrack_DissGraph} illustrates the dissipation of the upper crack path, which is below $D_{Mode I,lim}$, indicating that there is no crack branching.

	\begin{table}[H]
		\centering
		\begin{tabular}{p{3 cm}|p{2.5 cm}|p{2.5 cm}|p{2.5 cm}}
			\hline
			Mesh & 3280 FEs & 7320 FEs & 19800 FEs \\
			\hline
			Plane Stress & 4.89 min & 12.94 min & 33.52 min \\
			Plane Strain & 4.89 min & 11.37 min & 33.07 min \\ 
			\hline
		\end{tabular}
		\caption{Kalthoff's test: Computational time for the complete simulation 
		(with approximately 2000 increments, adjustable time step and 
		$\upDelta \tau_{max}=0.05 \mu$s).}
		\label{tab:KalthoffCompTime}
	\end{table}

\begin{figure}[H]
		\centering			
		\includegraphics[width=16cm]{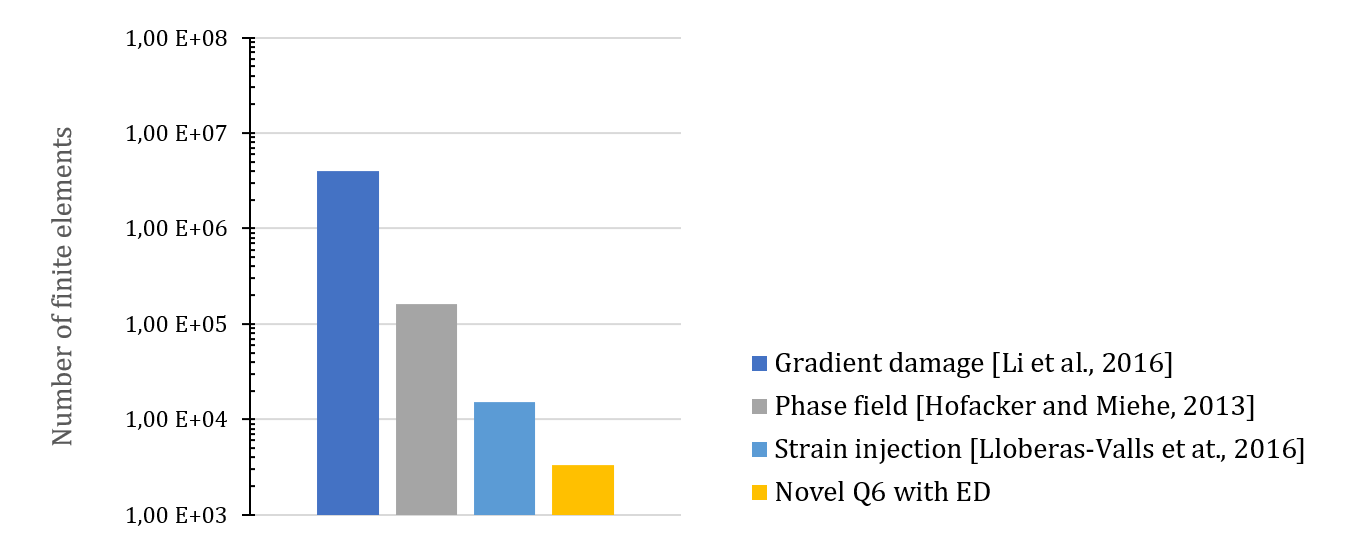}%
		\caption{Kalthoff's test. Mesh discretization used for different numerical methods: 
		Gradient method \cite{Li et al 2016-art} $\approx4.000.000$ FEs, 
		Phase field \cite{Hofacker Miehe 2013-art} $\approx160.000$ FEs, 
		Strain injection technique \cite{Lloberas-Valls et al 2016-art} $\approx15.000$ FE, 
		and novel Q6 with ED $3.280$ FE.}%
		\label{fig:Kalthoff_MeshDiscretisation}%
	\end{figure}

%	\noindent 
The computational time needed is presented in Table~\ref{tab:KalthoffCompTime} for the chosen meshes for both plane strain and plane stress models. The simulations were performed with the finite element computer code AceFEM \cite{AceFEM manual-web} on a personal computer with processor Intel Core i7-8550U CPU (1.80 GHz). As expected, both plane strain and plane stress models needed a similar amount of computational time for the same mesh. Fig.~\ref{fig:Kalthoff_MeshDiscretisation} compares mesh densities used by different numerical methods for fracture modelling in \cite{Li et al 2016-art}, \cite{Hofacker Miehe 2013-art} and \cite{Lloberas-Valls et al 2016-art}, to simulate Kalthoff's test. We can see that our approach, which relies on the Q6ED finite element with embedded discontinuity, requires the least number of elements (and thus much lower computational cost) to predict quite good approximations of the crack pattern, as demonstrated above.

\begin{figure}[H]
		\centering
		\includegraphics[width=\linewidth]{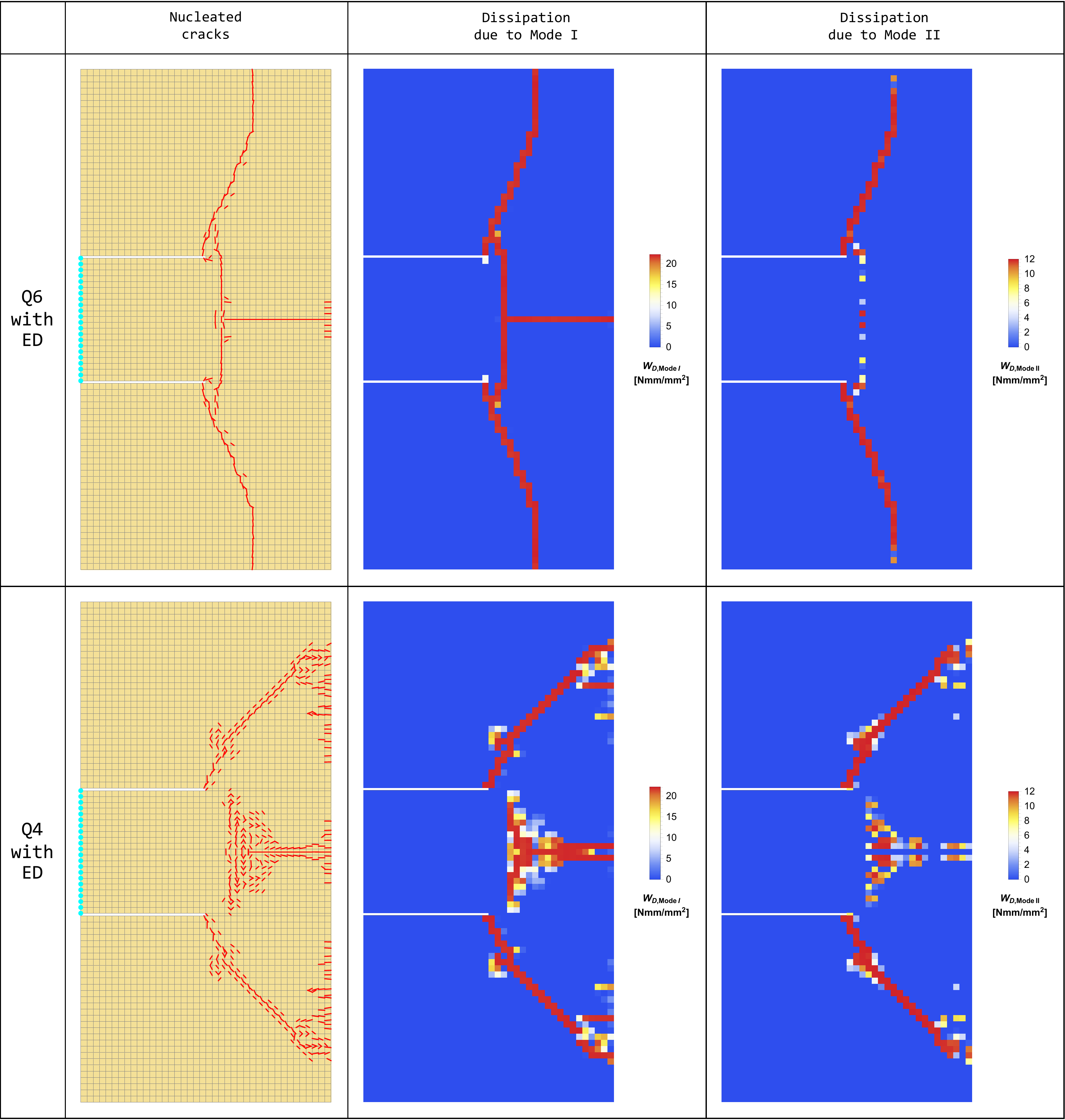}
		\caption{\correction{Kalthoff’s test: Comparing results of Q6ED (first row) with Q4ED (second row) for plain strain formulation and 3280 elements. The first column shows meshes with nucleated cracks and the last two columns display specific dissipated fracture energies due to cracking in Mode I and Mode II fashion.}} 
		\label{fig:Kalthoff_Q4_vs_Q6}
	\end{figure}

\correction{
	Kalthoff’s test is also computed by suppressing the incompatible modes (i.e. by using the Q4ED element). Fig. \ref{fig:Kalthoff_Q4_vs_Q6} shows nucleated cracks and dissipated energies for the coarse mesh (3280 elements) with plane strain elements Q6ED and Q4ED. Both formulations capture all three crack paths, but the Q6ED results are much better. Namely, for Q4ED, the upper and lower cracks propagate from the notch tips under the 60° angle towards the right edge. There are also many clustered nucleated cracks that smear dissipation of fracture energy.
	}  

\subsection{Single edge notched specimen (SENS) test}
\begin{figure}[H]
		\centering
		\includegraphics[width=0.55\linewidth]{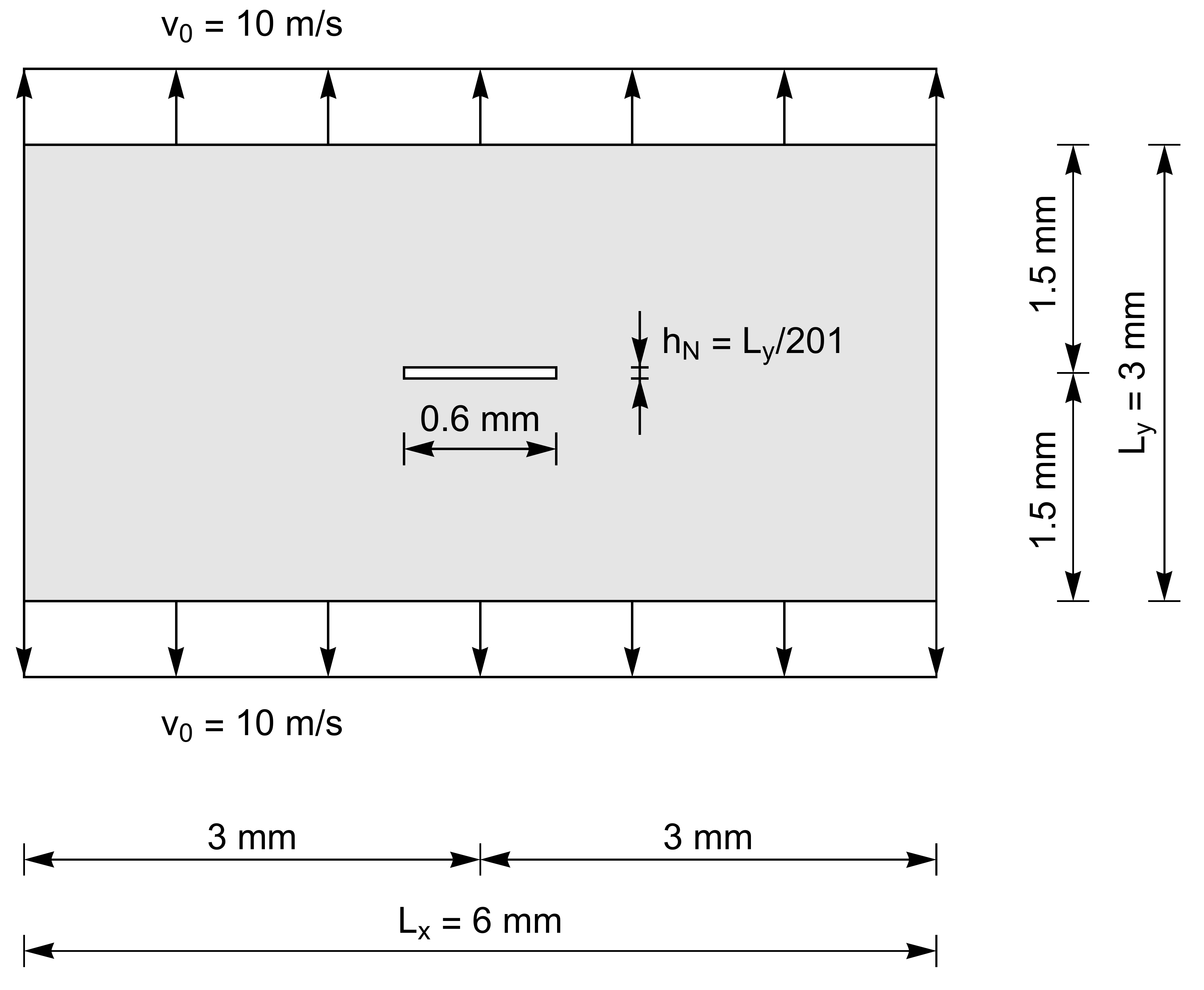}
		\caption{\correction{SENS test: Geometry and boundary conditions.}} 
		\label{fig:SENS_Geom}
	\end{figure}

\begin{figure}[H]
		\centering
		\includegraphics[width=\linewidth]{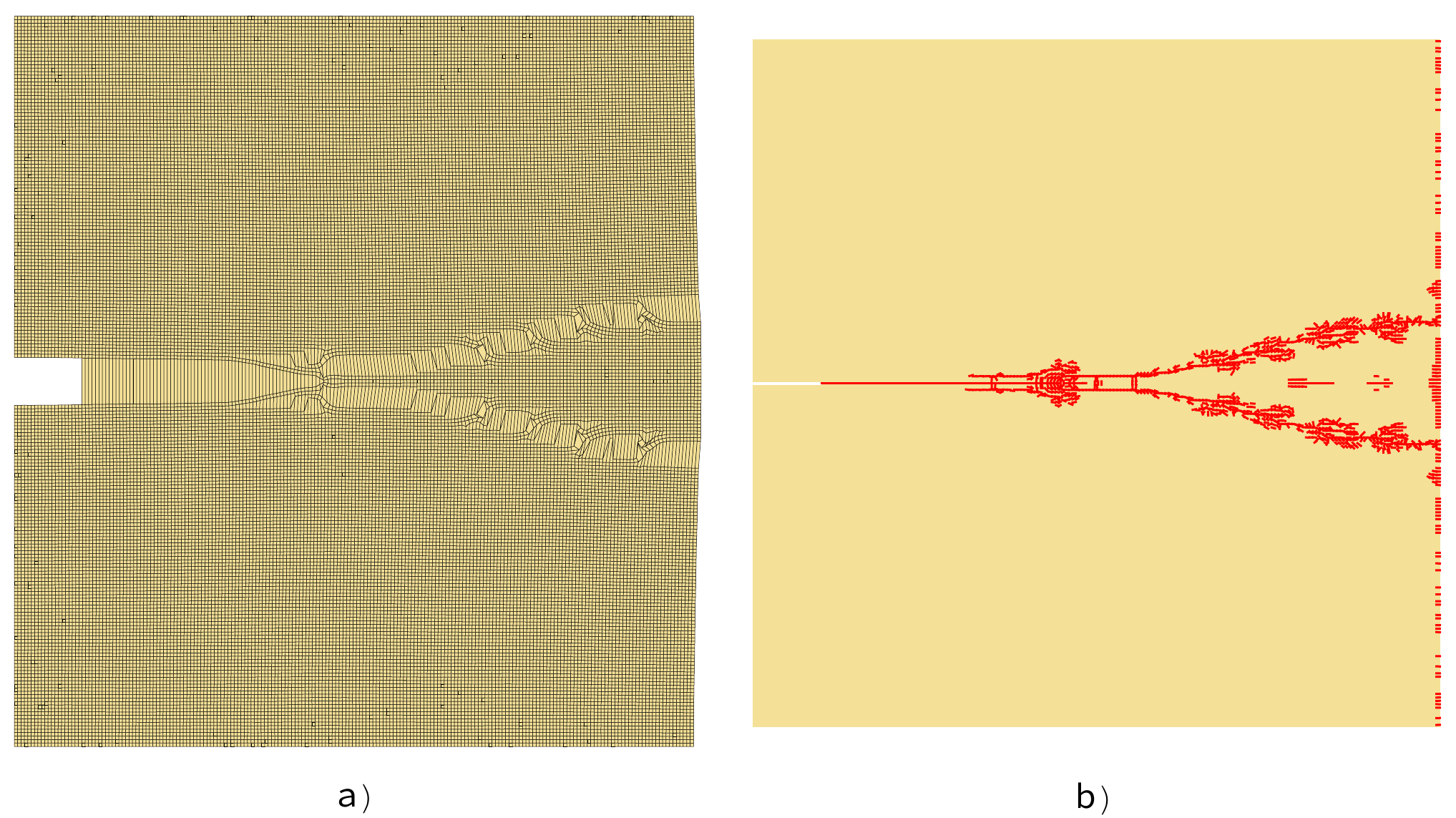}
		\caption{\correction{SENS test: a) Deformed mesh and b) nucleated cracks in the mesh at the end of simulation.}} 
		\label{fig:SENS_DefMesh_Cracks}
	\end{figure}

\begin{figure}[H]
		\centering
		\includegraphics[width=\linewidth]{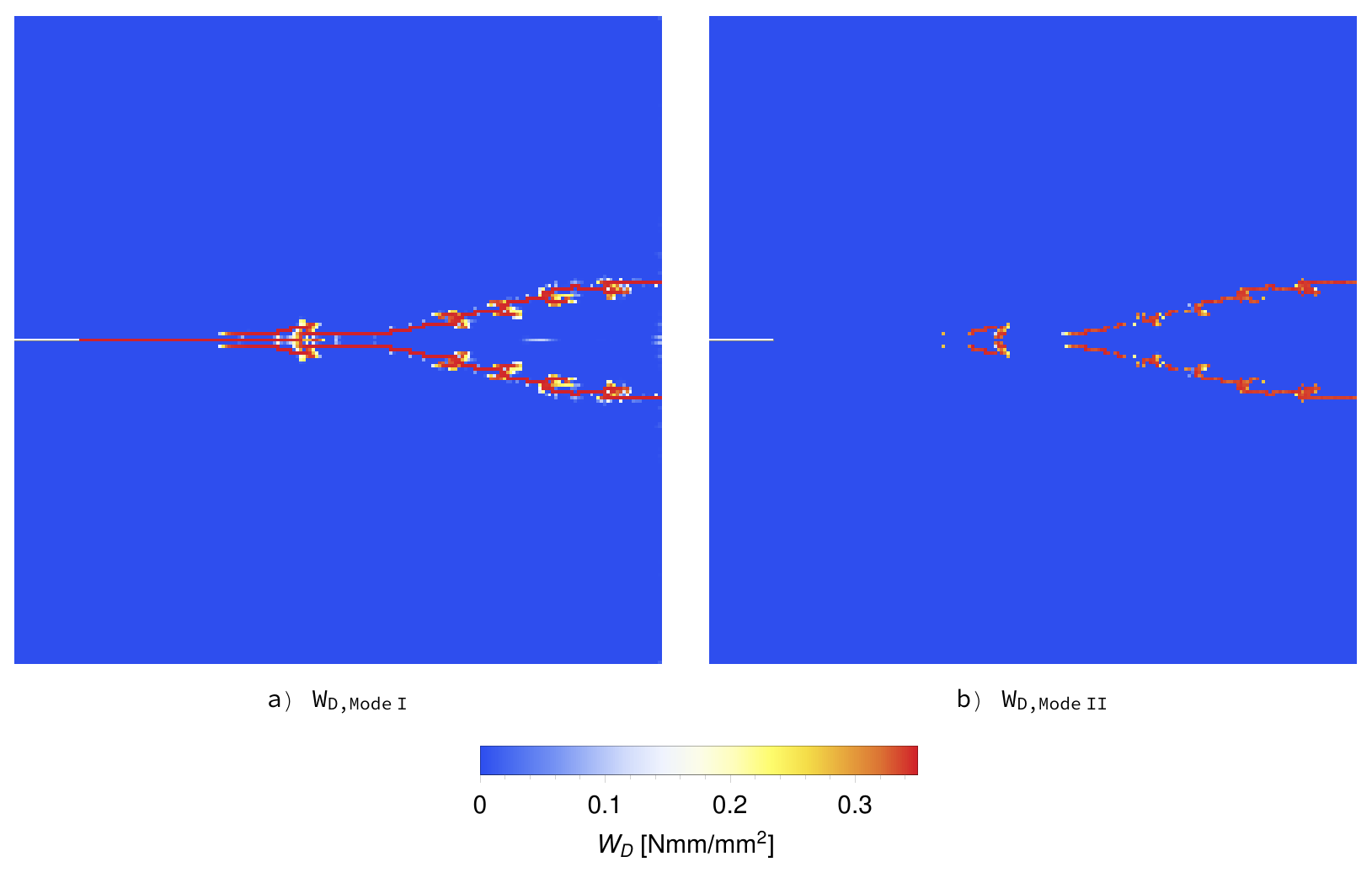}
		\caption{\correction{SENS test: Specific dissipated fracture energy for a) Mode I and b) Mode II.}} 
		\label{fig:SENS_Diss}
	\end{figure}

\begin{figure}[H]
		\centering
		\includegraphics[width=\linewidth]{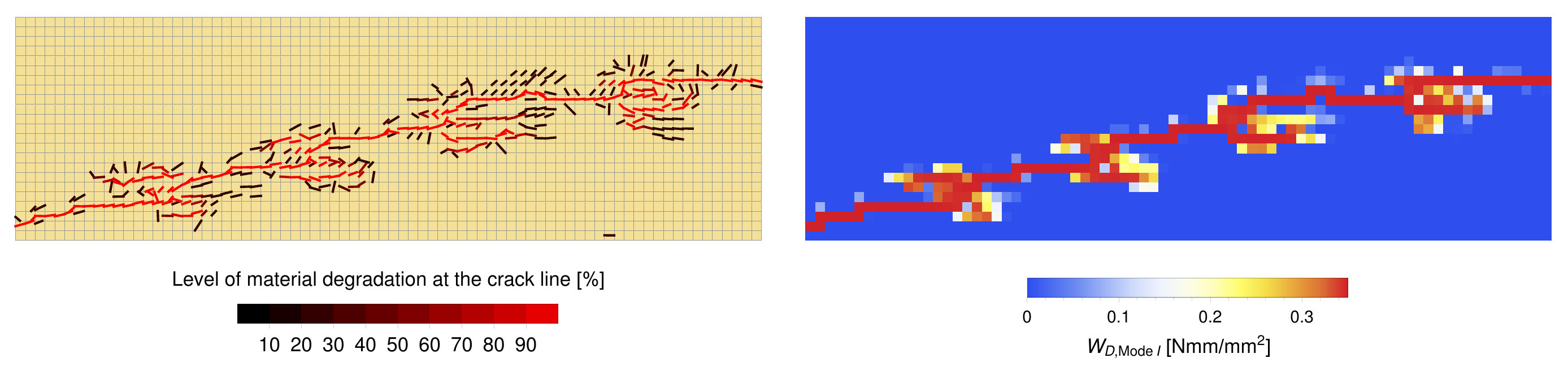}
		\caption{\correction{SENS test: Details of the mesh with the upper secondary crack path. Left: Nucleated cracks. Right: Specific dissipated fracture energy for Mode I.}} 
		\label{fig:SENS_DetailA}
	\end{figure}

\begin{figure}[H]
		\centering
		\includegraphics[width=\linewidth]{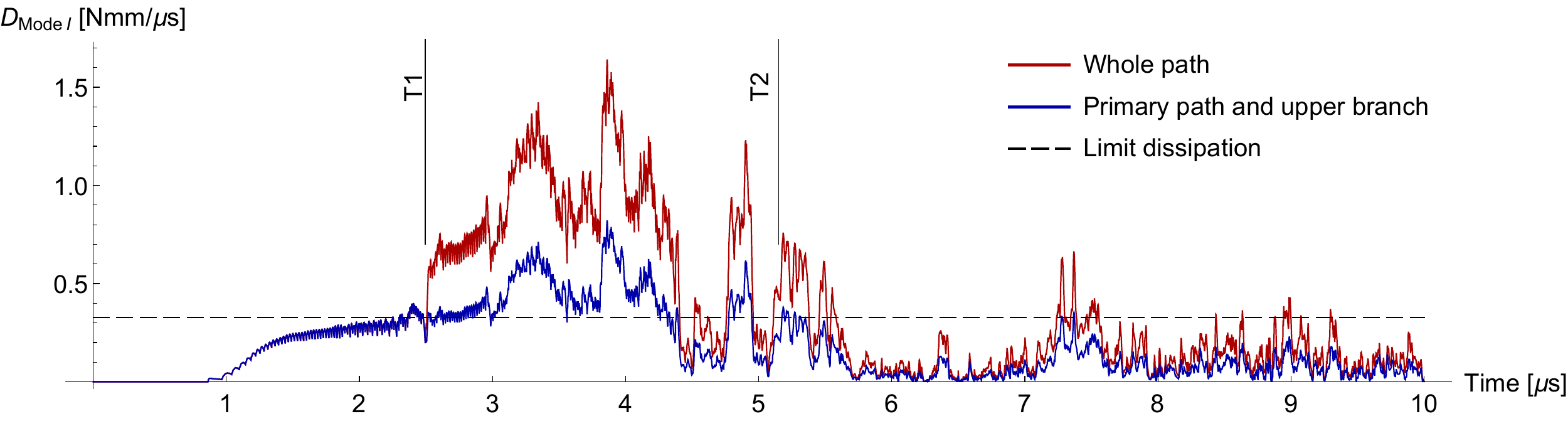}
		\caption{\correction{SENS test: Dissipation due to Mode I separation.}} 
		\label{fig:SENS_Diss_vel}
	\end{figure}

	\correction{Single edge notched specimen test from \cite{Armero Linder 2009-art} is a benchmark test for crack branching. A rectangular plate with a crack in the middle is subjected to tension. Fig. \ref{fig:SENS_Geom} shows the geometry and the boundary conditions. The velocity on top and bottom edges increases linearly from 0 at time t=0 to $v_{0}=10$ m/s at time t = $0.1 \, \upmu$s and remains constant for $10 \, \upmu$s. The material parameters are \cite{Armero Linder 2009-art}: Young's modulus $E=3240$ N/mm$^2$, Poisson's ratio $\nu=0.35$, mass density $\bar{\rho} = 1190 $ kg/m$^{3}$, Rayleigh's wave speed of $v_R = 938$ m/s, and tensile and shear strengths  $\sigma_{un} = \sigma_{um} = 129.6$ N/mm$^2$. Mode I and Mode II fracture energies are $G_{fn}= G_{fm}=0.35$ N/mm. Due to the symmetry, only one half of the plate is considered with appropriate boundary conditions. The mesh consists of 200 x 201 elements, with the existing crack modelled as a notch of thickness $h_{N} = L_{y}/201$.
	
	The results are presented in Figs.~\ref{fig:SENS_DefMesh_Cracks}--\ref{fig:SENS_DetailA}. Note that in \cite{Armero Linder 2009-art} coarser mesh with 40 x 41 elements was used. When we computed the test with such mesh, the crack line started at the notch tip and propagated in a straight line to the right edge. Our model requires a finer mesh, because the main crack branches into two secondary paths that propagate under the acute angle ($\sim 12°$) as can be seen in Figs. \ref{fig:SENS_DefMesh_Cracks} and \ref{fig:SENS_Diss}. The model response is comparable to the results obtained with the phase field method in \cite{Hofacker Miehe 2013-art}, but there the boundary conditions due to the symmetry were not taken into account.
	
	The primary crack branches into two secondary cracks: the upper and lower paths. What stands out in Fig. \ref{fig:SENS_DefMesh_Cracks}b is that each secondary path has four clusters of nucleated cracks. According to Fig. \ref{fig:SENS_Diss}, only few elements actually dissipate energy and form a crack line. Fig. \ref{fig:SENS_DetailA} displays details of the upper branch with the clustered cracks. One can see that in these clusters another branching occurs. Similar is reported in \cite{Armero Linder 2009-art}, where the crack path velocity is monitored and used in branching criterion. Additionally, Fig. \ref{fig:SENS_Diss}b shows that the fracture energy due to Mode II dissipates mostly in the areas, where crack paths kink.
	  
	Fig. \ref{fig:SENS_Diss_vel} displays the dissipation rate of fracture energy in Mode I for each time step, with two curves shown: the red curve is the dissipation in all elements, and the blue line represents the contribution of the finite elements from the primary path and the upper branch. One can see that the red and blue curves are identical until $t_{T1} = 2.5 \,\upmu$s, when the crack branching is initiated, after that the red curve is above the blue one. The secondary paths hit the edge of the plate at $t_{T2} = 5.15 \,\upmu$s, when the last element in the string of cracked elements dissipates $98\%$ of fracture energy $G_{fn}$. Although the plate is almost "broken" at $t_{T2}$, the analysis computation does not lose convergence and continues until $10 \,\upmu$s. The reason for this is in the exponential softening law, where traction stresses only approach to zero, so the finite element never separates completely and plate does not break into pieces. The computed results are compared with the limit dissipation rate based on the Rayleigh wave velocity $v_R$ and determined by Eq. \eqref{eq:LimDiss}. Fig. \ref{fig:SENS_Diss_vel} shows that the dissipation rate due to Mode I increases from zero to the limit dissipation at the instant of branching $t_{T1}$. Thereafter, the blue curve is around the limit dissipation rate. In the interval from $3 \,\upmu$s to $5 \,\upmu$s, the dissipation in the upper crack path is above the limit rate due to branching in four local regions presented in Fig. \ref{fig:SENS_DetailA}. Overall, the results indicate that the cracking mechanism in the Q6ED model can stick to the limitation from physics, as it is the Rayleigh wave velocity, without using any branching criterion.
}

\section{Conclusions}
\label{sec:Conclusions}
A novel quadrilateral finite element with embedded strong discontinuity (Q6ED) is presented and used for the simulation of fracture in quasi-brittle two-dimensional solids in both static and dynamic framework. The starting point for the embedded strong discontinuity finite element formulation is the quadrilateral with incompatible modes, sometimes denoted as Q6. It was chosen to treat the incompatible modes and the separation modes (due to the discontinuity embedding) in the same manner and to condense them on the element level. Moreover, the local equations from the incompatibility and separation modes are solved at the same time as the global equilibrium equations, which seems to be the easiest way for implementing  the embedded discontinuity. In this manner, one may stay with just one (standard) solution loop in the framework of the Newton-Raphson method and avoid a separate loop for the computation of the discontinuity parameters, which is very common (for various reasons) in the embedded discontinuity formulations. 	
%	\noindent 
Besides the above mentioned straightforward and simple implementation, a novel application of the nucleation criterion that determines when the crack is embedded into the element was proposed, and how to choose its position and orientation. It is based on the Gauss point stress state and does not use any stress averaging. 
%	\noindent 
It turns out that so designed embedded discontinuity quadrilateral Q6ED is able to compute very similar crack patterns as reported by experimental or numerical tests, without using a crack tracking algorithm. To the best of our knowledge, this is for the first time that an embedded-discontinuity formulation for a quadrilateral, which does not use any tracking algorithm, has proven successful. Indeed, the computed crack patterns are quite accurate, as has been shown by the set of computed numerical examples. Moreover, the presented formulation naturally allows for the simulation of crack branching as shown for Kalthoff's test \correction{and single edge notched specimen test}.

Furthermore, it is very fast in comparison with some other competitive approaches for fracture modelling (and also in comparison with other embedded strong discontinuity formulations for quadrilaterals, although this is not explicitly shown). It also dissipates a similar amount of the fracture energy for the computed examples as reported in the relevant references. The main drawback of the presented formulation seems to be the possibility of stress-locking (i.e.\ an incorrect transfer of stresses across the discontinuity), due to the application of only constant opening and sliding separation modes. This may occur in those elements that have a crack passing through opposite edges. This is exhibited, for example, in the four-point bending test seen in the force-displacement curve not approaching the abscissa for the fully developed crack, or in an increase of the force-displacement curve towards the end of the analysis in the double edge notched specimen test. However, our numerical experience is that for the discrete constitutive model that handles well the mixed mode opening, the exhibition of stress locking at the structural level, i.e.\ on the load-displacement curve and on the crack path, is likely to be mild for mixed mode problems. The possible inclusion of linear crack opening and sliding would have probably improved the performance of the element;	however, an efficient formulation for the linear sliding is still to be found.

	\section{Acknowledgements}
	A. Stanic and H.~G. Matthies acknowledge the support of the Deutsche Forschungsgemeinschaft (project PHAROS, DFG Geschäftszeichen GZ: \correction{MA 2236/28-1}).\\
	B. Brank acknowledges the support of the Slovenian Research Agency (project J2-1722).

\end{document}